\input amstex
\input amsppt.sty

\input epsf
\epsfverbosetrue  \magnification 1000\vsize=21
true cm \hsize=16.5 true cm \voffset=1.1 true cm \pageno=1
\NoRunningHeads \TagsOnRight

\def\p{\partial}
\def\ve{\varepsilon}
\def\f{\frac}
\def\na{\nabla}
\def\la{\lambda}
\def\al{\alpha}
\def\t{\tilde}
\def\vp{\varphi}
\def\O{\Omega}

\def\g{\gamma}
\def\G{\Gamma}
\def\si{\sigma}
\def\dl{\delta}
\def\a{O\bigl(b_0^{-\f{2}{\g-1}}\bigr)+O\bigl(b_0^{-2}\bigr)}

\def\o{\omega}

\def\ds{\displaystyle}

\topmatter
\topmatter \vskip 0.2 true cm \title{\bf Global multidimensional
shock waves for 2-D and 3-D unsteady potential flow equations}
\endtitle
\endtopmatter
\document

\vskip 0.2 true cm \footnote""{This work is supported by NSFC
(No.10931007, No.11025105, No.11001122),  the
Priority Academic Program Development of Jiangsu Higher Education
Institutions, and the DFG via the Sino-German project "Analysis of
PDEs and application".  This research was carried out
while Li Jun and Yin Huicheng were visiting the Mathematical Institute
of the University of G\"{o}ttingen.}\footnote""{Ingo Witt was
partly supported by the DFG via the Sino-German project ``Analysis of
PDEs and application.'' } \vskip 0.1 true cm \centerline{Li Jun \quad\quad Ingo Witt
\qquad
Yin Huicheng} \vskip 0.3 true cm

\centerline {\bf Abstract} \vskip 0.1 true cm Although local
existence of multidimensional shock waves has been established
in some fundamental references, there are few results on the global
existence of those waves except the ones for the unsteady potential
flow equations in $n-$ dimensional spaces ($n\ge 5$) or  in special
unbounded space-time domains with some artificial boundary
conditions. In this paper, we are concerned with both the local and global
multidimensional conic shock wave problems for unsteady potential flow
equations when a pointed piston (i.e., the piston at the initial
time degenerates into a single point) or an explosive wave expands fast in 2-D or 3-D
static polytropic gases. It is shown that a multidimensional shock wave
solution to such a class of quasilinear hyperbolic problems not only exists locally but also
exists globally in the whole
time-space and tends to a self-similar solution as $t\to\infty$.

\vskip 0.1 true cm

pseudo-differential operator,  improved Hardy-type inequality, modified
Klainerman's vector fields\vskip 0.1
true cm

\vskip 0.2 true cm

\centerline{\bf Contents}

\vskip 0.1 true cm

{$\S 1.$ Introduction}

{$\S 2.$ Analysis on the self-similar background solution}

{$\S 3$. Local-in-time existence}

\quad {$\S 3.1.$ A reformulation of problem (1.6) with (1.7)-(1.9) under some
nonlinear transformations}

\quad {$\S 3.2.$ Construction of an approximate solution to problem (3.10)
with (3.11)-(3.12)}

\quad {$\S 3.3.$ Choice of iteration scheme and proof of local existence}

\qquad {$\S 3.3.1.$ Choice of iteration scheme}

\qquad {$\S 3.3.2.$ Solvability and energy estimates of problem (3.48)}

\qquad {$\S 3.3.3.$ Solvability of problem (3.10)-(3.12) and proof of Theorem 3.1}

\qquad {$\S 3.3.4.$ Proof of an elementary Proposition}

{$\S 4$. Another reformulation of (1.6) with (1.7)-(1.9) and some preliminaries}

{$\S 5$. Proof of Theorem 1.1 in the case of $n=3$}

\quad {$\S 5.1.$ First order weighted energy estimate}

\quad {$\S 5.2.$ Higher order weighted energy estimates}

\quad {$\S 5.3.$ Proof of Theorem 1.1 for $n=3$}

{$\S 6$. Sketch on the proof of Theorem 1.1 for $n=2$}

Appendix A. Some basic computations

Appendix B. Modified background solution

References

\newpage

\centerline{\bf $\S 1.$ Introduction} \vskip 0.3
true cm

Although local existence results of multidimensional planar shock
waves have been established early in some fundamental references
(see [23-25] and the references therein), there are few results on
the global existence of those waves except the ones for the unsteady
potential flow equations in much higher space dimensions ($n\ge 5$)
or in special unbounded space-time domains with some artificial
boundary conditions (i.e., the Dirichlet boundary conditions of
potentials on some specially chosen boundaries, one can see [10-11] and
[36]). This paper concerns both the local and global
multidimensional shock wave problems when a pointed piston (i.e.,
the piston at the initial time degenerates into a single point) or an explosive wave
expands in the 2-D or 3-D isentropic irrotational gas. The piston
problem is a fundamental one in gas dynamics, in particular, the
expansive pointed  piston  problem is also closely related or similar to the
study of explosive waves in physics, one can see the references [1],
[4], [7-8], [18], [22], [26], [30-32]  and so on. Such a problem is
also one of the basic models in establishing the theory of weak
solutions to the multidimensional quasilinear hyperbolic equations
or systems. As described in pages 120 of [7]: A basic and typical
motion of gas is the one caused by a piston in a tube starting from
rest and suddenly moving with constant velocity $u_P$ into the quiet
gas. No matter how small $u_P$ is, the resulting motion cannot be
continuous. Generally speaking, if the piston recedes, then a
rarefaction wave will be caused, and otherwise, if the piston is
pushed, then a shock wave will be formed. In this paper, we will
study the multi-dimensional case of a piston problem as in [1], [4]
and so on, but there authors posed some symmetric properties on
solutions. That is, we suppose that there is a rest gas filling the
whole space outside a given pointed piston with expansive boundary.
With the development of time, the pointed piston gradually expands
its boundary into the air in two-dimensional spaces or
three-dimensional spaces. Subsequently, there will be a
multidimensional shock wave moving into the air away from the piston
(See the Figures 1-3 below). Mathematically, this is an initial
boundary value problem for  the 2-D or 3-D compressible Euler
system, which contains a free boundary (shock surface) and a moving
boundary (surface of expansive piston). For the rapidly expansive
piston in the air, we will establish both the local and global
existence of a multidimensional (2-D or 3-D) shock wave solution.
\vskip 0.3 true cm
$$
\epsfysize=80mm \epsfbox{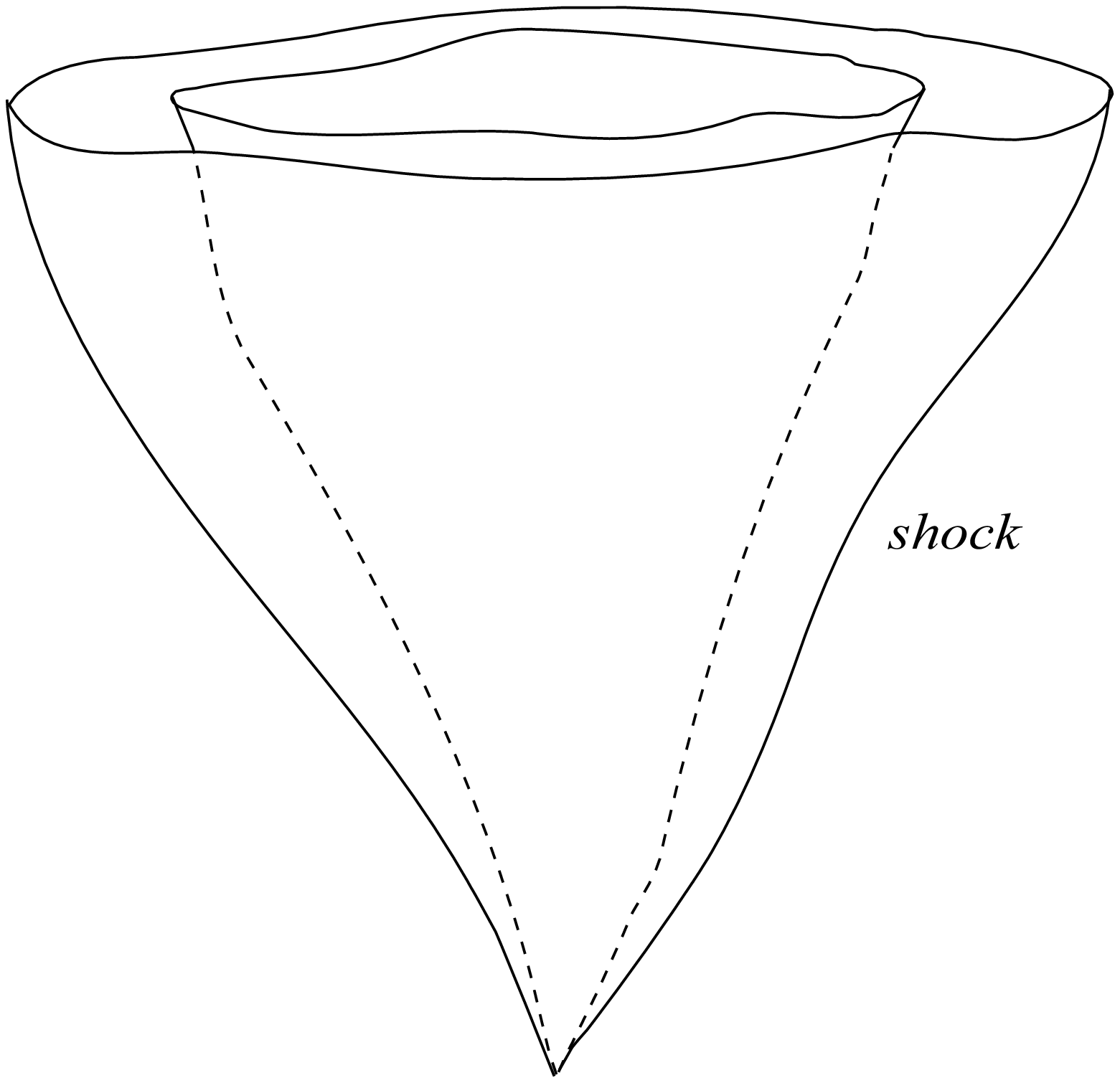}
$$
\centerline{\bf Figure 1. A conic shock is formed when a
pointed piston  expands in polytropic gas}

\vskip 0.5 true cm

$$
\epsfysize=70mm\epsfbox{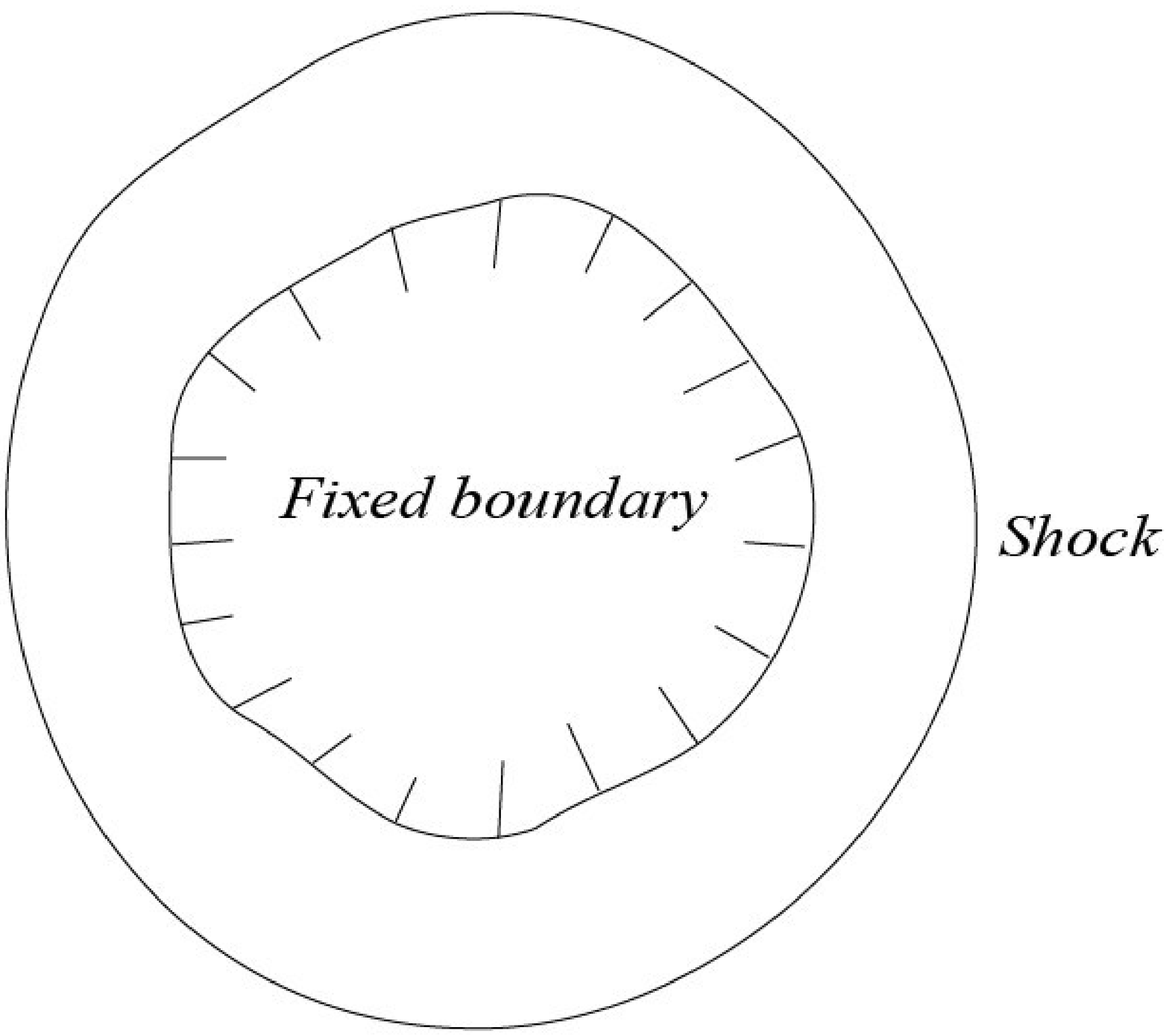}
$$
\centerline{\bf Figure 2. The picture of a 2-D conic shock at $t=1$}
\vskip 1.8 true cm
$$
\epsfysize=80mm\epsfbox{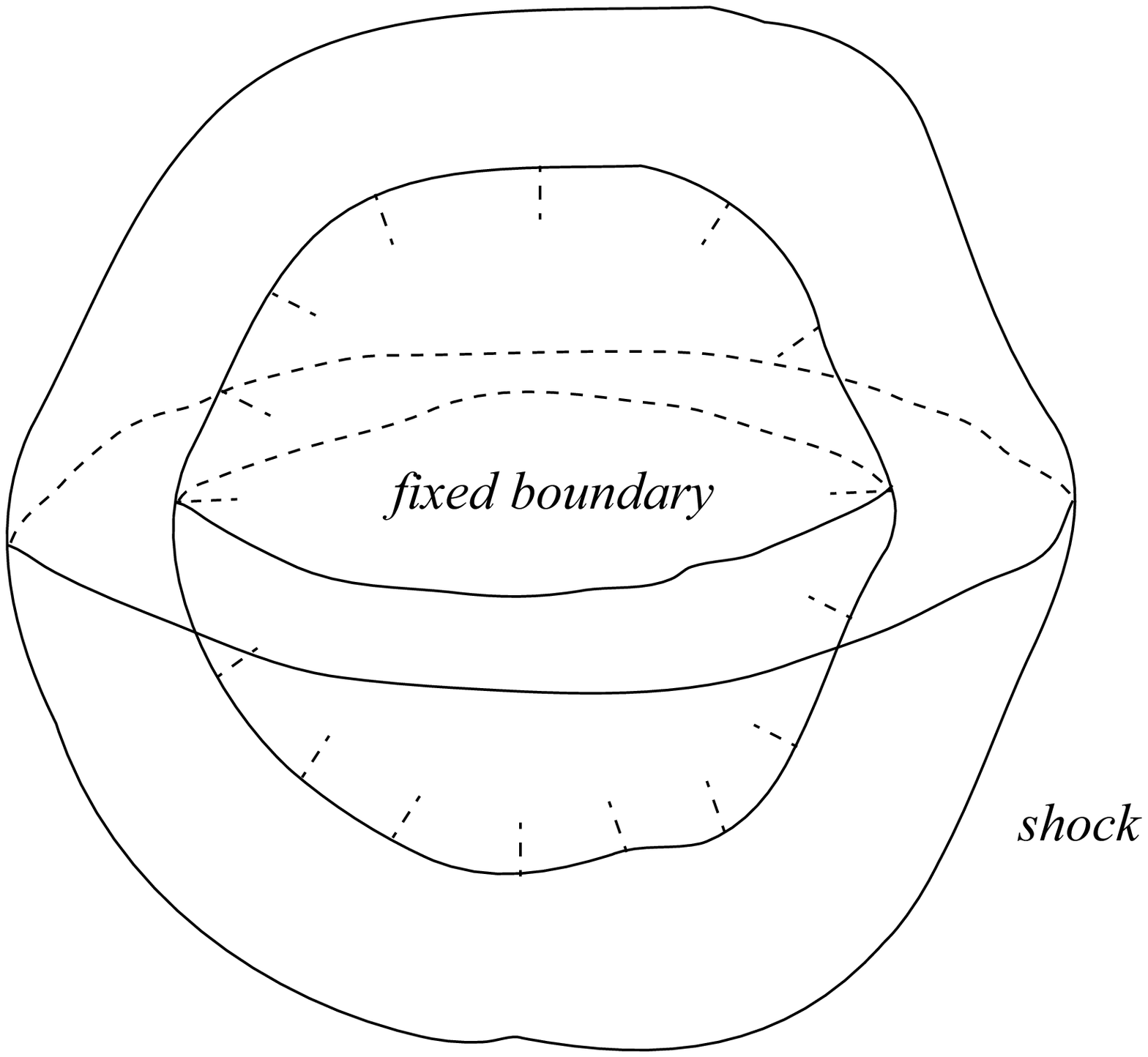}
$$
\centerline{\bf Figure 3. The picture of a 3-D conic shock at $t=1$}

\vskip 0.6 true cm

We will use $n$-dimensional ($n=2,3$) unsteady potential equation to
describe the motion of the polytropic gas in expansive pointed
piston problem (this model is also recommended in [4-5], [24],
[35], [38] and so on), where polytropic gas means that the pressure
$P$ and the density $\rho$ of the gas are described by the state
equation $P=A\rho^{\g}$ with $A>0$ a constant and the adiabatic
constant $\g$ satisfying $1<\g<3$ (for the air, $\g\simeq 1.4$). Let
$\Phi(t,x)$ be the potential of velocity $u=(u_1,\cdots, u_n)$ with
$x=(x_1, \cdots, x_n)$, i.e., $u_i=\p_i\Phi$ $(1\leq i\leq n)$, then
it follows from the Bernoulli's law that
$$\p_t\Phi+\f{1}{2}|\na_x\Phi|^2+h(\rho)=B_0,\tag1.1$$
here $h(\rho)=\ds\f{c^2(\rho)}{\g-1}$ is the specific enthalpy,
$c(\rho)=\sqrt{P'(\rho)}$ is the local sound speed,
$\na_x=(\p_1,\cdots, \p_n)$, $B_0=\ds\f{c^2(\rho_0)}{\g-1}$ is the
Bernoulli's constant of static gas with the constant density
$\rho_0$.

By (1.1) and the implicit function theorem due to
$h'(\rho)=\ds\f{c^2(\rho)}{\rho}>0$ for $\rho>0$, then the density
function $\rho(t, x)$ can be expressed as
$$\rho=h^{-1}\biggl(B_0-\p_t\Phi-\f{1}{2}|\na_x\Phi|^2\biggr)
\equiv H(\na\Phi),\tag1.2$$ here $\na=(\p_t, \na_x).$

Substituting (1.2) into the mass conservation equation
$\p_t\rho+\dsize\sum_{i=1}^{n}\p_i(\rho u_i)=0$ of gas yields
$$\p_t(H(\na\Phi))+\dsize\sum_{i=1}^n\p_i\bigl(H(\na\Phi)\p_i\Phi\bigr)=0.\tag1.3$$

More intuitively, for any $C^2$ solution $\Phi$, (1.3) can be
rewritten as the following second order quasilinear equation
$$
\p_t^2\Phi+2\dsize\sum_{k=1}^n\p_k\Phi\p_{tk}^2\Phi+\ds\sum_{i,j=1}^{n}\p_{i}\Phi\p_j\Phi\p_{ij}^2\Phi
-c^2(\rho)\Delta\Phi=0,\tag 1.4
$$
here $c(\rho)=c(H(\na\Phi))$, and the Laplace operator
$\Delta=\ds\sum_{i=1}^{n}\p_i^2$.

It is noted that (1.4) is strictly hyperbolic with respect to the
time $t$ when $\rho>0$ holds.

\medskip

For convenience to write and compute later on, the following
spherical coordinates are often used
$$(t, r, \o)=(t, |x|, \f{x}{|x|}).\tag 1.5$$

Under the coordinate transformation (1.5), we suppose that the
expansive path of the pointed piston in static gas is $\Sigma$:
$r=\sigma(t,\o)$ and denote that the potential before and
behind the resulting shock front $\Gamma$: $r=\zeta(t,\o)$ with $\zeta(0,\o)=0$ are written by
$\Phi^-(t,x)$ and $\Phi^+(t,x)$ respectively. And the corresponding
domains are denoted by $\O_-$ and $\O_+$ respectively. Since the gas
ahead of the piston is static, then $\Phi^-(t,x)\equiv 0$ can be chosen
in $\O_{-}$. In $\O_{+}$, $\Phi^+$ satisfies
$$\align
&\p_t^2\Phi^++2\ds\sum_{k=1}^{n}\p_k\Phi^+\p_{tk}^2\Phi^++\ds\sum_{i,j=1}^{n}\p_i\Phi^+
\p_j\Phi^+ \p_{ij}\Phi^+-c^2(\rho^+)\Delta\Phi^+=0,\tag
1.6\endalign$$ here $c(\rho^+)=c(H(\na\Phi^+))$.

On the surface $\Sigma$ of expansive piston, $\Phi^+$
satisfies the following solid  boundary condition
$$
\Cal B_{\si}\Phi^+\equiv \p_r\Phi^+-\f{1}{r^2}\ds\sum_{i=1}^{2n-3}Z_i\si\cdot
Z_i\Phi^+=\p_t\sigma,\tag 1.7
$$
where
$$Z_1=x_1\p_2-x_2\p_1,\ Z_2=x_2\p_3-x_3\p_2,
\ Z_3=x_3\p_1-x_1\p_3,$$
which form a basis of smooth vector fields tangent to the sphere $\Bbb S^2$. In particular,
$Z_1$ is a smooth vector field tangent to the unit circle $\Bbb S^1=\{(x_1, x_2): x_1^2+x_2^2=1\}$.

Meanwhile, on the shock surface $\Gamma$, by the
equation (1.3), the corresponding
Rankine-Hugoniot condition is
$$
H(\na\Phi^+)\p_r\Phi^+-(H(\na\Phi^+)-\rho_0)\p_t\zeta=\f{1}{r^2}H(\na\Phi^+)
\ds\sum_{i=1}^{2n-3}Z_i\zeta\cdot
Z_i\Phi^+\quad\text{on}\quad \Gamma.\tag 1.8$$

Moreover, the potential function $\Phi(t,x)$ is continuous across
$\Gamma$, namely,
$$\Phi^+=\Phi^-\equiv 0\quad\text{on}\quad \Gamma.\tag1.9$$

Additionally, the physical entropy condition holds on the shock surface
$$H(\na\Phi^+)>\rho_0\quad\text{on}\quad \Gamma.\tag 1.10$$

The main result in our paper can be stated as:

{\bf Theorem 1.1.} {\it For $n=2,3$, if the equation of $\Sigma$ is
$r=\si(t,\o)=tb(t,\o)$, here $b(t,\o)\in
C^{\infty}([0,+\infty)\times \Bbb S^{n-1})$,
$|Z^k(b(0,\o)-b_0)|\le\ds C_k\ve$ and
$\ds|\p_t^{k_0}Z^{k}(b(t,\o)-b_0)|\leq \f{C_{k_0
k}\ve}{(1+t)^{k_0}}$ with $k_0, |k|=k_1+k_2+k_3\in\Bbb N\cup\{0\}$, $Z\in\{Z_1, ..., Z_{2n-3}\}$ and a small
positive constant $\ve$, then for suitably large constant $b_0$,
there exists a positive constant $\ve_0$ depending on $b_0$ and
$\g$, such that the problem $(1.6)$ together with $(1.7)$-$(1.10)$
has a global shock solution $(\Phi^+(t,x),\zeta(t,\o))$ as
$0<\ve<\ve_0$. Moreover, $(\na\Phi^+,\ds\f{\zeta(t,\o)}{|x|})$ tend
to the corresponding ones, which are formed by the symmetrically
expansive pointed piston $r=b_0 t$ in static gas, with the decay
rate $(1+t)^{-m_0}$ for any positive numbers
$m_0<\ds\f{5}{4}-\ds\f{1}{4}\sqrt{\ds\f{\g+1}{2}}$ if $n=2$ and
$m_0<\ds\f{3}{2}-\ds\f{1}{4}\sqrt{\ds\f{\g+7}{2}}$  if $n=3$
respectively.}\medskip

{\bf Remark 1.1.} {\it From the expression $r=\si(t,\o)=tb(t,\o)$ of
the equation of $\Sigma$ and the assumptions on $b(t,\o)$ in
\rom{Theorem 1.1}, we know that the normal expansive velocity of the
pointed piston is $\p_t\si(t,\o)=b_0+O(\ve)$, which is a small
perturbation of the constant expansive speed $b_0$.}\medskip

{\bf Remark 1.2.} {\it It is noted that the nonlinear hyperbolic
equation $(1.6)$ is actually a second order quasi-linear wave
equation in three or two space dimensions. By a direct verification,
one can  know that $(1.6)$ does not fulfill the ``null-condition''
or admits the cubic order nonlinearity under some structural
transformations (one can see $[6], [16]$ or $[33]$-$[34]$).
Therefore, in terms of the results in $[2]$-$[3]$, $[13]$, $[15]$,
$[29]$, $[37]$ and so on, if there is no main shock for the equation
$(1.6)$, then the classical solution of $(1.6)$ must blow up in
finite time. Our result in Theorem $1.1$ means that the main
multidimensional shock can absorb all possible compressions of the
flow and prevent the formations of new shocks and other
singularities in expansive piston movement as in [19] or [36].}\medskip

{\bf Remark 1.3.} {\it Since BV spaces fail for the multidimensional
hyperbolic equations or systems as shown in $[27]$, then the Glimm
scheme method $($see $[22]$ and the references therein$)$ can not be
used to treat our really multidimensional problem in this paper. On
the other hand, there are no any symmetric assumptions on the
movement of the piston in our problem, then the equation $(1.4)$ can
not be reduced into a second order nonlinear elliptic equation with
a free boundary in a bounded domain as in $[4]$, where the solvability
is shown by some techniques from linear elliptic equations.}\medskip

{\bf Remark 1.4.}
 {\it For the solutions of $n$-dimensional linear wave equations,
it is well-known that their optimal decay rates are $\ds
(1+t)^{\f{1-n}{2}}$ as $t\to\infty$ \rom{(see [13])}. Compared with
this, since one can obtain the decay rates $\ds\f{1}{(1+t)^{m_0}}$
with $\ds\f12<m_0<1$ in 2-D case and $\ds\f34<m_0<1$ in 3-D case
respectively in Theorem $1.1$, we obtain a better decay in 2-D case
but a lower decay in 3-D case with respect to the problem
\rom{(1.6)} together with \rom{(1.7)-(1.10)}. However, due to some
special properties of the shock boundary conditions
\rom{(1.8)-(1.9)}, such decay rates will be enough to establish the
global existence in our nonlinear problem.}\medskip

{\bf Remark 1.5.} {\it Since the background solution is self-similar
in our problem, then the related nonlinear problem admits largely
variable coefficients as in \rom{[21]} and the related analysis are rather
involved.}
\medskip

{\bf Remark 1.6.} {\it \rom{(1.10)} naturally holds if we have
shown the solution of problem $(1.6)$ together with $(1.7)$-$(1.9)$
is a small perturbation of the background solution since the
background solution satisfies \rom{(1.10)} by the physical entropy
condition. Therefore, we can only focus on problem $(1.6)$ together
with $(1.7)$-$(1.9)$ from now on.}

\medskip

Let us comment on the proof of Theorem 1.1. As the first step, we
will establish the local existence of solution to problem (1.6)
together with (1.7)-(1.9) in $\O_+$. For this aim, we have to
overcome the difficulties induced by the unknown shock $\G$ and the
double conic point $(t,x)=(0, 0)$ between two conic surfaces
$\Sigma$ and $\Gamma$. We will apply a partial hodograph
transformation as in [25] and [4] to change the unknown domain
$\O_+$ into a fixed cylindrical domain $\t\O_+=(0, \infty) \times
(1, 2)\times\Bbb S^{n-1}$, while the equation (1.6) and boundary
conditions (1.7)-(1.9) become a rather complicated second order
nonlinear equation including the solution $\psi$ itself as well as
the derivatives of $\psi$ and three involved nonlinear boundary
conditions on the new resulting boundaries respectively. Based on this,
we start to construct a suitable approximate solution
$\psi_a^{\kappa_0}$ to this resulting nonlinear hyperbolic equation
by solving a series of second order elliptic equations, and
subsequently consider the related nonlinear equation together with
the corresponding boundary conditions on the function
$\psi-\psi_a^{\kappa_0}$, so that the weighted energy estimates are
derived and further the local existence of shock solution can be
established. In this process, except the construction of approximate
solution, we have to treat the well-posedness problems on the linear
hyperbolic equation with an inhomogeneous Neumann-type boundary
condition and an oblique derivative boundary condition, and then
continue to treat the nonlinear hyperbolic equation together with a
nonlinear Neumann-type boundary condition and a nonlinear uniform
oblique derivative boundary condition. Here we emphasize that some
fundamental methods in [5]  are not available for our problem. The
main reasons are: If we directly use a finite power expansion of $t$
to look for an approximate solution of (1.6) together with
(1.7)-(1.9) as in [5] (see (2.11), (2.1.6) and (2.1.10) in Section 2
of [5]), then  a series of equations only in the same domain (see
(2.1.4)-(2.1.5) in [5]) are obtained and thus a real approximate
solution can not be found since the crucial property of the free
boundary $\G$ is neglected. On the other hand, since the
Neumann-type boundary condition on the moving surface $\Sigma$ does
not satisfy the uniform Lopatinski condition, namely, the Local
Stability  Condition on $\Sigma$ introduced in [18] is not fulfilled
(however, the lines 2-3 from below on pages 177 of [5] give a
different assertion), then the well-posedness problem and further the
energy estimate on the resulting linearization problem can not be
established directly by the results in [18] other than claimed in
the proof procedure of Theorem 3 of [5] since we have to treat the
well-posedness problem on the second order hyperbolic equations with
an inhomogeneous Neumann-type boundary condition and a uniform
oblique derivative boundary condition in order to solve the
linearized problems from the iteration scheme. Thanks to our delicate analysis and
some ideas in [14] and [28], we can finally complete this task.

Based on the local existence result in the above, we will utilize
the continuous induction method to prove the global existence as in
[20] and so on. To achieve this,  we need to derive global
weighted energy estimates for the problem (1.6) with (1.7)-(1.9).
By such estimates, one then obtains the global existence,
stability, and the asymptotic behavior of the shock solution to the
perturbed nonlinear problem (1.6). The key methods in the
analysis to obtain weighted energy estimates are to chose the
appropriate multipliers and establish a new Hardy-type inequality on
the shock surface in terms of the special structures of the shock
boundary conditions. Finding such suitable multipliers is much more
delicate due to the following reasons: First, in order to obtain the
global existence, one needs to establish a global estimate
independent of the time $t$, of the potential function and its
derivatives on the boundaries as well as in the interior of $\O_+$.
This yields strict restrictions on the multiplier and makes the
computations involved. Second, for the three dimensional case, the
Neumann-type boundary condition (1.7) fulfilled by $\Phi^+$ yields
additional difficulties compared to [11] and [19], where [11] treats the case
of an artificial Dirichlet-type boundary condition for the potential
on a multi-dimensionally perturbed conic surface but away from the
conic point. The latter plays a key role in the analysis of [11]
since the corresponding Poincar\'e inequalities are available on the
shock surface and in the interior of the downstream domain
respectively, while this is not the case in the problem treated
here. Furthermore, it should be noted that the arbitrary and
artificial closeness between the shock surface and the fixed
boundary boundary also plays a crucial role in the analysis of [11],
which is also not the case for our problem. Meanwhile, the authors in
[19] only treat the 2-D shock problem of steady potential flow equation
if one takes the supersonic direction as the time direction, and where only the fixed
circular cone boundary is considered. Thanks to some careful
analysis together with an improved Hardy-type inequality derived by
utilizing the special structures of the shock boundary conditions
(see (4.11)-(4.12) below), we finally overcome all these
difficulties and obtain a uniform estimate of $\|\na_{t,
x}\Phi^+\|_{L^2(\O_+)}$. From this, higher-order energy estimates of
$\na_{t, x}\Phi^+$ can be established by making full use of modified
Klainerman's vector field and commutator arguments together with a
careful verification that some suitably higher-order derivative
combinations of the solution satisfy the Neumann-type boundary condition on
$\Sigma$. This finally derives Theorem~1.1.

The paper is organized as follows: In \S2, we derive some basic
estimates on the self-similar background solution which is formed by
the symmetric pointed expansion piston in the static polytropic gas.
This will be required to simplify the nonlinear problem (1.6)
together with (1.7)-(1.9) and look for the multipliers in related
energy estimates later on. The local-in-time existence of the
problem (1.6) with (1.7)-(1.9) is established in $\S 3$. In $\S4$,
we reformulate problem (1.6)-(1.9) by decomposing its solution as a
sum of the modified background solution and a small perturbation
$\dot\vp$ so that its main part can be studied in a convenient way.
In $\S 5$, we first establish a uniform weighted energy estimate for
the corresponding 3-D problem, where an appropriate multiplier is
also constructed. Based on such an energy estimate, we obtain a
uniform weighted energy estimate of $\na_{t,x}\dot\vp$ for the
nonlinear problem through establishing an improved Hardy-type
inequality. By the estimates derived in the first step, we continue
to establish uniform higher-order weighted energy estimates of
$\na_{t,x}\dot\vp$ in the case of $n=3$ through looking for the
suitably modified Klainerman's vector fields. Finally, the proof of
Theorem~1.1 for the case $n=3$ is completed by using Sobolev's
embedding theorem and continuous induction. In
$\S 6$, we give the sketch of the proof on Theorem 1.1 in the case
of $n=2$. Some basic computations are arranged in Appendix A.
Additionally, in order to deal with the Neumann-type boundary
condition on the curved boundary $r=\si(t,\o)$, we need to modify
the self-similar solution and obtain a modified background solution
$\Phi_a$ so that the boundary condition (1.7) can be fulfilled. This
will be given in Appendix B.

In what follows, we will use the following convention:

$\centerdot$ $C$ stands for a generic positive constant which does
not depend on any quantity except the adiabatic constant $\gamma$
($1<\g<3$).

$\centerdot$ $C(\cdot)$ represents a generic positive constant which
depends on its argument(s).

$\centerdot$ $O(\cdot)$ means that $|O(\cdot)|\leq C|\cdot|$ holds.


$\centerdot$ $dS$ stands for the surface measure in the
corresponding surface integral.

\vskip 0.5 true cm \centerline{\bf $\S 2.$ The analysis on the
self-similar background solution} \vskip 0.4 true cm

In this section, we will give some detailed properties on the
background self-similar solution which is formed by the symmetric
pointed expansive piston $r=b_0t$ in the static polytropic gas with suitably large
$b_0$. These properties will be applied in the later
analysis of $\S 3-\S 6$ below.

For the pointed expansive piston $r=b_0t$, there will
appear a conic shock $r=s_0 t$ $(s_0>b_0)$ in the static gas. Moreover the solution of
(1.3) with (1.1) behind the shock surface is self-similar, that is,
the states of density and velocity between the shock front and the
surface of piston have such forms: $\rho=\rho(s), u_i=
u(s)\ds\frac{x_i}{r} (i=1,\cdots, n)$ with $s=\ds\frac{r}{t}$. In
this case, the $n-$dimensional potential equation can be reduced to
a nonlinear ordinary differential system as follows
$$
\cases
\ds\rho'(s)=\frac{(n-1)(s- u)\rho u}{s\bigl((s-u)^2-c^2( \rho)\bigr)},\\
\ds u'(s)=\frac{(n-1)c^2(\rho)u}{s\bigl((s-u)^2-c^2(\rho)\bigr)}
\endcases\quad
\text{for}\quad b_0\leq s\leq s_0.\tag2.1
$$

By [1] or [4], the denominator
$(s-u)^2-c^2(\rho)<0$ holds for $b_0\leq s\leq s_0$. This means that
the system (2.1) makes sense.

On the shock front $r=s_0 t$, it follows from the Rankine-Hugoniot
conditions and Lax's geometric entropy conditions on the 2-shock
that
$$
\cases s_0[\rho]-[\rho u]=0,\\
s_0[u]-[\f12 u^2+h(\rho)]=0
\endcases\tag2.2
$$
and
$$
\cases
 u(s_0)-c(\rho(s_0))<s_0<u(s_0)+c(\rho(s_0)),\\
c(\rho_0)<s_0.
\endcases\tag2.3
$$

Additionally, the flow satisfies the fixed boundary condition on
$s=b_0$
$$u(b_0)=b_0.\tag2.4$$

It has been shown that the boundary value problem (2.1)-(2.4) can be
solved by [1] or [4].

\medskip

For suitably large $b_0$, some  properties on the background
solution can be given as follows:\medskip

{\bf Lemma 2.1.}  {\it For suitably large $b_0$ and $1<\g<3$, one has for
$b_0\le s\le s_0$,\medskip

$(i)$ $s_0=b_0\biggl(1+\a\biggr)$.

$(ii)$ $u(s)=b_0\biggl(1+\a\biggr)$.

$(iii)$
$\dsize\rho(s)=\bigl(\f{\g-1}{2A\g}\bigr)^{\f{1}{\g-1}}b_0^{\f{2}{\g-1}}
\biggl(1+\a\biggr)$.

$(iv)$ $u^2(s)-c^2(\rho(s))=\ds\f{3-\g}{2}b_0^2\biggl(1+\a\biggr)>0$.

$(v)$ $(s-u(s))^2-c^2(\rho(s))=-\ds\f{\g-1}{2}b_0^2\biggl(1+\a\biggr)<0$.

$(vi)$  $u(s)+c(\rho(s))-s=\sqrt{\ds\f{\g-1}{2}}b_0\biggl(1+\a\biggr)>0$.

$(vii)$
$u(s)-c(\rho(s))-s=-\sqrt{\ds\f{\g-1}{2}}b_0\biggl(1+\a\biggr)<0$.

$(viii)$ $\rho'(s)=O(b_0^{-1})$.

$(ix)$ $u'(s)=-(n-1)\biggl(1+\a\biggr)<0$.}

{\bf Proof.}  Set $\rho_+=\lim\limits_{s\to s_0-}\rho(s),
u_+=\lim\limits_{s\to s_0-}u(s)$ and $\al_0=\ds\f{\rho_+}{\rho_0}$.

It follows from (2.2) that
$$\f{A\g}{\g-1}\biggl((\rho_+)^{\g+1}-\rho_0^{\g-1}(\rho_+)^2\biggr)
+\f{1}{2}s_0^2(\rho_+-\rho_0)^2-s_0^2\rho_+(\rho_+-\rho_0)=0.\tag
2.5$$

Set
$$F(x)=\f{A\g}{\g-1}\biggl(x^{\g+1}-\rho_0^{\g-1}x^2\biggr)
+\f{1}{2}s_0^2(x-\rho_0)^2-s_0^2 x(x-\rho_0).$$

Then (2.5) implies that $F(\rho_0)=F(\rho_+)=0$. For $x\in
(0,\rho_0)$, a direct computation yields
$$\align
F'(x)&=\f{A\g(\g+1)}{\g-1}x^{\g}-\f{2A\g}{\g-1}\rho_0^{\g-1}x-s_0^2
x\\
&\leq x\bigl(c^2(\rho_0)-s_0^2\bigr).
\endalign$$

Due to condition (2.3), $F'(x)<0$ holds for $x\in (0,\rho_0)$. This, together
with $F(\rho_0)=F(\rho_+)=0$,  yields $\rho_+>\rho_0$.
In this case, it can be derived from (2.5) that
$$\f{\g}{\g-1}A\al_0^2\f{(\al_0^{\g-1}-1)}{\al_0^2-1}=\f{1}{2}
\rho_0^{-\f{1}{\g-1}}s_0^2,\tag 2.6$$
where $\al>1$.

Since the left hand side of (2.6) is bounded if $\al_0>1$ is
bounded, then for large $s_0$, $\al_0$ is also large. From this
fact, one has
$$\al_0=\f{1}{\rho_0}\biggl(\f{\g-1}{2A\g}\biggr)^{\f{1}{\g-1}}s_0^{\f{2}{\g-1}}
\biggl(1+\a\biggr).\tag 2.7$$

Substituting this into (2.2) yields
$$\cases
\rho_+=\biggl(\ds\f{\g-1}{2A\g}\biggr)^{\f{1}{\g-1}}s_0^{\f{2}{\g-1}}\bigl(1
+\a\bigr),\\
 u_+=s_0\bigl(1+\a\bigr).
\endcases\tag 2.8$$

Moreover, from $u'(s)\leq 0$ for $s\in [b_0, s_0]$, one has $u_+\leq
u(s)$, namely,
$$s_0(1-\f{1}{\al_0})\leq u(s)\leq  b_0\leq s_0.$$

Combining this with (2.7) yields (i) and (ii).

When $s$ is taken a function of $u$, it follows from (2.1) that
$$\f{dh(\rho)}{d u}=s-u,$$
which derives
$$h(\rho)=h(\rho_+)+\int_{u_+}^{u}(s-\tau)d\tau.$$

This, together with (i)-(ii) and (2.8), yields (iii) and further
(iv)-(ix). Therefore, the proof of Lemma 2.1 is completed.\qquad \qquad\qquad\qquad
\qquad\qquad\qquad\qquad\qquad\qquad\qquad\qquad\qquad\qquad\qquad\qquad\qquad
\qquad\qquad\qquad\qquad\qed

\medskip

{\bf Remark 2.1.} {\it Since the denominator of system $(2.1)$
is negative in the interval $[b_0, s_0]$ $($see \rom{Lemma 2.1
(v)}$)$, one can extend the background solution $(\hat\rho(s), \hat
u(s))$ of $(2.1)$-$(2.4)$ to the interval $[b_0-\tau_0,
s_0+\tau_0]$ for some small positive constant $\tau_0$ satisfying
$0<\tau_0\leq b_0^{-\frac{4}{\gamma-1}}(s_0-b_0)$. In the following
sections, we will still denote this extension of the background
solution to $\{(t, r)\colon t>0,\, (b_0-\tau_0) t\leq r\leq
(s_0+\tau_0)t\}$ by $(\hat\rho(s), \hat u(s))$, where
$s=\ds\f{r}{t}$. The corresponding extension of the potential will
be denoted by $\hat\Phi(t,x)$. Moreover, $\hat \Phi$ can be written as
$t\hat\phi(\ds\f{r}{t})$ with
$\hat\phi'(s)=\hat u(s)$.}

\vskip 0.4 true cm \centerline{\bf $\S 3.$ Local-in-time existence}
\vskip 0.4 true cm

In this section, we give the local-in-time existence of the problem
(1.6) with (1.7)-(1.9) for $n=3$. When the spatial dimensions are two,
the existence result can be proved analogously and even much simpler.

\vskip 0.4 true cm
\centerline {\bf $\S 3.1.$ Reformulation of problem (1.6) with
(1.7)-(1.9)  under some nonlinear transformations}\vskip 0.4 true cm

Let $(\Phi^+, \zeta(t,\o))$ be the solution of the problem (1.6)
with (1.7)-(1.9), and we set
$$\phi=\f{\Phi^+}{t},\quad s=\f{r}{t},\quad b(t,\o)=\f{\sigma}{t},\quad \chi(t,\o)
=\f{\zeta}{t}.\tag 3.1$$

Under the coordinate transformation (1.5) and the notations in
(3.1), it follows from a direct but tedious computation that the
problem (1.6) with (1.7)-(1.9) for $n=3$ can be written as
$$\cases
\ds\p_t^2\phi+\f{2(\p_s\phi-s)}{t}\p_{ts}^2\phi+\f{2}{ts^2}\ds\sum_{i=1}^{3}Z_i\phi\p_t
Z_i\phi+\f{(\p_s\phi-s)^2}{t^2}\p_s^2\phi+\f{2(\p_s\phi-s)}{t^2
s^2}\ds\sum_{i=1}^{3}Z_i\phi\p_s Z_i\phi\\
\quad\ds+\f{1}{t^2 s^4}\ds\sum_{i,j=1}^{3}Z_i\phi Z_j\phi Z_i
Z_j\phi-\f{c^2(\rho)}{t^2}\biggl(\p_s^2\phi+\f{2}{s}\p_s\phi+\f{1}{s^2}
\ds\sum_{i=1}^{3}Z_i^2\phi\biggr)
+\f{2}{t}\p_t\phi\\
\quad\ds+\f{2s-\p_s\phi}{t^2
s^3}\ds\sum_{i=1}^{3}(Z_i\phi)^2=0,\qquad \quad (t,s,\o)\in (0,+\infty)\times
(b(t,\o),\chi(t,\o))\times \Bbb S^2,\\
\ds\p_s\phi-\f{1}{s^2}\ds\sum_{i=1}^{3}Z_i b\cdot Z_i\phi=t\p_t
b+b\qquad \qquad \qquad \qquad \qquad \qquad \qquad \qquad \qquad
\quad \text{on}\quad
s=b(t,\o),\\
\ds \Cal H(\phi,\na\phi)(\p_s\phi)^2+(\Cal H(\phi,\na\phi)-\rho_0)(t
\p_t\phi-s\p_s\phi)+\Cal
H(\phi,\na\phi)\f{1}{s^2}\ds\sum_{i=1}^{3}(Z_i\phi)^2=0\quad
\text{on}\quad s=\chi(t,\o),\\
\ds\phi=0 \qquad \qquad \qquad \qquad \qquad \qquad \qquad \qquad
\qquad \qquad \qquad \qquad \qquad \qquad \qquad \quad
\text{on}\quad s=\chi(t,\o)
\endcases\tag 3.2$$with
$$\ds\Cal
H(\phi,\na\phi)=h^{-1}\bigl(B_0-\phi-t\p_t\phi+s\p_s\phi-\f{1}{2}(\p_s\phi)^2
-\f{1}{2s^2}\ds\sum_{i=1}^{3}(Z_i\phi)^2\bigr)\tag
3.3$$and
$$c^2(\rho)=(\g-1)\biggl(B_0-\phi-t\p_t\phi+s\p_s\phi-\f{1}{2}(\p_s\phi)^2
-\f{1}{2s^2}\ds\sum_{i=1}^{3}(Z_i\phi)^2\biggr). \tag 3.4$$

It is noted that the problem (3.2) is still a free boundary problem with the
shock surface as the unknown boundary. In order to prove the
local-in-time existence result, we will use the modified partial hodograph transformation
to straighten
$s=b(t,\o)$ and $s=\chi(t,\o)$ simultaneously. To this end, motivated by [5] and [24], we set
$$\psi=s-b(t,\o)-\f{\phi}{b_0},\tag 3.5$$
and take the modified partial hodograph transformation as follows
$$T=t,\quad R=\f{s-b}{\psi}+1,\quad \o=\o.\tag 3.6$$


In this case, the
boundaries $s=b(t,\o)$ and $s=\chi(t,\o)$ are changed into $R=1$ and $R=2$
respectively. Additionally, by the continuity
condition $\phi(t,\chi(t,\o),\o)=0$ in (3.2) and (3.5) we obtain
$\chi(t,\o)=b(t,\o)+\psi(t,2,\o)$, which means that $\chi(t,\o)$ can be determined
once the function $\psi(t,2,\o)$ is known.

Next we derive the nonlinear equation and corresponding boundary
conditions of $\psi$.

It follows from (3.5)-(3.6) and a direct computation that
$$\cases
\ds s=a_0(b,\psi),\\
\ds \p_s\phi=b_0 a_1(b,\psi,\na\psi)\cdot a_2(b,\psi,\na\psi),\\
\ds \p_t\phi=-b_0 a_1(b,\psi,\na\psi)\cdot a_3(b,\psi,\na\psi),\\
\ds Z_i\phi=-b_0 a_1(b,\psi,\na\psi)\cdot
a_4^{i}(b,\psi,\na\psi),\quad  i=1,2,3
\endcases\tag 3.7$$
with
$$\cases
\ds a_0(b,\psi)=b+(R-1)\psi,\\
\ds a_1(b,\psi,\na\psi)=\f{1}{\psi+(R-1)\p_R\psi},\\
\ds a_2(b,\psi,\na\psi)=\psi+(R-2)\p_R\psi,\\
\ds a_3(b,\psi,\na\psi)=\psi\p_T\psi+\psi\p_T b+(R-2)\p_T b\p_R\psi,\\
\ds a_4^{i}(b,\psi,\na\psi)=\psi Z_i\psi+\psi Z_i b+(R-2) Z_i
b\p_R\psi,\quad  i=1,2,3.
\endcases\tag 3.8$$

Set $\Bbb H(b,\psi,\na\psi)\triangleq\Cal
H(\phi,\na\phi)$, then substituting (3.7) into (3.3)-(3.4) yields
$$
\cases
\Bbb H(b,\psi,\na\psi)=h^{-1}\bigl(A_0(b,\psi,\na\psi)\bigr),\\
c^2(\rho)=(\g-1)A_0(b,\psi,\na\psi),
\endcases\tag 3.9$$
where
$$A_0(b,\psi,\na\psi)=B_0-b_0 (R-2)\psi+Tb_0 a_1\cdot a_3+b_0a_0\cdot a_1\cdot
a_2-\f{b_0^2}{2}(a_1\cdot
a_2)^2-\f{b_0^2}{2a_0^2}\sum_{i=1}^{3}(a_1\cdot a_4^i)^2.$$

At this time, the problem (3.2) can be rewritten as
$$\align
&\ds
A_1(b,\psi,\na\psi)\p_T^2\psi+A_2(b,\psi,\na\psi)\p_{TR}^2\psi+\sum_{i=1}^{3}
A_3^i(b,\psi,\na\psi)\p_T
Z_i\psi+A_4(b,\psi,\na\psi)\p_R^2\psi\\
&\qquad \ds+\sum_{i=1}^{3}A_5^i(b,\psi,\na\psi)\p_R
Z_i\psi+\sum_{i=1}^{3}\sum_{j=1}^{3}A_6^{ij}(b, \psi,\na\psi)Z_i
Z_j\psi+A_7(b,\psi,\na\psi)=0,\\
&\ds\qquad\qquad\qquad\qquad\qquad\qquad\qquad\qquad\qquad
\qquad\quad (T,R,\o)
\in (0,+\infty)\times (1,2)\times\Bbb S^2,\tag 3.10\\
&\ds \p_R\psi+\f{1}{b_0}(T\p_T
b+b-b_0)\psi-\f{\psi}{b^2}\sum_{i=1}^{3}Z_i b\cdot
Z_i\psi-\f{\psi-\p_R\psi}{b^2}\sum_{i=1}^{3}(Z_i b)^2=0\qquad\qquad\qquad \quad
\text{on}\quad R=1,\tag 3.11\\
&\ds \Bbb H(b,\psi,\na\psi)\psi-\f{1}{b_0 a_1}(\Bbb
H(b,\psi,\na\psi)-\rho_0)(T\p_T a_0+a_0)+\f{\Bbb
H(b,\psi,\na\psi)\psi}{(b+\psi)^2}\sum_{i=1}^{3}(Z_i a_0)^2=0\quad
\text{on}\quad R=2,\tag 3.12
\endalign$$
where
$$
\cases
A_k(b,\psi,\na\psi)=A_{k,0}(b,\psi,\na\psi)+\ds\f{A_{k,1}(b,\psi,\na\psi)}{T}
+\ds\f{A_{k,2}(b,\psi,\na\psi)}{T^2},\qquad
k=1, 2, 4, 7,\\
A_k^i(b,\psi,\na\psi)=A_{k,0}^i(b,\psi,\na\psi)+\ds\f{A_{k,1}^i(b,\psi,\na\psi)}{T}
+\ds\f{A_{k,2}^i(b,\psi,\na\psi)}{T^2},\qquad
k=3, 5,\quad 1\le i\le 3,\\
A_6^{ij}(b,
\psi,\na\psi)=A_{6,0}^{ij}(b,\psi,\na\psi)+\ds\f{A_{6,1}^{ij}(b,\psi,\na\psi)}{T}
+\ds\f{A_{6,2}^{ij}(b,\psi,\na\psi)}{T^2},\qquad
1\le i, j\le 3,\\
\endcases
$$
and all expressions of $A_{k,0}(b,\psi,\na\psi)$, ... , $A_{6,2}^{ij}(b,\psi,\na\psi)$
can be found in Lemma A.1-Lemma A.3 of Appendix A since their concrete expressions are required
in order to establish
related estimates.

Therefore, the problem (1.6) with (1.7)-(1.9) has been reformulated into the
problem (3.10)-(3.12). To establish the local-in-time existence result of
(1.6) with (1.7)-(1.9), we require to show the local existence of
(3.10) with (3.11)-(3.12), which can be stated as \medskip

{\bf Theorem 3.1.} {\it Under the assumptions in \rom{Theorem 1.1}, there exists a positive
constant $T^*$ independent of small $\ve$, such that the problem $(3.10)$
with $(3.11)$-$(3.12)$ has a unique solution $\psi(T,R,\o)\in
C^{\infty}([0, T^*]\times [1,2]\times\Bbb S^2)$, which satisfies for
$m\in\Bbb N\cup\{0\}$
$$\|\psi-\hat\psi\|_{C^{m}([0,T^*]\times [1,2]\times \Bbb S^2)}\leq C(m)\ve,$$
where the function $\hat\psi=\hat\psi(R)$ is obtained when we use
$b_0$ and $\hat\phi(s)$ give in \rom{Remark 2.1} instead of
$b(t,\o)$ and $\phi(t,s,\o)$ in \rom{(3.5)-(3.6)}.}\medskip

{\bf Remark 3.1.} {\it It is noted that there is no  initial data
information on $T=0$ of $\psi$ in the problem \rom{(3.10)} with
\rom{(3.11)-(3.12)}. Thus, solving \rom{(3.10)-(3.12)} is not
a standard procedure. In addition, by the proof procedure of Theorem 3.1,
we can derive that $T^*\ge\ds\f{C}{\ve}$ is actually permitted.}
\medskip

For later uses, we now establish some properties of
$\hat\psi(R)$.

{\bf Lemma 3.2.}  {\it $\hat\psi(R)$ admits the following estimates:

$(i).$ $\hat\psi=(s_0-b_0)(1+\a)>0,$

$(ii).$
$\hat\psi'(R)=\ds\f{1}{b_0}O\biggl((s_0-b_0)^2\biggr)>0,\quad
\hat\psi'(2)=\ds\f{2}{b_0}(s_0-b_0)^2(1+\a),$

$(iii).$  $\ds\hat\psi''(R)=\f{2}{b_0}(s_0-b_0)^2(1+\a).$}\medskip

{\bf Proof.} Due to $\hat\phi(s_0)=0$, then by (3.5) and mean value theorem there exists some
$\bar s\in (b_0, s_0)$ such that
$$
\hat\psi=s-b_0-\f{\hat\phi}{b_0}=s-b_0+\f{1}{b_0}\hat u(\bar s)(s_0-s).
$$

This, together with Lemma 2.1 (ii), yields Lemma 3.2 (i).

By the second equality in (3.7), one has
$$\left(b_0+(1-R)(b_0-\hat u(s))\right)\hat\psi'(R)=(b_0-\hat u(s))\hat\psi.\tag3.13$$

Noting $\hat u(b_0)=b_0$ and $b_0-\hat u(s)=\hat u(b_0)-\hat u(s)$, then combining this
with mean value theorem,
Lemma 2.1 (ix),  Lemma 3.2 (i) and (3.13) yields Lemma 3.2 (ii).

Similarly,  Lemma 3.2 (iii) can be obtained by taking the first order derivative on two
hand sides of (3.13)
and applying Lemma 2.1.\qquad\qquad\qquad\qquad\qquad\qquad\qquad\qquad\qquad\qquad
\qquad\qquad\qquad\qquad\qquad\qquad\qed

Since $\hat\psi(R)$ does not satisfy the fixed wall boundary
condition (3.11) for $b=b(t,\o)$, this will bring some troubles to
study the nonlinear problem (3.10) with Neumann-type boundary
condition (3.11). This difficulty can be overcome by choosing a new
function  $\hat\psi_{\si}(T,R,\o)$ instead of $\hat\psi(R)$ as
follows:
\medskip

{\bf Lemma 3.3.}  {\it Under the assumptions in \rom{Theorem 1.1},
we define a $C^{\infty}$ function $\hat\psi_{\si}(T,R,\o)$  as
follows
$$\hat\psi_{\si}(R,\o)=\ds\f{\hat\psi(R)}{b^2+\ds\sum_{i=1}^{3}(Z_i b)^2}\left(b^2
+\ds\sum_{i=1}^{3}(Z_i b)^2+(R-1)\ds\bigl(\ds\sum_{i=1}^{3}(Z_i
b)^2-\f{b^2}{b_0}(T\p_T b+b-b_0)\bigr) \right),$$ which satisfies

$(i)$  $\ds\p_R\hat\psi_{\si}+\f{1}{b_0}(T\p_T
b+b-b_0)\hat\psi_{\si}-\f{\hat\psi_{\si}}{b^2}\ds\sum_{i=1}^{3}Z_i
b\cdot
Z_i\hat\psi_{\si}-\f{\hat\psi_{\si}-\p_R\hat\psi_{\si}}{b^2}\ds\sum_{i=1}^{3}(Z_i
b)^2=0.$

$(ii)$ $\|\hat\psi_{\si}-\hat\psi\|_{C^m ([0,\infty)\times
[1,2]\times \Bbb S^2)}\leq C(m)\ve,\quad \forall\ m\in\Bbb
N\cup\{0\}.$}\medskip

{\bf Proof.}  Since (i) and (ii) can be verified directly with the
definition of $\hat\psi_{\si}(R,\o)$ and the assumptions on $b$ in Theorem
1.1, then we omit the details here.\qquad\qquad\qquad\qquad\qquad\qquad\qquad\qquad
\qquad\qquad\qquad\qquad\qquad\qed

In subsequent part, due to Remark 3.1, we start to construct an
approximate solution of (3.10) with (3.11)-(3.12) and give some
related estimates.

\vskip 0.5 true cm

{\bf $\S  3.2.$ Construction of an approximate solution to
problem (3.10) with (3.11)-(3.12)}

\vskip 0.5 true cm

We will use Taylor's formula to construct an approximate solution
of (3.10) with (3.11)-(3.12). For this end, we set
$$\psi(T,R,\o)=\psi_0(R,\o)+\sum_{l=1}^{K}T^{l}\psi_{l}(R,\o)+O(T^{K+1})\tag 3.14$$
and
$$b(T,\o)=d_0(\o)+\sum_{l=1}^{K}T^{l}b_l(\o)+O(T^{K+1}),\tag 3.15$$
where $K$ is a suitably  large positive integer, and $d_0(\o)=b(0,\o)$.

Substituting (3.14)-(3.15) into the
problem (3.10) with (3.11)-(3.12) and comparing  the
powers of $T$ in the resulting equalities, one can get a series of equations and boundary
conditions of $\psi_l$ for $0\le l\le K$.
This will be illustrated gradually in subsequent parts.

\vskip 0.2 true cm

{\bf Part 1. Determination of $\psi_0(R,\o)$}

\vskip 0.2 true cm

By comparing the coefficients of $T^{-2}$ and the coefficients of
$T^{0}$ in the resulting equalities from (3.10) and (3.11)-(3.12)
respectively by the expressions (3.14) and (3.15), one can arrive at
$$\cases
\ds A_{4,2}(d_0,\psi_0,\na\psi_0)\p_R^2\psi_0
+\sum_{i=1}^{3}A_{5,2}^i(d_0,\psi_0,\na\psi_0)\p_R Z_i\psi_0
+\sum_{i=1}^{3}\sum_{j=1}^{3}A_{6,2}^{ij}(d_0,\psi_0,\na\psi_0)Z_i Z_j \psi_0\\
\ds\qquad\qquad\qquad+A_{7,2}(d_0,\psi_0,\na\psi_0)=0,\quad\quad\quad\quad\quad\quad\quad
\quad\quad
\quad\quad\quad (R, \o)\in (1,2)\times \Bbb S^{2},\\
\ds\p_R\psi_0+\f{1}{b_0}(d_0-b_0)\psi_0-\f{\psi_0}{d_0^2}\ds\sum_{i=1}^{3}Z_i
d_0\cdot Z_i\psi_0
-\f{\psi_0-\p_R\psi_0}{d_0^2}\ds\sum_{i=1}^{3}(Z_i d_0)^2=0\quad \text{on}\quad R=1,\\
\ds
\Bbb H(d_0,\psi_0,\na\psi_0)\psi_0-\f{1}{b_0}(\Bbb H(d_0,\psi_0,\na\psi_0)-\rho_0)(\psi_0
+\p_R\psi_0)(d_0+\psi_0)\\
\ds\qquad\qquad\qquad+\f{\Bbb
H(d_0,\psi_0,\na\psi_0)\psi_0}{(d_0+\psi_0)^2}\ds\sum_{i=1}^{3}(Z_i
d_0+Z_i \psi_0)^2=0\quad\quad\quad\quad\quad\quad\quad
\text{on}\quad R=2.
\endcases\tag 3.16$$

Later on, the equation in (3.16) can be shown to be a second order nonlinear elliptic equation
(see Lemma 3.4 below).

On the other hand, with respect to the function $\hat\psi(R)$
given in Lemma 3.2, we have
$$\cases
A_{4,2}(b_0,\hat\psi,\na\hat\psi)\hat\psi''(R)+A_{7,2}(b_0,\hat\psi,\na\hat\psi)=0,\qquad\qquad
\quad R\in (1,2),\\
\hat\psi'(1)=0,\\
\ds \Bbb H(b_0,\hat\psi,\na\hat\psi)\hat\psi -\f{1}{b_0}(\Bbb
H(b_0,\hat\psi,\na\hat\psi)-\rho_0)
(\hat\psi+\hat\psi'(2))(b_0+\hat\psi)=0\quad \text{on}\quad
R=2.\endcases\tag 3.17$$

Set
$\bar\psi=\psi_0-\hat\psi$, then it follows from (3.16)-(3.17) that $\bar\psi$ satisfies
$$\cases
\ds
A_{4,2}(d_0,\psi_0,\na\psi_0)\p_R^2\bar\psi+\ds\sum_{i=1}^{3}A_{5,2}^{i}(d_0,\psi_0,\na\psi_0)\p_R
Z_i\bar\psi+\ds\sum_{i=1}^{3}\sum_{j=1}^{3}A_{6,2}^{ij}(d_0,\psi_0,\na\psi_0)Z_i
Z_j\bar\psi\\
\quad\quad\ds+B_1\p_R\bar\psi+\ds\sum_{i=1}^{3} B_2^{i}Z_i\bar\psi+
E_0\bar\psi=F_0(R,\o),\qquad\qquad\qquad\qquad \quad (R,\o)\in (1,2)\times \Bbb S^2,\\
D_{11}\p_R\bar\psi+D_{12}\bar\psi+\ds\sum_{i=1}^{3}D_{13}^{i}Z_i\bar\psi=G_1(R,\o)\qquad
\qquad\qquad\qquad\qquad\qquad\quad
\text{on}\quad R=1,\\
\ds
D_{21}^0\p_R\bar\psi+D_{22}^0\bar\psi+\ds\sum_{i=1}^{3}D_{23}^{0,i}Z_i\bar\psi=G_2(R,\o)
\qquad\qquad\qquad\qquad\qquad\qquad\quad \text{on}\quad R=2
\endcases\tag 3.18$$
with
$$\cases
F_0(R,\o)=\left(A_{4,2}(b_0,\hat\psi,\na\hat\psi)-A_{4,2}(d_0,\hat\psi,\na\hat\psi)\right)
\hat\psi''(R)+\left(A_{7,2}(b_0,\hat\psi,\na\hat\psi)-A_{7,2}(d_0,\hat\psi,\na\hat\psi)\right),\\
\ds G_1(R,\o)=\f{1}{b_0}
(b_0-d_0)\hat\psi+\f{\hat\psi-\p_R\hat\psi}{d_0^2}\ds\sum_{i=1}^{3}(Z_i
d_0)^2,\\
G_2(R,\o)=\ds\sum_{i=1}^{3}O(Z_i(d_0-b_0))+O(d_0-b_0)
\endcases\tag 3.19$$
and
$$\ds D_{11}=1+\f{1}{d_0^2}\ds\sum_{i=1}^{3}(Z_i d_0)^2,\quad
D_{12}=\f{1}{b_0}(d_0-b_0)-\f{1}{d_0^2}\ds\sum_{i=1}^{3}(Z_i
d_0)^2,\quad D_{13}^i(\psi_0)=-\f{\psi_0}{d_0^2}Z_i d_0\tag
3.20$$ and
$$\cases
\ds E_0(R,\o,\bar\psi,\na\bar\psi)
=\ds\int_0^1\p_{\psi}\bigl(\hat\psi''(R)
A_{4,2}+A_{7,2}\bigr)(d_0,\hat\psi+\theta\bar\psi,\na(\hat\psi+\theta\bar\psi))d\theta,\\
(D_{21}^0,D_{22}^0)(R,\o,\bar\psi,\na\bar\psi)
=\ds\int_0^1\na_{\p_R\psi, \psi}\biggl(\Bbb
H(\cdot)\psi-\f{1}{b_0}(\Bbb H(\cdot)
-\rho_0)(\psi+\p_R\psi)(d_0+\psi)\\
\qquad \qquad\qquad\qquad\qquad \qquad \ds+\f{\Bbb
H(\cdot)\psi}{(d_0 +\psi)^2}\ds\sum_{i=1}^{3}(Z_i
d_0)^2\biggr)(d_0,\hat\psi+\theta\bar\psi,\na(\hat\psi+\theta\bar\psi))d\theta.\\
\endcases\tag 3.21$$
In addition, the coefficients $B_1, B_2^i, D_{23}^{0,i}$ are smooth  functions on the
variables ($R, \o, \bar\psi, \na\bar\psi$),
whose concrete expressions are not required.

For the requirements to solve the nonlinear problem (3.18), we now give some estimates
on the coefficients
$A_4$, $A_5^i$ and $A_6^{ij}$ in the equation of (3.18) when their arguments $(d_0(\o),
\psi_0(R,\o),\na\psi_0(R,\o))$
are replaced by $(b_0,\hat\psi(R),\na\hat\psi(R))$.
\medskip

{\bf Lemma 3.4.} {\it For large $b_0$, we have
$$\cases
\ds A_{4,2}(b_0,\hat\psi,\na\hat\psi)=\ds-\f{b_0^2}{2(s_0-b_0)}(1+\a)<0,\\
\ds A_{5,2}^i(b_0,\hat\psi,\na\hat\psi)=0,\qquad i=1, 2, 3, \\
\ds
A_{6,2}^{ij}(b_0,\hat\psi,\na\hat\psi)=-\f{(\g-1)(s_0-b_0)}{2}(1+\a)\dl_{ij},\qquad
1\le i, j\le 3 ,\\
\endcases\tag 3.22$$
which means that the equation in \rom{(3.18)} is uniformly
elliptic.}

{\bf Proof.} By Lemma 3.2 and (3.8), a direct computation yields
$$\cases
a_0(b_0,\hat\psi)=b_0(1+\a),\\
a_1(b_0,\hat\psi,\na\hat\psi)=\ds\f{1}{s_0-b_0}(1+\a),\\
a_2(b_0,\hat\psi,\na\hat\psi)=(s_0-b_0)(1+\a),\\
a_3(b_0,\hat\psi,\na\hat\psi)=0,\\
a_4^i(b_0,\hat\psi,\na\hat\psi)=0,\qquad  i=1,2,3.\\
\endcases\tag 3.23$$

Then in terms of the concrete expressions of $A_{4,2}$, $A_{5,2}^i$
and $A_{6,2}^{ij}$ in Lemma A.3, with (3.23), we have
$$\align
A_{4,2}(b_0,\hat\psi,\na\hat\psi)&=\bigl(c^2(\rho)-(b_0 a_1 a_2-a_0)^2\bigr)\bigl((R-2)a_1-(R-1)a_1^2 a_2\bigr)\biggl|_{b=b_0;\psi=\hat\psi}\\
&=-\f{\g-1}{2(s_0-b_0)}b_0^2(1+\a)<0,\\
A_{5,2}^i(b_0,\hat\psi,\na\hat\psi)&=0,\qquad i=1,2,3,\\
A_{6,2}^{ij}(b_0,\hat\psi,\na\hat\psi)&=\ds-\f{c^2(\rho)}{a_0^2}\psi\delta_{ij}\biggl|_{b=b_0;\psi=\hat\psi}\\
&=\ds-\f{(\g-1)(s_0-b_0)}{2}\delta_{ij}(1+\a),\qquad  1\le i,j\le 3.
\endalign$$

Thus (3.22) is proved.\qquad\qquad\qquad\qquad\qquad\qquad\qquad\qquad\qquad
\qquad\qquad\qquad\qquad\qquad\qquad\qed

We now establish the solvability of the problem (3.18).

{\bf Lemma 3.5.}  {\it Under the assumptions in \rom{Theorem 1.1},
the problem $(3.18)$ has a unique smooth solution $\bar\psi$ which
satisfies
$$\|\bar\psi\|_{C^{3+m,\al}([1, 2]\times\Bbb S^2)}\leq C(m)\ve,\tag 3.24$$
where $0<\al<1$ is any fixed constant, and $m\in\Bbb N\cup\{0\}$.}

{\bf Proof.} At first, we claim that for any $\bar\psi_0$ satisfying
$\|\bar\psi_0\|_{C^{3+m,\al}([1, 2]\times\Bbb S^2)}\leq C(m)\ve,$
there exists a sequence $\{\bar\psi_l\}_{l=1,2,\cdots}$ such that
for
$V_l\equiv(\hat\psi+\bar\psi_{l-1},\na(\hat\psi+\bar\psi_{l-1}))$,
$\bar\psi_l$ satisfies
$$\cases
\ds
A_{4,2}(d_0,V_l)\p_R^2\bar\psi_l+\ds\sum_{i=1}^{3}A_{5,2}^{i}(d_0,V_l)\p_R
Z_i\bar\psi_l+\ds\sum_{i=1}^{3}\sum_{j=1}^{3}A_{6,2}^{ij}(d_0,V_l)Z_i
Z_j\bar\psi_l+B_{1}(d_0,V_l)\p_R\bar\psi_l\\
\qquad \quad\ds+\ds\sum_{i=1}^{3} B_{2}^{i}(d_0,V_l)Z_i\bar\psi_l+
E_0(d_0,V_l)\bar\psi_l=F_0(R,\o),\qquad \qquad \quad  (R,\o)\in (1,2)\times \Bbb S^2,\\
\ds D_{11}\p_R\bar\psi_l+D_{12}\bar\psi_{l}
+\ds\sum_{i=1}^{3}D_{13}^i(\hat\psi+\bar\psi_{l-1})Z_i\bar\psi_l=G_1(R,\o)\qquad\qquad
\qquad \quad
\text{on}\quad R=1,\\
\ds D_{21}^0(d_0,V_l)\p_R\bar\psi_l
+D_{22}^0(d_0,V_l)\bar\psi_l+\ds\sum_{i=1}^{3}D_{23}^{0,i}(d_0,V_l)Z_i\bar\psi_l=G_2(R,\o)
\qquad \quad \text{on}\quad R=2
\endcases\tag 3.25$$
and admits the following estimate
$$\|\bar\psi_l\|_{C^{3+m,\al}([1, 2]\times\Bbb S^2)}\leq C(m)\ve.\tag 3.26$$

We will use the induction method to prove the claim (3.26).

Assume that $\bar\psi_{l-1}$ satisfies (3.26), then it follows from Lemma A.4-Lemma A.5
in Appendix A that there exists a positive
constant $c_0$ such that
$$E_0(d_0,V_{l})>c_0>0,\quad
D_{2i}^0(d_0,V_{l})<-c_0<0\quad \text{for}\quad  i=1,2.\tag
3.27$$

Choosing a function
$\vartheta(R)=2+c_1 R$ with
$$c_1=\f{1}{4(1+2\|A_4\|_{L^{\infty}}+\|B_1\|_{L^{\infty}}+\|D_{21}^0\|_{L^{\infty}})}$$
and setting $\bar\Psi_l(R,\o)\equiv\vartheta(R)\bar\psi_l$, then $\bar\Psi_l(R,\o)$ satisfies
$$\cases
A_{4,2}(d_0,V_l)\p_R^2\bar\Psi_l+\ds\sum_{i=1}^{3}A_{5,2}^{i}(d_0,V_l)\p_R
Z_i\bar\Psi_l
+\ds\sum_{i,j=1}^{3}A_{6,2}^{ij}(d_0,V_l)Z_i Z_j\bar\Psi_l+\t B_{1}(d_0,V_l)\p_R\bar\Psi_l\\
\quad\quad+\ds\sum_{i=1}^{3} \t B_{2}^i(d_0,V_l)Z_i\bar\Psi_l+\t
E_0(d_0,V_l)\bar\Psi_l
=\vartheta(R)F_0(R,\o),\qquad \quad (R,\o)\in (1,2)\times \Bbb S^2,\\
\ds D_{11}\p_R\bar\Psi_l+\t D_{12}\bar\Psi_l +\ds\sum_{i=1}^{3}\t
D_{13}^{i}(\hat\psi+\bar\psi_{l-1})Z_i\Psi_l=\vartheta(R)G_1(R,\o)
\qquad\qquad \quad\quad
\text{on}\quad R=1,\\
\ds D_{21}^0(d_0,V_l)\p_R\bar\psi_l+\t D_{22}^0(d_0,V_l)\bar\psi_l
+\ds\sum_{i=1}^{3}\t D_{23}^{0,i}(d_0,V_l)Z_i\bar\psi_l
=\vartheta(R)G_2(R,\o)\quad \text{on}\quad R=2.
\endcases\tag 3.28$$

Then with the expression of $c_1$, (3.20), (3.27) and the
assumptions in Theorem 3.1, one has
$$\align
&\t E_0(d_0,V_{l})
=\f{2A_{4,2}}{\vartheta^2}(c_1)^2-\f{c_1 B_1}{\vartheta}+E_0
\geq\f{c_0}{2}>0,\\
&\t D_{12}
=-\f{c_1}{\vartheta}D_{11}+D_{12}
\leq -\f{c_1}{2}
<0,\\
&\t D_{22}^0(d_0,V_{l})
=\f{c_1}{\vartheta}D_{21}^0+D_{22}^0
\leq -\f{c_0}{2}<0.
\endalign$$

By Lemma 3.4 and Theorem 6.30-Theorem 6.31 in [9] that (3.28) has a
unique solution $\bar\Psi_l$ and then (3.25) has a unique solution
$\bar\psi_l$. With the expressions in (3.19) and the assumptions in
Theorem 3.1, $\bar\psi_l$ satisfies
$$
\align \|\bar\psi_l\|_{C^{3+m,\al}}&\leq
C(m)(\|F_0\|_{C^{2+m,\al}}+\|G_1\|_{C^{2+m,\al}}+\|G_2\|_{C^{2+m,\al}})
\biggl(1+\f{1}{b_0}\|\bar\psi_{l-1}
\|_{C^{3+m,\al}}\biggr)\\
&\leq C(m)\|d_0-b_0\|_{C^{2+m,\al}}\le C(m)\ve,\tag3.29
\endalign
$$
here we have used the assumption (3.26) in the case of $l-1$ and largeness of
$b_0$. Thus, the claim (3.26) is proved.

On the other hand, similar to proof of (3.29), we have for
small $\ve$
$$
\|\bar\psi_l-\bar\psi_{l-1}\|_{C^{2,\al}}\leq C
\ve\|\bar\psi_{l-1}-\bar\psi_{l-2}\|_{C^{2,\al}} \leq
\f{1}{2}\|\bar\psi_{l-1}-\bar\psi_{l-2}\|_{C^{2,\al}}.\tag3.30
$$

Combining (3.29) with (3.30) yields that there exists a function $\bar\psi\in C^{3+m,\al}$
such that $\bar\psi_l\to \bar\psi$ in $C^{2,\al}$ as $l\to\infty$. Furthermore, (3.24) holds.
Thus the proof on Lemma 3.5 is
completed.\qquad\qquad \qquad \qquad \qquad  \qquad\qed

\medskip

With Lemma 3.5, then $\psi_0=\hat\psi+\bar\psi$ is a unique solution of
(3.16). Based on this, next we continue to construct the
approximate solution of (3.10) with (3.11)-(3.12).\medskip

{\bf Part 2. Determination of $\psi_k$ $(k\geq 1)$}\medskip

Comparing the coefficients of $T^{k-2}$ and the coefficients of $T^{k}$ in the resulting equalities
from the equation (3.10) and
(3.11)-(3.12) respectively by the expressions
(3.14) and (3.15), then $\psi_k$ satisfies
$$\cases
\ds
A_{4,2}(d_0,V_1)\p_R^2\psi_k+\sum_{i=1}^{3}A_{5,2}^{i}(d_0,V_1)\p_R
Z_i\psi_k+\sum_{i=1}^{3}\sum_{j=1}^{3}A_{6,2}^{ij}(d_0,\psi_0,V_1)Z_iZ_j\psi_k\\
\ds\qquad+B_{k1}(d_0,V_1)\p_R\psi_k
+\sum_{i=1}^{3}B_{k2}^{i}(d_0,V_1)
Z_i\psi_k+E_k(d_0,V_1)\psi_k=F_k(\psi_l)_{0\leq l\le k-1},\\
\qquad\qquad\qquad\qquad\qquad\qquad\qquad\qquad\qquad\qquad\qquad\qquad (R,\o)\in (1,2)
\times \Bbb S^2,\\
\ds
D_{11}\p_R\psi_k+D_{12}\psi_k+\ds\sum_{i=1}^{3}D_{13}^{i}(\psi_0)Z_i\psi_k
=G_1^k(\psi_l)_{0\leq
l\le k-1}
\qquad\qquad \qquad\qquad\quad\quad \text{on}\quad R=1,\\
\ds D_{21}^k(d_0,V_1)\cdot\p_R
\psi_k+D_{22}^k(d_0,V_1)\cdot\psi_k+\ds\sum_{i=1}^{3}D_{23}^{k,i}(d_0,V_1)Z_i\psi_k
=G_2^k(\psi_l)_{0\leq
l\le k-1}\quad \text{on}\quad R=2
\endcases\tag 3.31$$
with
$$\align
E_k(d_0,V_1)=&k(k-1)\psi_0+\p_{\psi}\bigl(\hat\psi''(R)\cdot
A_{4,2}+A_{7,2}\bigr)(b_0,\hat\psi,\na\hat\psi)\\
&+k\p_{\p_T\psi}\bigl(\hat\psi''(R)
A_{4,1}+A_{7,1}+\f{\hat\psi''(R)}{T}
A_{4,2}+\f{A_{7,2}}{T}\bigr)(b_0,\hat\psi,\na\hat\psi)+O_k(\ve)\tag 3.32
\endalign$$
and
$$\cases
D_{21}^k(d_0,V_1)=\p_{\p_R\psi} \biggl(\Bbb
H(\cdot)\hat\psi-\f{1}{b_0}(\Bbb H(\cdot)-\rho_0)a_1
(T\p_T a_0+a_0)\biggr)(b_0,\hat\psi,\na\hat\psi)+O_k(\ve),\\
D_{22}^k(d_0,V_1)=\bigl(\p_{\psi}+k\p_{\p_T\psi}\bigr) \biggl(\Bbb
H(\cdot)\psi-\ds\f{1}{b_0}(\Bbb H(\cdot)-\rho_0) (T\p_T
a_0+a_0)a_1\biggr)(b_0,\hat\psi,\na\hat\psi)+O_k(\ve).
\endcases\tag 3.33$$
Here the error term estimates in (3.31)-(3.33) come from the assumptions on
$b(t,\o)$ in Theorem 1.1 and (3.24), and $O_k(\ve)$ stands for a generic
quantity satisfying
$|O_k(\ve)|\le C_k\ve$.

With respect to the problem (3.31), we have\medskip

{\bf Lemma 3.6.}  {The problem (3.31) has a unique smooth solution
$\psi_k$, which satisfies
$$\|\psi_k\|_{C^{3+m,\al}([1,2]\times \Bbb S^2)}\leq C(m,k)\ve.\tag 3.34$$}\medskip

{\bf Proof.} In terms of Lemma A.4-Lemma A.5 in Appendix A, we have
$$E_k(d_0,V_1)>c_0>0,\quad D_{21}^k(d_0,V_1)<-c_0<0,\quad
D_{22}^k(d_0,V_1)<-c_0<0.$$

Then similar to the treatment on the problem (3.25), one can derive that the linear
problem (3.31) has a unique
smooth solution $\psi_k$ satisfying (3.34). Consequently, the proof of Lemma
3.6 is completed.

\qquad \qquad\qquad \qquad \qquad\qquad\qquad \qquad\qquad\qquad \qquad\qquad\qquad
\qquad\qquad
\qquad \qquad\qquad\qquad \qquad\qquad\qed\medskip

{\bf $\S 3.3.$ Choice of iteration scheme and proof of local existence}\medskip

 {\bf $\S
3.3.1.$ Choice of iteration scheme}\medskip

By the construction of $\{\psi_l\}_{l\geq 0}$ in $\S 3.2$, for any
given integer $\kappa_0\geq 1$, we set
$\ds\psi_a^{\kappa_0}=\sum_{l=0}^{{\kappa_0}}T^{l}\psi_l$, which can
be regarded as the ${\kappa_0}-$order approximate solution of the problem
(3.10) with (3.11)-(3.12). By Lemma 3.5 and Lemma 3.6, there exists
$T^*>0$  depending only on ${\kappa_0}$ such that
$$\|\psi_a^{\kappa_0}-\hat\psi\|_{C^{m}([0,T^*]\times [1,2]\times\Bbb S^2)}\leq
C(m)\ve.\tag 3.35$$

Moreover, by (3.14)-(3.16), (3.31) and (3.35), $\psi_a^{\kappa_0}$
satisfies
$$\cases
&\ds
A_1(b,V_a^{\kappa_0})\p_T^2\psi_a^{\kappa_0}
+A_2(b,V_a^{\kappa_0})\p_{TR}^2\psi_a^{\kappa_0}+\sum_{i=1}^{3}A_3^i(b,V_a^{\kappa_0})\p_T
Z_i\psi_a^{\kappa_0}+A_4(b,V_a^{\kappa_0})\p_{R}^2\psi_a^{\kappa_0}\\
&\qquad \ds+\sum_{i=1}^{3}A_5^i(b,V_a^{\kappa_0})\p_{R}
Z_i\psi_a^{\kappa_0}+\sum_{i=1}^{3}\sum_{j=1}^{3}A_6^{ij}(b,V_a^{\kappa_0})Z_i
Z_j\psi_a^{\kappa_0}+A_7(b,V_a^{\kappa_0})\\
&\qquad=T^{{\kappa_0}-1}F^{\kappa_0},\qquad\qquad\qquad\qquad\qquad\qquad\qquad (T, R,\o)
\in (0, T^*)\times (1,2)\times \Bbb S^2,\\
&\ds \p_R\psi_a^{\kappa_0}+\f{1}{b_0}(T\p_T
b+b-b_0)\psi_a^{\kappa_0}-\f{\psi_a^{\kappa_0}}{b^2}\sum_{i=1}^{3}Z_i
b\cdot
Z_i\psi_a^{\kappa_0}-\f{\psi_a^{\kappa_0}-\p_R\psi_a^{\kappa_0}}{b^2}\sum_{i=1}^{3}(
Z_i b)^2=T^{{\kappa_0}+1}G_1^{\kappa_0}\\
&\qquad\qquad\qquad\qquad\qquad\quad\quad\quad\quad\quad\quad\quad\quad
\qquad\qquad\qquad\qquad\qquad\qquad\qquad  \text{on}\quad R=1,\\
&\ds \Bbb H(b,V_a^{\kappa_0})\psi_a^{\kappa_0}-\f{1}{b_0}(\Bbb
H(b,V_a^{\kappa_0})-\rho_0)a_2(b,V_a^{\kappa_0})(T\p_T
a_0+a_0)(b,V_a^{\kappa_0})\\
&\quad\quad\quad\quad\ds+\f{\Bbb
H(b,V_a^{\kappa_0})\psi_a^{\kappa_0}}{(b+\psi_a^{\kappa_0})^2}\sum_{i=1}^{3}(Z_i
a_0(b,V_a^{\kappa_0}))^2=T^{{\kappa_0}+1}G_2^{\kappa_0}\quad\quad\qquad\qquad\quad
\text{on}\quad R=2
\endcases\tag 3.36$$
with $V_a^{\kappa_0}=(\psi_a^{\kappa_0},\na\psi_a^{\kappa_0})$ and
$$\|(F^{\kappa_0}, G_1^{\kappa_0}, G_2^{\kappa_0})\|_{C^{m}([0,T^*]\times [1,2]\times\Bbb
S^2)}\leq C(m,{\kappa_0})\ve.\tag 3.37$$

Let $\dot\psi=\psi-\psi_a^{\kappa_0}$, $X=\ln T$ and $Y=R\o$. With the help of
(3.36), the problem (3.10) with (3.11)-(3.12) is equivalent to
$$\cases
\ds \Cal A_{1}(\dot\psi,\na\dot\psi)\p_X^2\dot\psi+2\Cal
A_{2}(\dot\psi,\na\dot\psi)\p_{XR}^2 \dot\psi+2\sum_{i=1}^{3}\Cal
A_{3}^i(\dot\psi,\na\dot\psi)\p_X Z_i \dot\psi+\Cal
A_4(\dot\psi,\na\dot\psi)\p_R^2\dot\psi\\
\ds\quad +2\sum_{i=1}^{3}\Cal A_5^i(\dot\psi,\na\dot\psi)\p_R Z_i
\dot\psi+\sum_{i=1}^{3}\sum_{j=1}^{3}\Cal
A_{6}^{ij}(\dot\psi,\na\dot\psi)Z_i Z_j\dot\psi=
\dot f\left(e^{({\kappa_0}+1)X}F^{\kappa_0},\dot\psi,\na_{X,Y}\dot\psi\right),\\
\qquad\qquad\qquad\qquad\qquad\qquad\qquad\qquad (X,R,\o)\in
(-\infty,
X_0]\times (1,2)\times \Bbb S^2,\\
\ds \Cal B_{11}\p_R\dot\psi+\sum_{i=1}^3\Cal B_{12}^i(\dot\psi)
Z_i\dot\psi=\dot
g_1\left(e^{({\kappa_0}+1)X}G_1^{\kappa_0},\dot\psi\right)
\qquad\qquad \quad \text{on}\quad R=1,\\
\ds \Cal B_{20}(\dot\psi,\na\dot\psi)\p_X\dot\psi+\Cal
B_{21}(\dot\psi,\na\dot\psi)\p_R\dot\psi+\sum_{i=1}^{3}\Cal
B_{22}^i(\dot\psi,\na\dot\psi)Z_i\dot\psi\\
\qquad\qquad\qquad\qquad\qquad\qquad=\dot
g_2\left(e^{({\kappa_0}+1)X}G_2^{\kappa_0},\dot\psi,\na_{X,Y}\dot\psi\right)\quad
\text{on}\quad R=2,
\endcases\tag 3.38$$
here $X_0=\ln T^*$, and the smooth functions $\dot f$, $\dot g_1$ and $\dot g_2$ satisfy
$$\dot f(0)=\dot g_1(0)=\dot g_2(0)=0,\quad \p_{\dot\psi}\dot g_1(0)=O(\ve),
\quad \na_{\p_X\dot\psi, \p_Y\dot\psi}\dot g_2(0)=0.\tag 3.39$$

In addition,
$$\cases
\ds (\Cal A_1, \Cal A_2, \Cal
A_3^i)=\f{\psi_a^{\kappa_0}+\dot\psi}{2(\g-1)A_0(b,V)}
(2A_1(b,V), e^{X}A_2(b,V), e^{X}A_3^i(b,V)),\qquad i=1,2,3,\\
\ds (\Cal A_4, \Cal A_5^i, \Cal
A_6^{ij})=\f{e^{2X}(\psi_a^{\kappa_0}+\dot\psi)}{2(\g-1)A_0(b,V)}(2A_4(b,V),
A_5^i(b,V), 2A_6^{ij}(b,V)),\qquad i, j=1,2,3
\endcases\tag 3.40$$
and
$$\Cal B_{11}=1+\f{1}{b^2}\sum_{i=1}^{3}(Z_i b)^2,
\qquad \Cal
B_{12}^i(\dot\psi)=-\f{1}{b^2}(\psi_a^{\kappa_0}+\dot\psi) Z_i b
\quad\text{for}\quad i=1,2,3 \tag 3.41$$ and
$$\cases
\ds \Cal B_{20}=-\f{1}{b_0 a_1}(\Bbb
H(\cdot)-\rho_0)\biggl|_{\psi=\psi_a^{\kappa_0}+\dot\psi} +\Cal
D_0\int_0^1\p_{\p_T\psi}\Bbb H(b,\psi_a^{\kappa_0}
+\theta\dot\psi,\na(\psi_a^{\kappa_0}+\theta\dot\psi))d\theta,\\
\ds \Cal
B_{21}=-\f{(T\p_T\psi_a^{\kappa_0}+\psi_a^{\kappa_0})}{b_0}(\Bbb
H(\cdot)-\rho_0) \biggl|_{\psi=\psi_a^{\kappa_0}+\dot\psi}+\Cal D_0
\int_0^1\p_{\p_R\psi}\Bbb
H(b,\psi_a^{\kappa_0}+\theta\dot\psi,\na(\psi_a^{\kappa_0}
+\theta\dot\psi))d\theta,\\
\ds \Cal B_{22}^i=\f{\psi_a^{\kappa_0}}{(b+\psi_a^{\kappa_0})^2}(Z_i
a_0(b,\cdot)+Z_i
a_0(b,\psi_a^{\kappa_0}))\biggl|_{\psi=\psi_a^{\kappa_0}+\dot\psi}+\Cal
D_0\int_0^1\p_{Z_i\psi}\Bbb
H(b,\psi_a^{\kappa_0}+\theta\dot\psi,\na(\psi_a^{\kappa_0}+\theta\dot\psi))d\theta
,\\
\endcases\tag 3.42$$
where $V=(\psi_a^{\kappa_0}+\dot\psi,\na(\psi_a^{\kappa_0}+\dot\psi))$,
and
$\Cal D_0=\biggl(\psi-\ds\f{1}{b_0
a_1}(T a_0+a_0)+\ds\f{\psi}{(b+\psi)^2}\sum_{i=1}^{3}(Z_i
a_0)^2\biggr)\biggl|_{\psi=\psi_a^{\kappa_0}}$.

With respect to the boundary condition on $R=2$ in problem (3.38),
we have\medskip

{\bf Lemma 3.7.} {\it If $\|\dot\psi\|_{C^1([0, T^*]\times [1,
2]\times\Bbb S^2)}\leq C\ve$, the boundary condition on $R=2$ in
$(3.38)$ satisfies the Local Stability Condition which is defined in
\rom{[24]}, that is, there exists a constant $\delta_0>0$
such that\medskip

$1)$ $\Cal B=(\Cal B_{20},\Cal B_{21}, \Cal B_{22}^{1}, \Cal
B_{22}^2, \Cal B_{22}^3)(\dot\psi,\na\dot\psi)$ is transversal to
the boundary $R=2$, namely, $|\Cal B_{21}|>\delta_0>0$.\medskip

$2)$ Denote $\Cal N=(\Cal A_2, \Cal A_4, \Cal A_{5}^1, \Cal A_5^2,
\Cal A_5^3)(\dot\psi,\na\dot\psi)$, then $\ds\t\Cal B=\f{1}{\Cal B_{21}}\Cal
B+\f{1}{|\Cal A_4|}\Cal N$ is a positive time-like direction,
namely,
\medskip

\quad $(i)$ $\ds\f{1}{\Cal B_{21}}\Cal B_{20}+\f{1}{|\Cal A_4|}\Cal
A_2>\delta_0>0$,\medskip

\quad $(ii)$ $\ds-\f{1}{\Cal A_4}\t\Cal B\ M\ \t\Cal
B^{T}>\delta_0>0$ with
$$M=\left(\matrix
\Cal A_1 & \Cal A_2 & \Cal A_3^1 & \Cal A_3^2 & \Cal A_3^3\\
\Cal A_2 & \Cal A_4 & \Cal A_5^1 & \Cal A_5^2 & \Cal A_5^3\\
\Cal A_3^1 & \Cal A_5^1 & \Cal A_6^{11} & \Cal A_6^{12} & \Cal
A_6^{13}\\
\Cal A_3^2 & \Cal A_5^2 & \Cal A_6^{21} & \Cal A_6^{22} & \Cal
A_6^{23}\\
\Cal A_3^3 & \Cal A_5^3 & \Cal A_6^{31} & \Cal A_6^{32} & \Cal
A_6^{33}
\endmatrix\right).$$

\quad $(iii)$
$$\ds\sum_{i=1}^{3}\left(|\Cal B_{22}^i|+|\Cal A_5^i|\right)\leq
C\ve,\quad |\Cal A_2|\leq \f{C}{b_0}(s_0-b_0)^2.$$}

{\bf Remark 3.2.}  {\it In terms of the result in \rom{[24]}, if the
Local Stability Condition on $R=2$ holds, then we can look for a
suitable multiplier to derive the energy estimates near the boundary
$R=2$. The details can be found in \rom{(3.66)-(3.68)}
below.}\medskip

{\bf Proof.} It suffices to verify that 1) and 2) hold when the variables $(b,
\psi^k, \dot\psi)$ in the corresponding coefficients of the second boundary
condition in $(3.38)$ are replaced by $(b_0, \hat\psi, 0)$
respectively.

By Lemma A.6 in Appendix A, we have
$$|\Cal B_{21}|=\bigl(\f{\g-1}{2A\g}\bigr)^{\f{1}{\g-1}}b_0^{\f{2}{\g-1}}(1+\a)>0,
\tag 3.43$$
which means that 1) is shown.

On the other hand, by Lemma A.1-A.3 in Appendix A
and Lemma 3.4, one has
$$\cases
&\Cal A_1=\ds\f{2(s_0-b_0)^2}{(\g-1)b_0^2}(1+\a)>0,\\
&\Cal A_2=\ds\f{\psi}{(\g-1)A_0(b,\psi,\na\psi)}(b_0 a_1
a_2-a_0)(1-(R-1)a_1\p_R\psi)\biggl|_{b=b_0,\psi=\hat\psi}\\
&\qquad\le -\ds\f{C}{b_0}(s_0-b_0)^2<0,\\
&\Cal A_2=0\qquad\qquad\qquad\qquad\qquad \quad \quad\text{on}\quad R=1,\\
&\Cal A_3^i=0, \qquad\qquad\qquad\qquad\qquad\qquad i=1,2,3,\\
&\Cal A_4=-(1+\a)<0,\\
&\Cal A_5^i=0, \qquad\qquad\qquad\qquad\qquad\qquad i=1,2,3,\\
&\Cal A_6^{ij}=\ds-\f{1}{b_0^2}\delta_{ij}, \qquad\qquad\qquad\qquad i,j=1,2,3.
\endcases\tag 3.44$$

Combining (3.44) with Lemma A.6 yields
$$\ds\f{1}{\Cal B_{21}}\Cal B_{20}+\f{\Cal A_2}{|\Cal A_4|}=\f{1}{b_0}(1+\a)>0\tag
3.45
$$
and
$$\align
&\ds-\f{1}{\Cal A_4}\t\Cal B M\t\Cal B^{-T}\\
&=-\f{\Cal A_1}{\Cal A_4}\biggl(-\f{\Cal A_4}{\Cal B_{21}}\Cal
B_{20}+\Cal A_2\biggr)^2(1+\a)\\
&=\f{\g-1}{2}(s_0-b_0)^2(1+\a)>0.\tag 3.46
\endalign$$

Collecting (3.43) and (3.45)-(3.46) yields Lemma 3.7. 1) and
(i)-(ii) in 2) if we select $\delta_0=\ds\f{\g-1}{4}(s_0-b_0)^2$.
Additionally, (iii) in 2) can be derived directly from (3.42) and
(3.44).\qquad\qquad\qquad\qquad\qquad \qquad\qquad
\qquad\qed\medskip

With respect to the boundary condition on $R=1$ in (3.38), we
have\medskip

{\bf Lemma 3.8.}  {\it The boundary condition on $R=1$ in $(3.38)$
is of inhomogeneous Neumann-type, that is,
$$\Cal A_2(\dot\psi)=0,\quad -\Cal
A_{4}(\dot\psi)=\Cal B_{11}(\dot\psi),\quad -\Cal
A_5^i(\dot\psi)=\Cal B_{12}^{i}(\dot\psi)\quad \text{on}\quad
R=1.$$}\medskip

{\bf Proof.} This can be verified directly by Lemma A.1-Lemma A.3 in Appendix A
and the expressions in (3.40)-(3.41), then we
omit the details here.\qquad\qquad\qquad\qquad\qquad\qquad\qquad\qquad\qquad\qquad
\qquad\qquad \qquad \qquad\qed\medskip

For $\|\dot\psi\|_{C^{1}([0, T^*]\times [1, 2]\times\Bbb S^2)}\leq
C\ve$, it follows from (3.44) that the second order quasilinear
equation in  (3.38) is strictly hyperbolic with respect to the
variable $X$ when $X\leq X_0$. Meanwhile, by Lemma 3.7-Lemma 3.8, we
know that the first boundary condition in (3.38) is inhomogeneous
Neumann-type and the second one admits the Local Stability
Condition. Since the Neumann boundary condition does not satisfy the
uniform Lopatinski condition and is very sensitive to the
perturbation in deriving the well-posedness of solution, we have to
choose a suitable iteration  scheme to solve the linearized problem
of (3.38) and further establish the solvability of nonlinear problem
(3.38).

Defining
$\dot\psi_{-1}=\dot\psi_0\equiv\hat\psi_{\si}-\psi_a^{\kappa_0}$, where $\hat\psi_{\si}$
is given Lemma 3.3.
We assume that
$\dot\psi_{n-1}(n\in\Bbb N)$ has been constructed. Set $\t\psi_{n-1}=\dot\psi_{n-1}+\psi_a^{\kappa_0}$.
Motivated by the expressions in (3.8), we
define
$$\cases
\ds a_{0,n}=b+(R-1)\t\psi_{n-1},\\
\ds a_{1,n}=\f{1}{\t\psi_{n-1}+(R-1)\p_R\t\psi_{n-1}},\\
\ds a_{2,n}=\t\psi_{n-1}+(R-2)\p_R\t\psi_{n-1},\\
\ds a_{3,n}=\t\psi_{n-1}\p_T\t\psi_{n-1}+\t\psi_{n-1}\p_T
b+(R-2)\p_T
b\cdot \p_R\t\psi_{n-1},\\
\ds a_{4,n}^i=\t\psi_{n-2} Z_i\t\psi_{n-1}+\t\psi_{n-1} Z_i
b+(R-2)Z_i b\cdot\p_R\t\psi_{n-1},\qquad i=1, 2, 3,
\endcases\tag 3.47$$
here we emphasize that the appearance of the term $\t\psi_{n-2} Z_i\t\psi_{n-1}$ (other than
$\t\psi_{n-1} Z_i\t\psi_{n-1}$)
in the expression of $a_{4,n}^i$ is due to the requirement of
Neumann type boundary condition on $R=1$ in the iteration process of solving (3.38) (one can see
the concrete explanations in Remark 3.3 below).

Let $\dot\psi_{n}$ be determined by the
following problem
$$\cases
\Cal L_n(\dot\psi_n)\equiv\Cal A_{1,n}\p_X^2\dot\psi_n+2 \Cal
A_{2,n}\p_{XR}^2\dot\psi_n+2\ds\sum_{i=1}^{3}\Cal A_{3,n}^{i}\p_X
Z_i\dot\psi_n+\Cal A_{4,n}\p_R^2\dot\psi_n+2\ds\sum_{i=1}^{3}
\Cal A_{5,n}^i\p_R Z_i\dot\psi_n\\
\qquad+\ds\sum_{i,j=1}^{3}\Cal A_{6,n}^{ij}Z_i Z_j\dot\psi_n=\dot
f(e^{({\kappa_0}+1)X}F^{\kappa_0},\dot\psi_{n-1},\na\dot\psi_{n-1}),\quad
(X,R,\o)\in (-\infty,X_0]\times (1,2)\times
\Bbb S^2,\\
\Cal B_1^n (\dot\psi_n)\equiv\Cal
B_{11}^n\p_R\dot\psi_n+\ds\sum_{i=1}^{3}\Cal B_{12}^{i,n}
Z_i\dot\psi_n=\dot
g_1(e^{({\kappa_0}+1)X}G_1^{\kappa_0},\dot\psi_{n-1})
\qquad\qquad \qquad \quad \text{on}\quad R=1,\\
\Cal B_2^n(\dot\psi_n)\equiv\Cal B_{20}^n\p_X\dot\psi_n+\Cal
B_{21}^n\p_R\dot\psi_n+ \ds\sum_{i=1}^{3}\Cal B_{22}^{i,n}
Z_i\dot\psi_n=\dot
g_2(e^{({\kappa_0}+1)X}G_2^{\kappa_0},\dot\psi_{n-1},\na\dot\psi_{n-1})\\
\qquad\qquad\qquad\qquad\qquad\qquad\qquad\qquad\qquad\qquad\qquad\qquad\qquad
\qquad\qquad\qquad\quad \text{on}\quad R=2,
\endcases\tag 3.48$$
where the coefficients in the operator $\Cal L_n(\dot\psi_n)$ and the boundary operators
$\Cal B_i^n(\dot\psi_n)$$(i=1,2)$ are given in terms of
the corresponding coefficients in (3.38), whose arguments $(\dot\psi, a_0, a_i,
a_4^{i})$ $(i=1,2,3)$ are  replaced by $(\dot\psi_{n-1}, a_{0,n}, a_{i,n},
a_{4,n}^i)$ $(i=1,2,3)$ respectively.

{\bf Remark 3.3.}  {\it For  the expressions in $(3.47)$, it follows from a
direct computation that the
boundary condition on $R=1$ in $(3.48)$ is of Neumann type, namely, we have on $R=1$
$$\Cal A_{2, n}=0,\quad -\Cal A_{4, n}=\Cal B_{11}^{n},\quad
-\Cal A_{5, n}^{i}=\Cal B_{12}^{i, n},\qquad i=1,2,3.$$}\medskip

 {\bf $\S
3.3.2.$ Solvability and energy estimates of problem (3.48)}\medskip

With some modifications on the notations in [24], we will use some
weighted Sobolev spaces in this section. For any smooth function
$u(X,R,\o)$ which vanishes at $X=-\infty$, we define the following
norms of $u(X,Y)$ in the domain $\{(X,Y): X\in (-\infty,a], Y\in
(1,2)\times \Bbb S^2\}$ for the constants $\la\in\Bbb N\cup\{0\}$,
$\eta>0$ and $a\in\Bbb R$:
$$\align
&|u|_{\lambda,\eta,X}^2=\ds\sum_{\tau_0+|\tau|=\lambda}\iint_{1\leq
|Y|\leq 2}e^{-2\eta X}\eta^{2\tau_0}|\na_{X,Y}^{\tau}u(X,Y)|^2 dY,\tag 3.49\\
&\|u\|_{\la,\eta, a}^2=\sum_{\tau_0+|\tau|=\la}\int_{-\infty}^{a}\ds\iint_{1\leq
|Y|\leq 2}e^{-2\eta X}\eta^{2\tau_0}|\na_{X,Y}^{\tau}u(X,Y)|^2
dY dX.\tag 3.50
\endalign$$

On the boundaries $R=i$ $(i=1,2)$, we define the boundary norms as follows
$$\ds\langle
u\rangle_{\la,\eta, a}^2=\int_{-\infty}^{a}
\sum_{\tau_0+\tau_1+\tau_2=\la}e^{-2\eta
X}\eta^{2\tau_0}\|\p_X^{\tau_1}u(X,\cdot)\|^2_{H^{\tau_2+1/2}(\p
B_1(0))} dX\tag 3.51$$ and
$$ \ll
u\gg_{\la,\eta,i,a}^2=\ds\int_{-\infty}^{a}\sum_{\tau_0+\tau_1+\tau_2=\la}e^{-2\eta
X}\eta^{2\tau_0}\|\p_{X}^{\tau_1}u(X,\cdot)\|_{H^{\tau_2}(\p
B_i(0))}^2  dX, \quad\ i=1,2.\tag 3.52$$

Based on the notations given in (3.49)-(3.52), we define the
weighted Sobolev spaces $H_{\lambda+1,a}^{\eta}$ and $\Cal
H_{\lambda+1,a}^{\eta}$ in the domain $(-\infty,a]\times (1,2)\times
\Bbb S^2$ as
$$\align
H_{\lambda+1,a}^{\eta}=&\biggl\{u\in
H^{\lambda+1}:|||u|||_{\la+1,\eta,a}
=\sup\limits_{-\infty<X<a}|u|_{\la+1,\eta,X}
+\eta\|u\|_{\la+1,\eta,a}\\
&\qquad +\ds\sum_{i=1}^{2}\bigl(\ll u\gg_{\la+1,\eta,i,a}+\ll \p_R
u\gg_{\la, \eta,i,a}\bigr)<\infty\biggr\}\endalign$$ and
$$\align
\Cal H_{\lambda+1,a}^{\eta}=&\biggl\{u\in
H^{\lambda+1}:\overline{|||}u\overline{|||}_{\la+1,\eta, a}
=\sup\limits_{-\infty<X<a}|u|_{\la+1,\eta,X}
+\eta\|u\|_{\la+1,\eta,a}\\
&\qquad +\ll u\gg_{\la+1,\eta,2,a}+\ll \p_R u\gg_{\la,
\eta,2,a}<\infty\biggr\}.\endalign$$\medskip

To establish the solvability and energy estimates of problem (3.48), at first,
we consider the following initial-boundary
problem by some ideas in [14] and [28]:
$$\cases
L(u)\equiv e_1\p_t^2 u+2e_2\p_{t r}^2
u+2\ds\sum_{i=1}^{3}e_3^i\p_{t}Z_i u-e_4\p_r^2
u-2\ds\sum_{i=1}^{3}e_5^i\p_{r}Z_i u
-\ds\sum_{i,j=1}^{3}e_6^{ij}Z_i Z_j u\\
\qquad =f(t,r,\o)\qquad\qquad\qquad\qquad\qquad\qquad\quad \text{in}
\quad\frak D_0=(-\infty, 0]\times (1,2)\times\Bbb S^2,\\
B_i(u)=d_{i1}\p_r u+\ds\sum_{j=1}^{3}d_{i2}^jZ_j u+d_{i3}\p_t
u=g_i\qquad\quad \text{on}\quad \frak B_i=(-\infty,
0]\times\{i\}\times \Bbb S^2,\quad i=1,2,
\endcases\tag 3.53$$
where we have the following assumptions\medskip

{\bf (A$_1$)} {\it The operator $L$ is strictly hyperbolic with
respect to the time $t$, and fulfill that there exists two positive
constants $\lambda_1<\lambda_2$ such that for any $ (t,x)\in\frak
D_0$ and $\xi\in\Bbb R^4$,
$$\lambda_1\leq e_1(t,x)\leq \lambda_2,\quad \lambda_1|\xi|^2\leq
e_4\xi_0^2+2\ds\sum_{i=1}^{3}
e_5^i\xi_0\xi_i+\ds\sum_{i,j=1}^{3}e_6^{ij}\xi_i\xi_j\leq\lambda_2|\xi|^2.
$$}

{\bf (A$_2$)} {\it The boundary condition on $\frak B_1$ is of
Neumann type, namely,
$$|d_{11}+e_4|+\ds\sum_{i=1}^{3}|d_{12}^i+e_5^i|+|d_{13}-e_2|=0\quad \text{on}
\quad \frak B_1.$$}

{\bf (A$_3$)} {\it The boundary condition on $\frak B_2$ satisfies
the Local Stability Condition. Moreover, there exists a positive
constant $\lambda_0$ which is smaller than $\lambda_1$, such that
$$|d_{21}-e_4|+\ds\sum_{i=1}^{3}|d_{22}^i-e_5^i|\leq \lambda_0,
\quad d_{23}+e_2\geq\lambda_1\quad \text{on}\quad \frak B_2.$$}

{\bf (A$_4$)} {\it For any integer $\lambda>5$,
$$\align
&a(\lambda)\equiv \ds\sum_{i=1,2,4}\|e_i\|_{\lambda, \eta, 0}
+\ds\sum_{l=3,5}\ds\sum_{i=1}^{3}\|e_l^i\|_{\la, \eta, 0}
+\ds\sum_{i,j=1}^{3}\|e_6^{ij}\|_{\lambda, \eta,
0}+\ds\sum_{i=1,3}\ds \langle d_{1i}\rangle_{\lambda,\eta,0}+
\ds\sum_{i=1}^{3}\langle d_{12}^{i}\rangle_{\lambda, \eta, 0}\\
&\qquad\qquad+\ds\sum_{i=1,3}\ds \ll d_{2i}\gg_{\lambda,\eta, 2,0}+
\ds\sum_{i=1}^{3}\ll d_{22}^i\gg_{\lambda, \eta , 2, 0}+1< \infty.
\endalign$$}\medskip

{\bf Proposition 3.9.}  {\it Under the assumptions
\rom{(A$_1$)-(A$_4$)}, there exists a positive constant $\eta_0$,
for $\eta\geq \eta_0$, the integer $\lambda>5$, and
$$||f||_{\lambda,\eta,0}
+\langle g_1\rangle_{\lambda,\eta,0}+\ll g_2\gg_{\lambda,\eta, 2,0}<
+\infty,$$ then the problem $(3.53)$ has a unique solution $u\in
\Cal H_{\lambda+1,0}^{\eta}$ and satisfies for $0\leq
\lambda'\leq \lambda$
$$\overline{|||}u\overline{|||}_{\la'+1,\eta, 0}
\leq C\left(\f{1}{\eta}\|f\|_{\la',\eta,0} +\langle
g_1\rangle_{\la',\eta,0} +\ll g_2\gg_{\la',\eta,2,0}\right).\tag
3.54$$}

To establish Proposition 3.9, we now give some necessary
preparations.

It is easy to know that there exist constants $b_l\ (1\leq
l\leq\la)$ such that $\ds\sum_{l=1}^{\la}b_l(-l)^{k}=1$ for $0\leq
k\leq \la-1$. For $v(t,\cdot)\in H^{\la}(-\infty,
0]$, we define $\breve v(t,\cdot)$ as
$$\breve v(t,\cdot)\equiv\left\{\matrix v(t,\cdot) \quad\text{for $t<0$},\\
\ds\sum_{l=1}^{\la}b_l v(-lt,\cdot) \quad\text{for $t>0$}.\endmatrix\right.\tag
3.55$$ A direct computation yields
$$\|\breve v(t,\cdot)\|_{H^{\la}(\Bbb R)}\leq C\|v(t,\cdot)\|_{H^{\la}(-\infty, 0]}.
\tag 3.56$$

In terms of the definition (3.55), we define the resulting operators
$\breve L$ and $\breve B_l$ $(l=1,2)$ from (3.53) as
$$\cases
\breve L\equiv \breve e_1\p_t^2+2\breve
e_2\p_{tr}^2+2\ds\sum_{i=1}^{3}\breve e_3^i\p_{t}Z_i -\breve
e_4\p_r^2-2\ds\sum_{i=1}^{3}\breve e_5\p_r Z_i
-\ds\sum_{i,j=1}^{3}\breve e_6^{ij} Z_i Z_j,\\
\breve B_i\equiv\breve d_{i1}\p_r+\ds\sum_{j=1}^{3}\breve d_{i2}^j
Z_j+\breve d_{i3}\p_t, \qquad \ i=1, 2,
\endcases\tag 3.57$$
moreover, with (3.56), the corresponding assumptions (A$_1$)-(A$_4$)
still hold.

A truncated function $\chi_1(t)\in C^{\infty}(\Bbb R)$ with
$0\leq \chi_1(t)\leq 1$ is defined as
$$\chi_1(t)=\left\{\matrix 1\quad \text{for $t<\f{1}{4}$},\\
0\quad \text{for $t>\f{3}{4}$}.\endmatrix\right.\tag 3.58$$

In order to prove Proposition 3.9,  we now study the following modified problem of (3.53)
for $l\geq 0$
$$\cases
\t L(\t u)\equiv\bigl(\t e_1\p_t^2 +2\t e_2\p_{t r}^2
+2\ds\sum_{i=1}^{3}\t e_3^i\p_{t}Z_i -\t e_4\p_r^2
-2\ds\sum_{i=1}^{3}\t e_5^i\p_{r}Z_i
-\ds\sum_{i,j=1}^{3}\t e_6^{ij}Z_i Z_j\bigr)(\t u)\\
\qquad =\chi_1(t)\breve f
\qquad\qquad \qquad\qquad \quad\qquad\qquad \text{in}\quad \t \frak D_0=\Bbb R\times (1,2)
\times \Bbb S^2,\\
\t B_1^{l}(\t u)\equiv\bigl(\t d_{11}\p_r+\ds\sum_{i=1}^{3}\t
d_{12}^iZ_i+\t d_{13}\p_t+l\p_t\bigr)\t u= \chi_1(t)\breve g_1
\quad\text{on}\quad \t\frak B_1=\Bbb R\times \{1\}\times \Bbb S^2,\\
\t B_2(\t u)\equiv\bigl(\t d_{21}\p_r+\ds\sum_{i=1}^{3}\t
d_{22}^iZ_i+\t d_{23}\p_t\bigr)\t u= \chi_1(t)\breve g_2
\qquad \quad\text{on}\quad \t\frak B_2=\Bbb R\times \{2\}\times \Bbb S^2\\
\endcases\tag 3.59$$
with
$$\cases
\t L=\chi_1(t)\breve L+(1-\chi_1(t))(\p_t^2-\p_r^2-\ds\sum_{i=1}^{3}Z_i^2),\\
\t B_1^l=\chi_1(t)\breve B_1-(1-\chi_1(t))\p_r+l\p_t,\\
\t B_2=\chi_1(t)\breve B_2+(1-\chi_1(t))\p_r.
\endcases\tag 3.60$$

With respect to problem (3.59), we have

{\bf Lemma 3.10.}  {\it For any fixed $l>0$, if $\t u_l\in
H_{\la+1,\infty}^{\eta}$ is a solution of \rom{(3.59)}, then for any
$T\in\Bbb R$ and $0\leq\la'\leq\la$, we have
$$||| \t u_l|||_{\la'+1,\eta,T}^2\leq C(l,\la')(
\|\chi_1(t)\breve f\|_{\la',\eta,T}^2 +\ds\sum_{i=1}^{2}\ll
\chi_1(t)\breve g_i\gg_{\la',\eta,i, T}^2).\tag 3.61$$}

{\bf Proof.} We will look for a suitable differential operator of
first order $Q=Q_0\p_{t}+Q_1\p_r+\ds\sum_{i=1}^{3}Q_2^i Z_i$ so that
the required norms can be dominated by integrating $2e^{-2\eta t}Q(\t u_l)\t L(\t u_l)$
over the domain
$(-\infty, T]\times (1,2)\times \Bbb S^2$.

Choose a $C^{\infty}$ cut-off function $\Upsilon(r)=\left\{\matrix
1,\quad r<\f{5}{4},\\
\quad\\
0,\quad r>\f{7}{4},
\endmatrix\right.$  and set
$$V_1=\Upsilon(r)\t u_l,\quad\qquad\quad
V_2 =(1-\Upsilon(r))\t u_l.\tag 3.62$$

Then $V_1$ satisfies
$$\cases
\t L(V_1)\ds=F_1\qquad\qquad\qquad\qquad\qquad\qquad\quad \text{in}\quad \t\frak D_0,\\
\t B_{1}^{l}(V_1)=\chi_1(t)\breve g_1\qquad\qquad\qquad\qquad\qquad\quad \text{on}
\quad \t\frak B_1,\\
\text{supp}V_1\subset \{r\leq \f{7}{4}\}
\endcases\tag 3.63$$
with $F_1(t,r,\o)=\Upsilon(r)\chi_1(t)\breve f+[\t L, \Upsilon(r)]\t
u_l.$

On the other hand, $V_2$ satisfies
$$\cases
\t L(V_2)=F_2\qquad\qquad\qquad\qquad\quad \text{in}\quad \t\frak D_0,\\
\t B_2(V_2)=\chi_1(t)\breve g_2\qquad\qquad\qquad\quad \text{on}\quad \t\frak B_2,\\
\text{supp}V_2\subset \{r\geq \f{5}{4}\},
\endcases\tag 3.64$$
with $F_2(t,r,\o)=(1-\Upsilon(r)\chi_1(t)\breve f+[\t L,
1-\Upsilon(r)]\t u_l.$

With two operators $\ds
Q_i=Q_{i0}\p_t+Q_{i1}\p_r+\sum_{j=1}^{3}Q_{i2}^j Z_j$ $(i=1,2)$ to be
determined later on, it follows from (3.63)-(3.64) that
$$\aligned &\int_{-\infty}^{T}\int_{(1,2)\times\Bbb
S^2}2e^{-2\eta t}
Q_i(V_i)\t L(V_i) d xdt\\
=&e^{-2\eta T}\int_{(1,2)\times \Bbb S^2}H_{i0}(T,x)dx
+2\eta\int_{-\infty}^{T}\int_{(1,2)\times\Bbb S^2}e^{-2\eta t}H_{i0}(t,x)dx dt\\
&+(-1)^i \int_{-\infty}^{T}\int_{|x|=i}e^{-2\eta t}H_{i1}dS
dt+\int_{-\infty}^{T}\int_{(1,2)\times\Bbb S^2}e^{-2\eta t}H_{i2}
dxdt
\endaligned\tag 3.65$$
with
$$\align
H_{i0}&=\t e_1 Q_{i0}(\p_t V_i)^2+2\t e_1 Q_{i1}\p_t V_i\p_r
V_i+2\ds\sum_{j=1}^{3}\t e_1
Q_{i2}^j\p_t V_i Z_j V_i\\
&+(2\t e_2 Q_{i1}+\t e_4 Q_{i0})(\p_r V_i)^2+\ds\sum_{j=1}^{3}(2\t
e_2 Q_{i2}^{j}+2\t e_3^j Q_{i1}+2\t e_5^j Q_{i0})\p_r V_i
Z_j V_i\\
&+\ds\sum_{j=1}^{3}\sum_{k=1}^{3}(2\t e_3^j Q_{i2}^k+\t e_6^{jk}
Q_{i0}) Z_j V_i Z_k V_i,\\
H_{i1}&=(2\t e_2 Q_{i0}-\t e_1 Q_{i1})(\p_t V_i)^2-2\t e_4
 Q_{i0}\p_t V_i\p_r V_i\\
&+\ds\sum_{j=1}^{3}(2\t e_2 Q_{i2}^j-2\t e_3^j Q_{i1} -2\t e_5^j
Q_{i0})\p_t V_i Z_j V_i-\t e_4
Q_{i1}(\p_r V_i)^2\\
&-\ds\sum_{j=1}^{3}2\t e_4 Q_{i2}^j\p_r V_i Z_j
V_i-\ds\sum_{j=1}^{3}\sum_{k=1}^{3}(2\t e_5^j
Q_{i2}^k-\t e_6^{jk} Q_{i1})Z_j V_i Z_k V_i,\\
\endalign$$
and $H_{i2}$ stands for the quadratic polynomial of $ \na_{t,x}V_i$
and is independent of $\eta$, whose precise expression is not required.

At first, we treat the case of $i=2$ in (3.63).

By (A$_3$), (3.57) and (3.60), Lemma 3.7 and Remark 3.2, we can
choose a first order operator $Q_2$ such that
$$\cases
H_{20}(t,x)\geq c_0|\na V_2|^2(t,x),\\
H_{21}\geq c_0|\na V_2|^2-c_1|\chi_1(t)\breve g_2|^2\quad
\text{on}\quad \t\frak B_2,
\endcases\tag 3.66$$
where $c_0, c_1$ are some positive constants.

In addition, with the properties of $H_{22}$ and $F_2$ in (3.64), one has
$$
\ds\bigl| H_{22}\bigr|+\bigl|Q_2(V_2)\cdot \t L (V_2)\bigr| \leq
C\bigl(\kappa_0\eta|\na\t u_l| +\f{1}{\kappa_0 \eta}|\chi_1(t)\breve
f|^2+\f{1}{\kappa_0\eta}|\t u_l|^2\bigr),\tag 3.67$$ with
the constant $\kappa_0>0$  being determined later on.

Substituting (3.66)-(3.67) into (3.65) yields
$$\align
&|\na V_2|_{0,\eta,T}^2+\eta||\na V_2||_{0,\eta, T}^2+\ll
\na V_2\gg_{0,\eta,2, T}^2\\
\leq &\ds C\left(\f{1}{\kappa_0\eta}||\chi_1(t)\breve
f||_{0,\eta,T}^2 +\kappa_0\eta||\na\t u_l||_{0,\eta,T}^2
+\f{1}{\kappa_0\eta}||\t u_l||_{0,\eta,T}^2\ds+\ll \chi_1(t)\breve
g_2\gg_{0,\eta,2, T}^2\right).\tag 3.68
\endalign$$

For the case of $i=1$ and $l>0$, we also know that the
boundary condition on $r=1$ in (3.63) also satisfies the Local Stability
Condition, then similar to (3.68), we have
$$\align
&|\na V_1|_{0,\eta,T}^2+\eta||\na V_1||_{0,\eta, T}^2+\ll
\na V_1\gg_{0,\eta,1, T}^2\\
\leq &\ds C(l)\left(\f{1}{\kappa_0\eta}||\chi_1(t)\breve
f||_{0,\eta,T}^2 +\kappa_0\eta||\na\t u_l||_{0,\eta,T}^2
+\f{1}{\kappa_0\eta}||\t u_l||_{0,\eta,T}^2\ds+\ll \chi_1(t)\breve
g_1\gg_{0,\eta,1, T}^2\right)\tag 3.69
\endalign$$

Combining (3.68)-(3.69) with the definition (3.62) yields for $l>0$
$$\align
&|\na \t u_l|_{0,\eta, T}^2+\eta\|\na\t u_l\|_{0,\eta, T}^2
+\ds\sum_{i=1}^{2}\ll \na\t u_l\gg_{0,\eta, i,T}^2\\
\leq &C(l)\biggl(\f{1}{\kappa_0\eta}\|\chi_1(t)\breve
f\|_{0,\eta,1}^2+\kappa_0\eta\|\na\t
u_l\|_{0,\eta,T}^2+\ds\sum_{i=1}^{2}\ll \chi_1(t)\breve
g_i\gg_{0,\eta,i,1}^2\\
&\quad\ds+|\t u_l|_{0,\eta, T}^2+\eta\|\t u_l\|_{0,\eta,
T}^2+\ds\sum_{i=1}^{2}\ll \t u_l\gg_{0,\eta, i, T}^2\biggr).\tag
3.70
\endalign$$

On the other hand, according to Hardy-type inequality, one has
$$\align
&|\t u_l|_{0,\eta, T}^2
=-2\eta\|\t u_l\|_{0,\eta, T}+2\int_{-\infty}^{T}\int_{1\leq
|x|\leq 2}e^{-2\eta t}\p_t\t u_l\cdot\t u_l dxdt\\
&\qquad \leq \f{1}{\kappa_0\eta}\|\t
u_l\|_{0,\eta,T}^2+\kappa_0\eta\|\na\t u_l\|_{0,\eta, T}^2.\tag3.71
\endalign$$

Similarly,
$$\eta\|\t u_l\|_{0,\eta,T}^2+\ll \t
u_l\gg_{0,\eta,2,T}^2\leq C(l)\left(\kappa_0\eta\|\na\t
u_l\|_{0,\eta, T}^2+\f{1}{\eta}\ll \p_t\t
u_l\gg_{0,\eta,2,T}^2\right).\tag 3.72$$

Substituting (3.71)-(3.72) into (3.70) yields for $\kappa_0=\f{1}{2C(l)}$
$$
|||\t u_l|||_{1,\eta,T}^2\leq C(l)\left(\f{1}{\eta}\|\chi_1(t)\breve
f\|_{0,\eta,T}^2+\ds\sum_{i=1}^{2}\ll\chi_1(t)\breve
g_i\gg_{0,\eta,i,T}^2\right),\tag 3.73$$ which means (3.61) holds
for $\la'=0$ with (3.55)-(3.56).

To obtain the higher order energy estimates of  $\t u_l$, we take
the tangential
differential operators $Z_i$ $(i=1,2,3)$ and $\p_t$ on each equality in (3.59)
and then obtain as in (3.73)
$$\align
|||\p\t u_l|||_{1,\eta, T}^2\leq& C(l)\biggl(\f{1}{\eta}||\p
\left(\chi_1(t)\breve
f\right)||_{0,\eta,T}^2+\ds\sum_{i=1}^{2}\ll\p\left(\chi_1(t)\breve
g_i\right)\gg_{0,\eta,i,T}^2\\
&+\f{1}{\eta}||[\t L,\p]\t
u_l||_{0,\eta,T}^2+\ll[\t B_1^l,\p]\t u_l\gg_{0,\eta,1,T}^2
+\ll [\t B_2,\p]\t u_l\gg_{0,\eta,2,T}^2\bigr)\\
\leq& C(l)\left(\f{1}{\eta}\|\chi_1(t)\breve
f\|_{1,\eta,T}^2+\ds\sum_{i=1}^{2}\ll\chi_1(t)\breve
g_i\gg_{1,\eta,i,T}^2\right),\tag3.74
\endalign$$
where $\p$ stands for $Z_i$ $(i=1,2,3)$ or $\p_t$.

By the definition of space $ H_{\lambda+1,T}^{\eta}$, in order to
obtain (3.61) for $\lambda'=1$, it suffices to estimate
$$\sup\limits_{t\in (-\infty, T]}|\p_r^2 \t u_l|_{0,\eta,t}^2
+\eta\|\p_r^2 \t u_l\|_{0,\eta,T}^2,$$ since $\p_R \t u_l$ is a
linear combination of $\p \t u_l$ and $\t B_1^l\t u_l$ or $\t B_2\t
u_l$ on $R=i(i=1,2)$.

It is noted that we have by the
equation in (3.59)
$$\sup\limits_{t\in (-\infty, T]}|\p_r^2 \t u_l|_{0,\eta,t}^2
+\eta\|\p_r^2 \t u_l\|_{0,\eta,T}^2 \leq C(|||\p\t
u_l|||_{1,\eta,T}^2 +\sup\limits_{t\in (-\infty,
T]}|\chi_1(t)\breve f|_{0,\eta,t}^2+\eta||\chi_1(t)\breve
f||_{0,\eta,T}^2).\tag 3.75$$

In addition, one has
$$\sup\limits_{t\in (-\infty, T]}|\chi_1(t)\breve
f|_{0,\eta,t}^2+\eta\|\chi_1(t)\breve f\|_{0,\eta,T}^2\leq
\f{C}{\eta}\|\p_t(\chi_1(t)\breve f)\|_{0,\eta,T}^2.\tag3.76$$

Collecting (3.74)-(3.76) yields (3.61) for $\la'=1$. Analogously, we
can complete the proof of Lemma 3.10 for
$0\leq\la'\leq\la$.\qquad\qquad \qquad\qquad\qquad\qquad\qquad\qquad
\qquad\qquad\qquad\qquad\qquad\qquad\qquad\qquad\qquad
\qquad\qquad\qquad\qed\medskip

In order to solve (3.53), we have to establish a uniform estimate independent of
$l$ in Lemma 3.10.
Since the boundary condition on $\t\frak B_1$ does not satisfy the Local Stability
Condition for
$l=0$, then the estimates in Lemma 3.10 can not be used directly in this case.
To overcome this difficulty, we will apply for some ideas in [28]. \medskip

{\bf Lemma 3.11.} {\it If
$\t u_l\in\Cal H_{\la+1,\infty}^{\eta}$ is a solution of the
problem $(3.59)$, then one has for $0\leq\la'\leq\la$
$$\overline{|||}\t u_l\overline{|||}_{\la'+1,\eta, \infty}^2\leq C\left(
\|\chi_1(t)\breve f\|_{\la',\eta,\infty}^2 +\langle \chi_1(t)\breve
g_1\rangle_{\la',\eta,\infty}+\ll \chi_1(t)\breve
g_2\gg_{\la',\eta,2, \infty}^2 \right),\tag
3.77$$
where $C>0$ is independent of $l$.}

{\bf Proof.} For the problem (3.59) and
$l\geq 0$, we choose $Q_{1}=\p_t$ in the process of deriving (3.65). At this time,
a direct computation
yields
$$
H_{10}(t,x)\geq c_0|\na V_1|^2(t,x),\tag 3.78$$ where $c_0$ is
a positive constant.

Additionally, in terms of the boundary condition on $r=1$ in (3.63) and the
assumption $(A_2)$, we have on $r=1$
$$\align
H_{11}&=2\p_t V_1(\t d_{11}\p_r V_1+\ds\sum_{i=1}^3\t d_{12}^i Z_i V_1
+\t d_{13}\p_t V_1)\\
&=2\p_t V_1(-l\p_t V_1+\chi_1(t)\breve g_1)\\
&\leq 2\p_t V_1\cdot\chi_1(t)\breve g_1.
\endalign$$

Thus,
$$\ds\int_{-\infty}^{T}\int_{|x|=1}e^{-2\eta t}H_{11} dSdt
\leq  C\left(\kappa_0\int_{-\infty}^{T}e^{-2\eta t}\|\p_t
V_1(t,\cdot)\|_{H^{-1/2}(\p B_1(0))}^2 dt+\f{1}{\kappa_0}\langle
\chi_1(t)\breve g_1\rangle_{0,\eta,T}^2\right).\tag 3.79
$$

In addition, similar to (3.67), we have
$$
\bigl| H_{12}\bigr|+\ds\bigl|Q_1(V_1)\cdot \t L(V_1)\bigr| \leq
C\bigl(\kappa_0\eta|\na\t u_l|^2 +\f{1}{\kappa_0
\eta}|\chi_1(t)\breve f|^2+\f{1}{\kappa_0 \eta} |\t
u_l|^2\bigr).\tag 3.80$$

Substituting (3.78)- (3.80) into (3.65) for the case of $i=1$ yields
$$\align
&|\na V_1|_{0,\eta,T}^2+\eta||\na V_1||_{0,\eta, T}^2\\
\leq &\ds C\bigl(\f{1}{\kappa_0\eta}||\chi_1(t)\breve f||_{0,\eta,
T}^2 +\kappa_0\eta||\na\t u_l||_{0,\eta, T}^2
+\f{1}{\kappa_0\eta}||\t u_l||_{0,\eta}^2\\
&+\f{1}{\kappa_0}\langle \chi_1(t)\breve g_1\rangle_{0,\eta,T}^2
+\kappa_0\int_{-\infty}^{T}e^{-2\eta t} ||\p_t
V_1(t,\cdot)||_{H^{-1/2}(\p B_1(0))}^2 dt\bigr).
\endalign$$

Combining this with (3.68) and (3.71) shows
$$\align
&|\na\t u_l|_{0,\eta, \infty}^2+\eta\|\na\t u_l\|_{0,\eta,
\infty}^2+\ll\na\t u_l\gg_{0,\eta,2, \infty}^2\\
\leq &  C\biggl(\f{1}{\kappa_0\eta}\|\chi_1(t)\breve
f\|_{0,\eta,\infty}^2+\kappa_0\eta\|\na\t
u_l\|_{0,\eta,\infty}^2+\f{1}{\kappa_0}\langle\chi_1(t)\breve
g_1\rangle_{0,\eta, \infty}^2+\ll\chi_1(t)\breve
g_2\gg_{0,\eta,2,\infty}^2\\
&\ds+\kappa_0\int_{\Bbb R}e^{-2\eta t}\|\p_t\t
u_l(t,\cdot)\|_{H^{-1/2}(\p B_1(0))}^2 dt\biggr).\tag 3.81
\endalign$$

According to an elementary Proposition 3.17 given in $\S 3.3.4$ below, one has
$$\ds\int_{\Bbb R}e^{-2\eta t}\|\p_t\t u_l\|_{H^{-1/2}(\p
B_1(0))}^2 dt\leq  C\left(\|\chi_1(t)\breve
f\|_{0,\eta,\infty}^2+\|\t u_l\|_{1,\eta,\infty}^2\right).\tag3.82$$

Substituting (3.82) into (3.81) with $\kappa_0=\f{1}{2 C}$ yields
(3.77) for $\la'=0$.

Analogously, the estimates in the case of $\la'>0$ in (3.77) can
be obtained as in Lemma 3.10. Thus, we complete the proof of
Lemma 3.11.\qquad\qquad\qquad\qquad\qquad\qquad\qquad\qquad\qquad\qquad\qquad
\qquad\qquad\qquad
\qquad\qed\medskip

Based on Lemma 3.10 and Lemma 3.11, we now start to  prove the solvability
of the problem (3.59).\medskip

{\bf Lemma 3.12.}  {\it
For any $l\geq 0$, the problem $(3.59)$ has a unique solution $\t
u_l\in\Cal H_{\la+1,\eta,\infty}$ which satisfies the estimate
$(3.77)$.}\medskip

{\bf Proof.} The proof is divided into the following three steps.

{\bf Step 1.} First we consider the case which all the coefficients
mentioned in (A$_4$) are smooth. Moreover, we assume that there
exist two constants $T_2<T_3$ such that
$$\text{supp}f\subset (T_2,T_3)\times (B_2(0)\setminus B_1(0)),
\quad \text{supp}g_i\subset (T_2,T_3)\times \p B_i(0),\quad \ i=1,2.\tag 3.83$$

For $l>0$, we consider the following problem
$$\cases
\t L(\t v_l)=\chi_1(t)\breve f\qquad\qquad \quad\quad \text{in}\quad
(T_2-1,T_3)\times
(1,2)\times \Bbb S^2,\\
\t B_1^{l}(\t v_l)=\chi_1(t)\breve g_1\qquad\qquad \quad
\text{on}\quad (T_2-1,T_3)\times
\{1\}\times\Bbb S^2,\\
\t B_2(\t v_l)=\chi_1(t)\breve g_2\qquad\qquad \quad \text{on}\quad
(T_2-1,T_3)\times\{2\}\times\Bbb S^2,\\
\t v_l(T_2-1,r,\o)=\p_t\t v_l(T_2-1,r,\o)=0.
\endcases\tag 3.84$$

By Theorem 3 in Page 142 of [14], the problem (3.84) has a unique
solution $\t v_l$ which admits
$$\cases
\t v_l(t,r,\o)\in H^{\la+1}((T_2,T_3)\times (B_2(0)\setminus B_1(0))),\\
\t v_l(t,1,\o)\in H^{\la+1}((T_2,T_3)\times \p B_1(0)),\quad \t
v_l(t,2,\o)\in
H^{\la+1}((T_2,T_3)\times \p B_2(0)),\\
\p_r\t v_l(t,1,\o)\in H^{\la}((T_2,T_3)\times \p B_1(0)),\quad
\p_r\t v_l(t,2,\o)\in H^{\la}((T_2,T_3)\times \p B_2(0)).
\endcases$$
Moreover, it follows from (3.83) and the uniqueness of solution to (3.84)
that $\t v_l\equiv 0$ for $t\in (T_2-1, T_2)$.

Let $\t u_l$ be defined as
$$\t u_l=\left\{\matrix \t v_l\quad\text{for}\quad t\in (T_2,T_3),\\
0\quad \text{for}\quad t\leq T_2.\endmatrix\right.$$

Then $\t u_l\in H^{\la+1}((-\infty,T_3)\times (B_2(0)\setminus
B_1(0)))$ is a solution of  (3.59) for $t\in (-\infty,
T_3)$ and
$$\cases
&\t u_l(t,1,\o)\in H^{\la+1}((T_2,T_3)\times \p B_1(0)),\quad
\t u_l(t,2,\o)\in
H^{\la+1}((T_2,T_3)\times \p B_2(0)),\\
&\p_r\t u_l(t,1,\o)\in H^{\la}((T_2,T_3)\times \p B_1(0)),\quad
\p_r\t u_l(t,2,\o)\in
H^{\la}((T_2,T_3)\times \p B_2(0)).
\endcases$$

With the help of Lemma 3.10, we know that $\t u_l$ can be extended
to $t\in\Bbb R$ and satisfies (3.59) as well as the estimate (3.61) for any
$l>0$.\smallskip

{\bf Step 2.} We consider the case that the related coefficients
have the corresponding regularities given by (A$_4$).

In this case, there exist smooth functions $\breve
e_{i,\delta} (i=1, 2, 4), \breve e_{i,\dl}^{j} (i=3, 5; j=1, 2, 3), \breve
e_{6,\dl}^{ij} (i, j=1, 2, 3)$ and $\breve d_{ij}^{\dl} (i=1, 2; j=1, 3),
\breve d_{i2}^{j,\dl} (i=1, 2; j=1, 2, 3),\breve f_{\dl}, \breve
g_{i,\dl} (i=1,2)$ for $\delta>0$ such that
$$\align
&\ds\sum_{i=1,2,4}^{3}\|\breve
e_{i,\delta}-\breve
e_i\|_{\la,\eta,1}+\ds\sum_{i=3,5}\sum_{j=1}^{3}\|\breve
e_{i,\dl}^{j}-\breve e_i^j\|_{\la,\eta,1}+\ds\sum_{i,j=1}^{3}\|\breve
e_{6,\dl}^{ij}-\breve e_6^{ij}\|_{\la,\eta,1}+\ds\sum_{i=1,3}\langle
\breve d_{1i}^{\dl}-\breve d_{1i}\rangle_{\la,\eta,1}\\
&\quad\quad+\ds\sum_{i=1}^{3}\langle \breve d_{12}^{i,\dl}-\breve
d_{12}^i\rangle_{\la,\eta,1}+\ds\sum_{i=1,3}\ll \breve
d_{2i}^{\dl}-\breve d_{2i}\gg_{\la,\eta,2,1}+\ds\sum_{i=1}^{3}\ll
\breve d_{22}^{i,\dl}-\breve d_{22}^i\gg_{\la,\eta,2,1}\\
&\quad\quad+\|\breve f_{\dl}-\breve f\|_{\la,\eta,1}+\langle \breve
g_{1,\dl}-\breve g_1\rangle_{\la,\eta,1}+\ll \breve g_{2,\dl}-\breve
g_2\gg_{\la,\eta,2,1}\to 0\qquad\text{as $\dl\to 0$},\tag 3.85
\endalign$$
moreover $\breve f_{\dl}$ and $\breve g_{i,\dl}$ $(i=1,2)$ have compact support
with respect to the variable $t$.

At this time, there exists a positive constant $\delta_l=\delta(l)$
such that when $0<\delta<\delta_l$, (A$_1$) and (A$_3$)-(A$_4$) are
also satisfied with each function replaced by the corresponding
smooth function and $\la_0, \lambda_1, \lambda_2,a(\lambda)$
replaced by $2\la_0, \lambda_1/2, 2\lambda_2, 2a(\lambda)$
respectively. Meanwhile, (A$_2$) will be replaced by

{\bf (A$_2'$)}\quad $|\breve d_{11}^{\dl}+\breve
e_{4,\dl}|+\ds\sum_{i=1}^{3}|\breve d_{12}^{i,\dl}+\breve
e_{5,\dl}^i|+|\breve d_{13}^{\dl}-\breve e_{2,\dl}|\ll
l\quad\text{for $(t,r,\o)\in\Bbb R\times \{1\}\times\Bbb S^2$}.$

For $l>0$ and $0<\delta<\delta_l$, we consider the following problem
$$\cases
\t L^{\dl}(\t u_{l\dl})=\chi_1(t)\breve f_{\dl}\qquad\qquad\qquad \text{in}
\quad \t\frak D_0,\\
\t B_{1\dl}^l(\t u_{l\dl})=\chi_1(t)\breve g_{1\dl}\qquad\qquad\quad
\text{on}\quad \t\frak
B_1,\\
\t B_{2\dl}(\t u_{l\dl})=\chi_1(t)\breve g_{2\dl}\qquad\qquad\quad
\text{on}\quad \t\frak B_2,
\endcases\tag 3.86$$
where the operators $\t L^{\dl}, \t B_{1\dl}^l, \t B_{2\dl}$ admit the analogous
forms of (3.60), whose coefficients are replaced by the ones in (3.85) respectively.

By Step 1, (A$_2'$) and Lemma 3.10, the problem (3.86) has a unique
solution $\t u_{l\dl}\in H_{\la+1,\infty}^{\eta}$, which satisfies
for any $T\in\Bbb R$,
$$|||\t u_{l\dl}|||_{\la'+1,\eta,T}^2\leq C(l)(\|\chi_1(t)\breve f_{\dl}\|_{\la',\eta, T}^2
+\ds\sum_{i=1}^{2}\ll\chi_1(t)\breve g_{i\dl}\gg_{\la',\eta,i,T}^2).$$

This, together with (3.85), yields that there exist a sequence $\{\dl_{i}\}_{i\in\Bbb
N}\subset (0,\dl_l)$ with $\lim\limits_{i\to\infty}\dl_i=0$ and
a function $\t u_l\in H_{\la+1,\infty}^{\eta}$ such that
$$\t u_{l\dl_i}\rightharpoonup \t u_l\quad \text{in}\quad H^{\lambda+1}(\t\frak D_0)
\cap H^{\lambda+1}(\t\frak B_{1})\cap H^{\lambda+1}(\t\frak
B_{2}),\quad \text{and} \quad \t u_{l\dl_i}\to \t u_l\quad
\text{in}\quad C^2(\Bbb R\times [1,2]\times \Bbb S^2).$$

This shows that $\t u_l$ is a classical solution to the problem
(3.59), which satisfies (3.61) and further admits the uniqueness for any $l>0$.\smallskip

{\bf Step 3.} Since $H_{\la+1,\infty}^{\eta}\subset\Cal
H_{\la+1,\infty}^{\eta}$, then $\t u_l$ also satisfies (3.77) for
any $l>0$. Thus, there exist a positive number sequence $l_i(i\in\Bbb N)$ with
$\lim\limits_{i\to\infty}l_i=0$ and a function $\t u\in \Cal
H_{\la+1,\infty}^{\eta}$ such that
$$\t u_{l_i}\rightharpoonup \t u\quad \text{in}\quad \Cal H_{\la+1,\infty}^{\eta},
\quad \t u_{l_i}\to \t
u\quad \text{in}\quad C^2(\Bbb R\times [1,2]\times \Bbb S^2).$$

Similar to the argument in Step 2, one can derive that $\t u$ is a unique solution
of problem
(3.58) with $l=0$, which satisfies (3.77). Thus,  the proof of Lemma
3.12 is finished.\qquad\qquad\qquad\qquad\qquad\qquad \qed\medskip

In the end of this section, we start to show Proposition 3.9.

{\bf Proof of Proposition 3.9.}  By Lemma
3.11-Lemma 3.12, the problem (3.53) has a solution $u\in\Cal
H_{\la+1,0}^{\eta}$, which satisfies for $0\leq
\la'\leq \la$
$$\align
\overline{|||}u\overline{|||}_{\la'+1, \eta, 0}\leq& C\left(\f{1}{\eta}\|\chi_1(t)\breve
f\|_{\la',\eta,\infty}+\langle \chi_1(t)\breve
g_1\rangle_{\la',\eta,\infty}+\ll \chi_1(t)\breve
g_2\gg_{\la',\eta,\infty}\right)\\
\leq & C\left(\f{1}{\eta}\|f\|_{\la',\eta,0}+\langle
g_1\rangle_{\la',\eta,0}+\ll g_2\gg_{\la',\eta,0}\right).
\endalign$$
Namely, (3.54) is proved. The remainder is to show the
uniqueness of the solution to the problem (3.53). Suppose that (3.53)
has two solutions $u_i\in\Cal H_{\la+1,0}^{\eta} (i=1,2)$, then
$\dot u=u_1-u_2$ satisfies
$$\cases
L(\dot u)=0\quad \text{in}\quad \frak D_0,\\
B_i(\dot u)=0\quad \text{on}\quad \frak B_i,\qquad i=1, 2.
\endcases$$

Utilizing the notations in Lemma 3.10-Lemma 3.11, we can define $\dot
u_1=\Upsilon(r)\dot u$ and $\dot u_2=(1-\Upsilon(r))\dot u$ and
replace $V_i$ $(i=1,2)$ by $\dot u_i$ in the expressions of $H_{ij}$
$(j=0,1,2)$. Under the assumptions (A$_1$)-(A$_4$), it follows from
Lemma 3.10 that for $\la'=4$ and $T=0$
$$\overline{|||}\dot u\overline{|||}_{5,\eta,0}=0,$$
which implies $\dot u=0$. Thus the proof of Proposition 3.9 is completed.
\qquad\qquad\qquad\qquad\qquad\qquad\qquad \qed\medskip

{\bf $\S 3.3.3.$ Solvability of problem (3.10)-(3.12) and proof of Theorem 3.1}
\medskip

In order to solve the nonlinear problem (3.10)-(3.12), we will use
the Newton's iteration. First, we take the approximate solution
$\psi_a^{\kappa_0}$ (mentioned in  (3.35)) with large $\kappa_0$ as
the starting point of the iteration, and set
$\dot\psi_{-1}=\dot\psi_0=\hat\psi_{\si}-\psi_a^{\kappa_0}$, then we
use the modified Newton's iteration scheme (see (3.48)) to modify
$\dot\psi$ gradually to obtain the precise solution. It is noted
that  $\psi_a^{\kappa_0}$ is an approximate solution with error
$e^{({\kappa_0}+1)X}$ near $X=\infty$, and the factor
$e^{({\kappa_0}+1)X}$ will play a crucial role in canceling  the
singularity appeared in the weight of the norm $\Cal
H_{\la+1,X_0}^{\eta}$. Due to (3.37) and (3.39), for the fixed large
$\kappa_0$, one can select $T^*$ suitably small such that for $X\leq
X_0=\ln T^*$ and $\la>5$, on has in (3.48)
$$
\align
&\f{1}{\eta}\|\dot f(e^{({\kappa_0}+1)X}F^{\kappa_0},0)\|_{\la,\eta,X_0}
+\langle \dot
g_1(e^{({\kappa_0}+1)X}G_1^{\kappa_0},0,0)\rangle_{\la,\eta,X_0}\\
&\qquad\qquad \qquad  +\ll \dot
g_2(e^{({\kappa_0}+1)X}G_2^{\kappa_0},0,0)\gg_{\la, \eta, 2,
X_0}\leq C\ve\tag3.87
\endalign$$
and
$$\overline{|||}\dot\psi_{i}\overline{|||}_{\la+1, \eta, X_0}\leq C\ve,
\qquad  i=-1,0.$$

Suppose that
$\overline{|||}\dot\psi_{n-i}\overline{|||}_{\la+1,\eta,X_0}\leq\ve$
holds for $i=1,2$, then it follows from Proposition 3.9 and the
smallness of $\ve$ that the problem (3.48) has a unique solution
$\dot\psi_n\in\Cal H_{\la+1,X_0}^{\eta}$, which satisfies
$$\overline{|||}\dot\psi_n\overline{|||}_{\la+1, \eta, X_0}\leq C\ve,
\qquad  n\in\Bbb N.\tag 3.88$$

To prove the convergence of $\{\dot\psi_{n}\}_{n\in\Bbb N}$, we take
$$\Delta_n\dot\psi=\dot\psi_{n+1}-\dot\psi_n,$$
which satisfies
$$\cases
\Cal L_{n+1}(\Delta_n\dot\psi)=\dot F_{n+1},  \quad (X,R,\o)\in
(-\infty,
X_0]\times (1,2)\times\Bbb S^2,\\
\Cal B_1^{n+1}(\Delta_n\dot\psi)=\dot G_1^{n+1}  \qquad
\qquad \qquad\qquad \qquad \qquad  \text{on}\quad R=1,\\
\Cal B_2^{n+1}(\Delta_n\dot\psi)=(\Cal B_2^n-\Cal
B_2^{n+1})(\dot\psi_n)+\dot G_2^{n+1}  \qquad \quad \text{on}\quad
R=2
\endcases\tag 3.89$$
with
$$\cases
&\dot F_{n+1}=(\Cal L_n-\Cal L_{n+1})(\dot\psi_n)
+\dot f\left(e^{({\kappa_0}+1)X}F^{\kappa_0},\dot\psi_,\na\dot\psi_n\right)-\dot f\left(e^{({\kappa_0}+1)X}F^{\kappa_0},\dot\psi_{n-1},\na\dot\psi_{n-1}\right),\\
&\dot G_1^{n+1}=(\Cal B_1^n-\Cal
B_1^{n+1})(\dot\psi_n)+\dot g_1\left(e^{({\kappa_0}+1)X}G_1^{\kappa_0},\dot\psi_n\right)-\dot g\left(e^{({\kappa_0}+1)X}G_1^{\kappa_0},\dot\psi_{n-1}\right),\\
&\dot G_2^{n+1}=(\Cal B_2^n-\Cal
B_2^{n+1})(\dot\psi_n)+\dot
g_2\left(e^{({\kappa_0}+1)X}G_2^{\kappa_0},\dot\psi_n,\na\dot\psi_n\right)-\dot
g_2\left(e^{({\kappa_0}+1)X}G_2^{\kappa_0},\dot\psi_{n-1},\na\dot\psi_{n-1}\right).
\endcases$$

By Proposition 3.9 and (3.37),  it follows from (3.89) that
$$\align
\overline{|||}\Delta_n\dot\psi\overline{|||}_{\la, \eta, X_0}\leq
&C\left(\f{1}{\eta}\|\dot F_{n+1}\|_{\la-1,\eta, X_0}+\langle \dot
G_{n+1}^1\rangle_{\la-1,\eta,X_0}+\ll\dot G_{n+1}^2\gg_{\la-1,\eta,
2,X_0}\right)\\
\leq
&C\ve\ds\sum_{i=1}^{2}\overline{|||}\Delta_{n-i}\dot\psi\overline{|||}_{\la, \eta, X_0}\\
\le&\f14\ds\sum_{i=1}^{2}\overline{|||}\Delta_{n-i}\dot\psi\overline{|||}_{\la, \eta, X_0},\tag3.90
\endalign$$for small $\ve$.

Combining (3.90) with (3.88) yields that $\{\dot\psi_n\}_{n\in\Bbb N}$ is
convergent  to a function
$\dot\psi$ in $\Cal H_{\la+1,X_0}^{\eta}$. Thus,
$\dot\psi\in\Cal H_{\la+1,X_0}^{\eta}$ is a solution to (3.38), and $\dot\psi+\psi_a^{\kappa_0}$ is a
solution to (3.10) with (3.11)-(3.12). Consequently, the proof
of Theorem 3.1 is completed.\qquad\qquad\qquad\qquad\qquad\qquad\qquad
\qquad\qquad\qquad\qquad\qquad\qquad\qquad\qquad \qed\medskip

{\bf $\S 3.3.4.$ A proof of an elementary Proposition}\medskip

In this part, we will establish an elementary Proposition used in (3.82).

For the second order equation in the domain $D=\Bbb R\times \Bbb R^3_+=\{(t,x): t\in\Bbb R,
x_1, x_2\in\Bbb R, x_3>0\}$
$$
P(u)\equiv A_0(t,x)D_t^2 u+2\ds\sum_{i=1}^{3}A_i(t,x)D_i D_t
u-\ds\sum_{i,j=1}^{3}A_{ij}(t,x)D_{ij}^2u=w(t,x),\qquad \tag 3.91$$
where $D_t=\ds\f{\p_t}{\sqrt{-1}}$, $D_i=\ds\f{\p_i}{\sqrt{-1}}$,
$A_i, A_{ij}\in C^{2,\al}(\bar D)$, and there exists two positive constants
$\Lambda_1<\Lambda_2$ such that
$$\cases
&\Lambda_1\leq A_0(t,x)\leq \Lambda_2,\quad \ds\sum_{i=0}^{3}\|A_i\|_{C^{2,\al}(\bar D)}
+\ds\sum_{i,j=1}^{3}\|A_{ij}\|_{C^{2,\al}(\bar D)}\leq
\Lambda_2\\
&\Lambda_1 |\varsigma|^2\leq
\ds\sum_{i,j=1}^{3}A_{ij}(t,x) \varsigma_i\varsigma_j\leq \Lambda_2|\varsigma|^2\qquad
\text{for any}\quad \varsigma\in\Bbb R^3.\\
\endcases\tag3.92
$$

Denote by $x'=(x_1,x_2)$, and introduce the following notations  for an integer
$\tau\geq 0$, $s\in\Bbb R$,
$\eta\geq 1$, and functions $u(t,x)$, $v(t,x)$
$$\align
&|u|_{\tau,\eta}^2=\ds\sum_{k+|k'|= \tau}\ds\int_{\Bbb R}e^{-2\eta
t}\eta^{2(\tau-k-|k'|)}\|D_t^k D_x^{k'}u(t,\cdot)\|_{L^{2}(\Bbb R_+^3)}^2\, dt,\\
&\langle u\rangle^2_{s,\eta}=\ds\int_{\Bbb R}e^{-2\eta t}\|
u(t,x',0)\|_{H^{s}(\Bbb R^2)}^2\, dt,\\
&(u,v)_{\eta}=\ds\int_{D}e^{-2\eta t}u\cdot v\,
dx dt,\\
&\langle u,v\rangle_{\eta}=\ds\int_{\p D}e^{-2\eta
t}u(t,x',0)\cdot v(t,x',0)\, dx' dt,
\endalign$$
where $\p D=\{(t, x): (t,x')\in\Bbb R^3, x_3=0\}$.

We now prove such a crucial conclusion. \medskip

{\bf Lemma 3.13.}  {\it Under the assumptions \rom{(3.92)}, if
$$\text{supp}u\subset\{(t,x',x_3):t\in \Bbb R, |x'|< M,
0\leq x_3<M\}$$ for some constant $M>0$, and
$|u|_{1,\eta}<+\infty$, then there exists a positive constant
$C=C(\Lambda_1,\Lambda_2)$ such that
$$<D_t u>_{-1/2,\eta}^2\leq
C\bigl(\eta^{-1}|Pu|_{0,\eta}^2+\eta|u|_{1,\eta}^2\bigr).$$}\medskip

To prove Lemma 3.13, we shall use some notations and ideas in
micro-local analysis.

Set  $\Bbb X=(t,x'), \Bbb Y=(s,y'),
\Xi=(\tau,\xi')$ with $\Bbb X, \Bbb Y, \Xi\in \Bbb R\times \Bbb R^2$
and $\lambda_{\eta}(\Xi)=(|\Xi|^2+\eta^2)^{1/2}$ with $\eta\ge 1$. The symbol class
$S_{\eta}^{\delta}$ is defined as
$$S_{\eta}^{\delta}=\{\theta(\Bbb X,\Xi, \Bbb Y, \eta)\in C^{\infty}(\Bbb R^{10}):
|D_{\Bbb X}^{k_1}D_{\Xi}^{k_2}D_{\Bbb Y}^{k_3}\theta|\leq
C_{k_1k_2k_3}\lambda_{\eta}(\Xi)^{\delta-|k_2|},\quad\forall\,
k_1,k_2, k_3\in (\Bbb N\cup\{0\})^3\}.$$

The corresponding  weighted pseudo-differential operator
$\Theta\in\Psi_{\eta}^{\dl}$ with symbol $\theta\in S_{\eta}^{\dl}$
is defined as
$$\Theta(\Bbb X, D_{\Bbb X}, \Bbb Y,
\eta)u=(2\pi)^{-3}\ds\int_{\Bbb R^3}\int_{\Bbb R^3}e^{i(\Bbb
X\cdot\Xi-\Bbb Y\cdot\Xi)}\theta(\Bbb X, \Xi, \Bbb Y, \eta)
e^{(t-s)\eta} u(\Bbb Y)d\Xi d\Bbb Y,$$ for any $u\in \Cal
S_{\eta}(\Bbb R^3)=\{u\in \Cal{D}'(\Bbb R^3): e^{-t\eta}u\in\Cal
S(\Bbb R^3)\}$.\medskip

{\bf Lemma 3.14.}  {\it Let $\Theta$ be
the weighted pseudo-differential operator with the symbol
$\theta$.

$(i)$ Let $\si>0$, $\theta\in\Cal S_{\eta}^{-\si}$, and put
$$N_{\si}^{0}(\theta)=\sup\limits_{\Bbb X, \Xi, \Bbb
Y,\eta}\max\limits_{|\beta|\leq n+1} |D_{\Xi}^{\beta}\theta(\Bbb X, \Xi, \Bbb
Y, \eta)|\lambda_{\eta}(\Xi)^{\si+|\beta|},$$ then
$$|\Theta u|_{0,\eta}\leq C(\si)N_{\si}^{0}(\theta)|u|_{0,\eta}.$$

$(ii)$ For $\theta\in\Cal S_{\eta}^{0}$ and  $\beta\in (0,1)$, put
$$\align
N_{\beta}(\theta)&=\sup\limits_{\Bbb X, \Xi, \Bbb Y, \Bbb Z,
\eta}\max_{|k|\leq n+1}\biggl\{|D_{\Xi}^{k}\theta(\Bbb X, \Xi,
\Bbb Y,
\eta)|\lambda_{\eta}(\Xi)^{|k|}\\
&\qquad +\ds\f{|D_{\Xi}^{k}\theta(\Bbb X, \Xi, \Bbb Y,
\eta)-D_{\Xi}^{k}\theta(\Bbb X, \Xi, \Bbb Z,
\eta)|\lambda_{\eta}(\Xi)^{|k|}}{|\Bbb Y-\Bbb
Z|^{\beta}}\\
&\qquad +\ds\f{|D_{\Xi}^{k}\theta(\Bbb X, \Xi, \Bbb Y,
\eta)-D_{\Xi}^{k}\theta(\Bbb Z, \Xi, \Bbb Y,
\eta)|\lambda_{\eta}(\Xi)^{|k|}}{|\Bbb X-\Bbb Z|^{\beta}}\biggr\},
\endalign$$
then
$$|\Theta u|_{0,\eta}\leq C(\beta)N_{\beta}(\theta)|u|_{0,\eta}.$$}

{\bf Proof.}  One can see page 59 in [28], we omit the proof here.\qquad\qquad
\qquad\qquad \qquad\qquad \qed\medskip

{\bf Lemma 3.15.}  {\it Assume $a(X)\in C^{2,\al}(\Bbb R^3)$ and the
function $\phi_0(\Xi,\eta)\in\Cal S_{\eta}^{\dl}$ with $0<\dl\leq
1$. Let $\Phi_0$ be the weighted pseudo-differential operator with
symbol $\phi_0$ and put the commutator
$[a,\Phi_0]=a(\Phi_0\cdot)-\Phi_0( a\cdot)$. Then
$$|[a, \Phi_0]u|_{0,\eta}\leq C(\dl,\|a\|_{C^{2,\al}})|u|_{0,\eta}.$$}\medskip

{\bf Proof.} Denote by $t=x_0$ and $\tau=\xi_0$, then we have
$$\align
[a, \Phi_0]u&=\int_{\Bbb R^3}\int_{\Bbb R^3}e^{i(\Bbb
X\cdot\Xi-\Bbb Y\cdot\Xi)}\phi_0(\Xi,\eta)(a(\Bbb X)-
a(\Bbb Y))e^{(t-s)\eta}u(\Bbb Y)d\Xi d\Bbb Y\\
&=\int_{\Bbb R^3}\int_{\Bbb R^3}e^{i(\Bbb X\cdot\Xi-\Bbb
Y\cdot\Xi)}\phi_0(\Xi,\eta)(\Bbb X-\Bbb Y)\cdot(\int_{0}^{1}\na
a(\Bbb Y+\kappa(\Bbb X-\Bbb Y))d\kappa)e^{(t-s)\eta}u(\Bbb
Y)d\Xi d\Bbb Y\\
&=-i\int_{\Bbb R^3}\int_{\Bbb R^3}e^{i(\Bbb X\cdot\Xi-\Bbb
Y\cdot\Xi)}\na_{\Xi}\phi_0(\Xi,\eta)\cdot(\int_{0}^{1}\na
a(\Bbb Y+\kappa(\Bbb X-\Bbb Y))d\kappa)e^{(t-s)\eta}u(\Bbb Y)d\Xi
d\Bbb Y\\
&=\undersetbrace I_1\to{-i\int_{\Bbb R^3}\int_{\Bbb R^3}e^{i(\Bbb
X\cdot\Xi-\Bbb Y\cdot\Xi)}\na_{\Xi}\phi_0(\Xi,\eta)\cdot\na
a(\Bbb
Y)e^{(t-s)\eta}u(\Bbb Y)d\Xi d\Bbb Y}\\
&\qquad \undersetbrace I_2\to{-\int_{\Bbb R^3}\int_{\Bbb R^3}e^{i(\Bbb
X\cdot\Xi-\Bbb Y\cdot\Xi)}\t\phi_0(\Bbb X, \Xi, \Bbb Y, \eta)
e^{(t-s)\eta}u(\Bbb Y)d\Xi d\Bbb Y}\tag 3.93
\endalign$$
with
$$\t\phi_0(\Bbb X, \Xi, \Bbb Y, \eta)=\ds\sum_{l,k=0}^{2}\p_{\xi_l\xi_k}^2
\phi_0(\Xi,\eta)
\biggl(\int_{0}^{1}\int_{0}^{1}\p_{y_l y_k}^2a(\Bbb
Y+\kappa_1\kappa_2(\Bbb X-\Bbb Y))d\kappa_1 d\kappa_2\biggr).\tag
3.94$$

Since $\na_{\Xi}\phi_0(\Xi,\eta)\in \Cal S_{\eta}^{\dl-1}$ and is
independent of the variables $\Bbb X$ and $\Bbb Y$, then by Lemma
3.14 and $a\in C^2$, one can arrive at
$$|I_1|\leq C(\dl)|\na a\cdot u|_{0,\eta}\leq
C(\dl)||a||_{C^1}|u|_{0,\eta}.\tag 3.95$$

Due to $a\in C^{2+\al}$, then by Littlewood-Paley
decomposition theory, we have
$$\int_{0}^{1}\int_{0}^{1}\p_{y_l
y_k}^2a(\Bbb Y+\kappa_1\kappa_2(\Bbb X-\Bbb Y))d\kappa_1
d\kappa_2=\ds\sum_{j=-1}^{\infty}a_{lkj}(\Bbb X, \Bbb Y),\qquad
l,k=0,1,2$$ with
$$a_{lkj}(\Bbb X, \Bbb Y)\in C^{\infty}(\Bbb R^3\times\Bbb
R^3),\quad |a_{lkj}|\leq C(\|\na^2
a\|_{C^{\al}})2^{-j\al}.\tag 3.96$$

Now, with (3.94) and (3.96), $\t\phi_0(\Bbb X, \Xi, \Bbb Y, \eta)$ can
be rewritten as
$$\t\phi_0(\Bbb X, \Xi, \Bbb Y,
\eta)=\ds\sum_{j=-1}^{\infty}\ds\sum_{l,k=0}^{2}\t\phi_0^{lkj}(\Bbb
X, \Xi, \Bbb Y, \eta)\tag 3.97$$ with
$$\t\phi_0^{lkj}(\Bbb X, \Xi, \Bbb Y,
\eta)=\p_{\xi_l\xi_k}^2\phi_0(\Xi,\eta)a_{lkj}(\Bbb X, \Bbb
Y)\in\Cal S_{\eta}^{\dl-2},\qquad j=-1,\cdots,\infty; \quad
l,k=0,1,2.$$

By Lemma 3.14 and (3.96), for  $j=-1,\cdots,\infty;\qquad
l,k=0,1,2$, we have
$$\biggl|\int_{\Bbb R^3}\int_{\Bbb R^3}e^{i(\Bbb X\cdot\Xi-\Bbb
Y\cdot\Xi)}\t\phi_0^{lkj}(\Bbb X, \Xi, \Bbb Y,
\eta)e^{(t-s)\eta}u(\Bbb Y)d\Bbb Y\biggr|_{0,\eta}\leq C(\dl,\|\na^2
a\|_{C^{\al}})2^{-j\al}|u|_{0,\eta}.
$$

Combining this with (3.94) and (3.97) yields
$$|I_2|_{0,\eta}\leq C(\dl,\|\na^2
a\|_{C^{\al}})|u|_{0,\eta}.\tag 3.98$$

Thus, it follows from (3.93), (3.95) and (3.98) that
$$|[a, \Phi_0]u|_{0,\eta}\leq C(\dl,\|
a\|_{C^{2,\al}})|u|_{0,\eta},$$ and then the proof of Lemma 3.15 is
completed.\qquad\qquad\qquad\qquad\qquad\qquad\qquad\qquad
\qquad\qquad\qquad\qquad \qquad\qed\medskip

With respect to the equation (3.91), we can easily establish the following a priori
estimate.\medskip

{\bf Lemma 3.16.}  {\it Under the assumptions in $(3.92)$ and \rom{Lemma
3.13}, if $|u|_{1,\eta}+<u>_{1,\eta}<+\infty$, then
$$\langle D_t u\rangle_{0,\eta}\leq C\bigl(\f{1}{\eta}|P
(u)|_{0,\eta}+\ds\sum_{j=1}^{2}\langle D_j u\rangle_{0,\eta}\bigr).$$
}\medskip

{\bf Proof.}  Since the proof procedure is just only a routine $L^2-$energy
estimate and it can be done by following the proof procedure of Lemma 3.11
step by step, then we omit it here.\qquad\qquad\qquad\qquad\qquad
\qquad\qquad\qquad\qed\medskip

Next we start to show Lemma 3.13.

{\bf Proof of Lemma  3.13.} Define $\phi_1(\Xi,\eta)\in
C^{\infty}(\Bbb R^4)$ with $0\leq \phi_1\leq 1$ and
$$\phi_1(\Xi,\eta)=\left\{\matrix 1\quad \text{for}\quad |\Xi|^2+\eta^2\geq 1,\\
\quad\\
0\quad \text{for}\quad |\Xi|^2+\eta^2\leq \ds\frac{1}{2}.
\endmatrix\right.$$

For a small number $\nu\in (0,1)$ which will be determined later on, we set
$\phi_2(\Xi,\eta)\in C^{\infty}\left(\Bbb R^4\setminus\{0\}\right)$
so that $0\leq\phi_2\leq 1$ and
$$\phi_2(\Xi,\eta)=\left\{\matrix\ds 1\quad \text{for}\quad
\nu(\tau^2+\eta^2)\geq 2|\xi'|^2,\\
\quad\\
\ds 0\quad\text{for}\quad\nu(\tau^2+\eta^2)\leq |\xi'|^2.
\endmatrix\right.\tag 3.99$$

Let $\Phi_1, \Phi_2$ and $\Phi_3$ be the weighted pseudo-differential
operators with symbols $1-\phi_1$, $\phi_1\phi_2$ and $\phi_1(1-\phi_2)$
respectively. A direct computation shows that
$$\Phi_1\in \Psi_{\eta}^{-\infty},\quad \Phi_2\in \Psi_{\eta}^{0},\quad
\Phi_3\in \Psi_{\eta}^{0}.\tag 3.100$$

Choose $\chi_i(x)\in C^{\infty}(\Bbb R^3)(i=1,2)$ so that
$$\chi_1(x)=\left\{\matrix 1\quad\text{for}\quad |x|\leq 2M,\\
\quad\\
0\quad\text{for}\quad |x|\geq \ds\f{9}{4}M,
\endmatrix\right.\quad\quad \chi_2(x)=
\left\{\matrix 1\quad\text{for}\quad |x|\leq \ds\f{5}{2}M,\\
\quad\\
0\quad\text{for}\quad |x|\geq 3M.\endmatrix\right.$$

In this case, by the assumption in
Lemma 3.12, one has
$u=\chi_1 u$ and $u=\ds\sum_{i=1}^{3}\Phi_i u$.

At first, we estimate $<D_t\left(\chi_1\Phi_2 u\right)>_{-1/2,\eta}$.

It follows from (3.91) and a
direct computation that
$$P(\chi_1\Phi_2 u)=w_1\quad in\quad  D\tag 3.101$$
with
$$\align
w_1=&\chi_1\Phi_2(\chi_2 P(u))-\bigl\{(D_3^2\chi_1)\Phi_2
u+2(D_3\chi_1)\Phi_2 (D_3 u)\bigr\}\\
&-\ds\sum_{i=1}^{2}(A_{3i}+A_{i3})\bigl\{(D_{3i}\chi_1)\Phi_2
u+(D_3\chi_1)D_i(\Phi_2 u)+(D_i\chi_1)\Phi_2 (D_3 u)\bigr\}\\
&-\ds\sum_{i,j=1}^{2}A_{ij}\bigl\{(D_i\chi_1)D_j\Phi_2
u+(D_j\chi_1)D_i\Phi_2 u+(D_{ij}\chi_1)\Phi_2 u\bigr\}\\
&+2\ds\sum_{i=1}^{3}A_i (D_i\chi_1)D_t\left(\Phi_2
u\right)+\chi_1\biggl\{\ds\sum_{i=1}^{2}[\Phi_2 D_i,\chi_2
A_{3i}]D_3 u\\
&+\ds\sum_{i,j=1}^{2}[\Phi_2 D_i, \chi_2 A_{ij}]D_j
u-2\ds\sum_{i=1}^{3}[\Phi_2 D_t, \chi_2 A_i]D_i u-[\Phi_2 D_t,
\chi_2 A_0]D_t u\biggr\}\\
&+\chi_1\Phi_2\biggl\{-\ds\sum_{i=1}^{2}(D_i(\chi_2 A_{3l}))D_3
u-\ds\sum_{i,j=1}^{2}(D_i(\chi_2 A_{ij}))D_j u\\
&+2\ds\sum_{i=1}^{3}(D_t(\chi_2 A_i))D_i u+(D_t(\chi_2 A_0))D_t
u\biggr\}.
\endalign$$
Here $[\cdot, \cdot]$ means the commutator.

Then utilizing the results in Lemma 3.14-Lemma 3.15 and (3.100)
yields
$$|w_1|_{0,\eta}\leq C(\Lambda_1,\Lambda_2)(|P(u)|_{0,\eta}+|u|_{1,\eta}).\tag 3.102$$

In addition, by Lemma 3.16 and (3.101)-(3.102), we have
$$
<D_t(\chi_1\Phi_2 u)>_{0,\eta}^2\leq
C(\Lambda_1,\Lambda_2)\bigl(\ds\sum_{i=1}^{2}<D_i(\chi_1\Phi_2
u)>_{0,\eta}^2 +\eta |u|_{1,\eta}^2+\eta^{-1}
|P(u)|_{0,\eta}^2\bigr).\tag 3.103$$

On the other hand, it follows from (3.99) that
$$\ds\sum_{i=1}^{2}<D_i(\chi_1\Phi_2 u)>_{0,\eta}^2\leq C\nu\bigl(
<D_t(\chi\Phi_2 u)>_{0,\eta}^2+\eta^2|u|_{1,\eta}^2\bigr).$$

Substituting this into (3.103) and taking $\ds\nu=\frac{1}{2C
C(\Lambda_1,\Lambda_2)\eta}$ yields
$$<D_t(\chi_1\Phi_2 u)>_{0,\eta}^2\leq C(\Lambda_1,\Lambda_2)\bigl(\eta
|u|_{1,\eta}^2+\eta^{-1} |P(u)|_{0,\eta}^2\bigr).\tag 3.104$$

Next, we estimate $\bigl\langle\chi_1 \langle D'\rangle^{-1/2}\Phi_3
u\bigr\rangle_{0,\eta}$.

Set
$$\phi_4(\Xi,\eta)=(|\Xi|^2+\eta^2)^{1/4}(1+|\xi'|^2)^{-1/4}
\phi_1(\Xi,\eta)(1-\phi_2(\Xi,\eta)),$$ and let $\Phi_4$ be the
weighted pseudo-differential operator with symbol $\phi_4$. By
(3.99), one has
$$\phi_4\in \Cal S_{\eta}^{0},\quad \langle D'\rangle^{-1/2}\Phi_3
u=\Lambda^{-1/2}\Phi_4 u\tag 3.105$$ with $\Lambda^{-1/2}$ be the
weighted pseudo-differential operator with symbol
$(|\Xi|^2+\eta^2)^{-1/4}$.

In addition, by (3.103), one has
$$P(\chi_1\langle D'\rangle^{-1/2}\Phi_3 u)=w_2\quad in\quad
D,$$ where $w_2$ is a function given by replacing $\Phi_2$ by
$\Lambda^{-1/2}\Phi_4$  in the expression of $w_1$.

Applying for Lemma 3.14-Lemma 3.15 with (3.100) and (3.105) to $w_2$
yields
$$|w_2|_{0,\eta}\leq C(\Lambda_1,\Lambda_2)\bigl(|P(u)|_{0,\eta}+|u|_{1,\eta}
\bigr).\tag 3.106$$

It follows from Lemma 3.16 and (3.106) that
$$\langle D_t(\chi_1\langle D'\rangle^{-1/2}\Phi_3
u)\rangle_{0,\eta}^2\leq
C(\Lambda_1,\Lambda_2)(\ds\sum_{i=1}^{2}\langle
D_i(\chi_1\langle D'\rangle^{-1/2}\Phi_3
u)\rangle_{1,\eta}^2+\eta |u|_{1,\eta}^2+\eta^{-1}
|P(u)|_{0,\eta}^2).\tag 3.107$$

By (3.99), one has
$$\ds\sum_{i=1}^{2}\bigl\langle D_i(\chi_1\langle
D'\rangle^{-1/2}\Phi_3 u)\bigr\rangle_{0,\eta}^2\leq C\langle
u\rangle_{1/2,\eta}^2\leq C|u|_{1,\eta}^2.\tag 3.108$$

Due to
$\bigl\|\left[\chi_1, \langle D'\rangle^{-1/2}\right]v\bigr\|_{L^2(\Bbb R^2)}\leq C\|\langle
D'\rangle^{-1/2}v\|_{L^2(\Bbb R^2)}$, then we arrive at
$$\align
\langle D_t(\chi_1\langle D'\rangle^{-1/2}\Phi_3
u)\rangle_{0,\eta}&\geq \langle \langle
D'\rangle^{-1/2}( D_t\left(\chi_1\Phi_3
u\right))\rangle_{0,\eta}-\langle [\chi_1,
\langle D'\rangle^{-1/2}]\left( D_t(\Phi_3 u)\right)\rangle_{0,\eta}\\
&\geq \langle D_t\left(\chi_1\Phi_3
u\right)\rangle_{-1/2,\eta}-C\langle \langle
D'\rangle^{-1/2}\left( D_t(\Phi_3 u)\right)\rangle_{0,\eta}.
\endalign$$

Noticing that $|\xi'|^2\geq\nu(\tau^2+\eta^2)$ holds on
$\text{Supp}\phi_1(1-\phi_2)$, then one obtains for $\ds\nu=\f{1}{2C
C(\Lambda_1,\Lambda_2)\eta}$
$$\langle \langle D'\rangle^{-1/2}\left(D_t\left(\Phi_3
u\right)\right)\rangle_{0,\eta}^2\leq
C(\Lambda_1,\Lambda_2)\eta|u|_{1,\eta}^2.$$

Substituting these estimates into (3.107) derives
$$\langle D_t\left(\chi_1\Phi_3 u\right)\rangle_{-1/2,\eta}^2\leq
C(\Lambda_1,\Lambda_2)\bigl(\eta
|u|_{1,\eta}^2+\eta^{-1}|P(u)|_{0,\eta}^2\bigr).\tag 3.109$$

Because of $\text{supp}(1-\phi_1)\subset\{(\Xi,\eta)\in\Bbb R^4:
|\Xi|^2+\eta^2\leq 1\}$ and (3.100), then
$$\langle D_t \left(\chi_1\Phi_1 u\right)\rangle_{-1/2,\eta}^2\leq \langle
u\rangle_{1/2,\eta}^{2}\leq C|u|_{1,\eta}^2.\tag 3.110$$

It is noted that $D_t u=\ds\sum_{i=1}^{3}D_t(\chi_1\Phi_i u)$, then in terms of (3.104) and
(3.109)-(3.110), we have
$$\langle D_t u\rangle_{-1/2,\eta}^2\leq C(\Lambda_1,\Lambda_2)\bigl(\eta
|u|_{1,\eta}^2+\eta^{-1} |P(u)|_{0,\eta}^2\bigr),$$ which means that
the proof of Lemma 3.13 is
completed.\qquad\qquad\qquad\qquad\qquad\qquad\qquad\qquad\qquad\qquad\qed\medskip

{\bf Proposition 3.17 (An elementary Proposition)} {\it If $\O\subset\Bbb R^3$ is a
smooth bounded
domain with compact boundary,
with respect to the following equation in $\frak D=\Bbb R\times \O$
$$
\t P(u)=\t A_0(t,x)D_t^2 u+2\ds\sum_{i=1}^{3}\t A_i(t,x)D_i D_t u
-\ds\sum_{i,j=1}^{3}\t A_{ij}(t,x)D_{ij}u=\t w(t,x),
$$
where $\t A_i(0\leq i\leq 3)$ and $\t A_{ij}(1\leq i,j\leq 3)$
satisfy the assumption \rom{(3.92)} in $\frak D$. Then we have the
following a priori estimate
$$\align
&\int_{\Bbb R}e^{-2\eta t}\|D_t u(t,\cdot)\|_{H^{-1/2}(\p\O)}^2 dt\\
\leq &C(\Lambda_1,\Lambda_2)\left(\f{1}{\eta} \int_{\Bbb R}e^{-2\eta
t}\|\t P(u)\|_{L^{2}(\O)}^2 dt +\eta\ds\sum_{k+|k'|=1}\int_{\Bbb
R}e^{-2\eta t} \eta^{2(1-k-|k'|)}\|D_t^k D_x^{k'}u\|_{L^2(\O)}^2
dt\right).\tag 3.111\endalign$$}

{\bf Proof.} Since $\O$ is a smooth bounded domain in $\Bbb R^3$
with compact boundary, (3.111) follows from the skills of partition
of unity and local flattening of the boundary of $\O$, and Lemma
3.13, then Proposition 3.17 is
proved.\qquad\qquad\qquad\qquad\qquad\qquad\qquad\qquad
\qquad\qquad\qquad\qquad\qquad\qquad
\qquad\qquad\qquad\qquad\qquad\qquad\qquad\qed

\vskip 0.3 true cm \centerline{\bf $\S 4.$ Another reformulation of (1.6)
with (1.7)-(1.9) and some preliminaries} \vskip 0.4 true cm

In this section, we reformulate the
problem (1.6) with (1.7)-(1.9) in another form, which will be required
to establish the global weighted energy
estimate and further prove the global existence in subsequent sections.

As in [17], we denote certain partial Klainerman's vector fields
by
$$Z=\{Z_j: 0\leq j\leq 2n-3\}, \tag 4.1$$
where $n=2, 3$, $Z_0=t\p_t+r\p_r$, and $ Z_i$ $(1\leq i\leq 2n-3)$
is given in (1.7).

In addition, for notational convenience, for $\nu>0$, $l\in\Bbb N$
and $m\in \Bbb Z$, we define the following  space
$$O_{m}^{l}(\nu)=\{u(t,r,\o)\in C^{l}(\O_+): |\na_{t,r}^{l_0}
Z_1^{l_1}Z_2^{l_2}Z_3^{l_3}u|\leq C_{l}
(1+t)^{m-l_0}\nu \text{ with }|l_0|+\ds\sum_{i=1}^{3}l_i=l\}.\tag
4.2$$

Under the coordinate transformation (1.5), the equation (1.6)
becomes
$$
\p_t^2\Phi^{+}+2\dsize\sum_{k=1}^n\p_k\Phi^{+}\p_{tk}^2\Phi^{+}
+\ds\sum_{i,j=1}^{n}\p_{i}\Phi^{+}\p_j\Phi^{+}\p_{ij}^2\Phi^{+}-c^2(\rho^{+})
\biggl(\p_r^2+\frac{n-1}{r}\p_r+\frac{1}{r^2}\Delta_{\Bbb
S^{n-1}}\biggr) \Phi^{+}=0\quad \text{in}\quad \Omega_+,\tag 4.3
$$
where $\Delta_{\Bbb S^{n-1}}$ is the Laplace-Beltrami operator on
$(n-1)-$dimensional spherical surface. Namely,
$$\Delta_{\Bbb S^1}=Z_1^2,\quad
\Delta_{\Bbb S^{2}}=\ds\sum_{i=1}^{3}Z_i^2.$$

Since the back ground solution $\hat\Phi(t,x)$ given in Remark 2.1
does not satisfy the boundary condition (1.7),  we have to modify it
so that the new resulting background solution $\Phi_a(t,x)$ in
$\O_+$ satisfies (1.7) as well as  other required properties. This
is achieved by choosing a smooth function $f_a=f_a(t,r,\o)$ in
$\O_+$ with $f_a\in O_{-1}^{\infty}(\ve)$ and setting
$\Phi_a(t,r,\o)=(1+f_a(t,r,\o))\hat\Phi(t,x)$, the details can be
found in Appendix B.

Let $(\Phi^+(t,r,\o), \zeta(t,\o))$  be the solution of the problem (4.3) with
(1.7)-(1.10) and $(\dot\vp(t,r,\o), \xi(t,\o))$ be the perturbation of the
modified background solution, that is, $\dot\vp=\Phi^+-\Phi_a,
\ds\xi=\f{\zeta}{t}-s_0$. We now start to reformulate the nonlinear problem (4.3) with
(1.7)-(1.9). For notational convenience, from now on, we neglect all
the superscripts $``+"$.

By a direct computation, (4.3) can be reduced to
$$
\Cal L\dot\vp+\Cal P\dot\vp=R_0(t,x)\quad \text{in}\quad
\Omega_{+},\tag 4.4
$$
where the operators $\Cal L$ and $\Cal P$ have the forms
$$\cases
\ds\Cal
L=\p_{t}^2+2P_1(s)\p_{tr}^2+P_2(s)\p_{r}^2-\f{1}{r^2}P_3(s)\Delta_{\Bbb
S^{n-1}}
+\ds\f{1}{r}P_4(s)\p_t+\f{1}{r}P_5(s)\p_r,\\
\ds\Cal P=\ds\sum_{i,j=1}^{n}f_{ij}(t,x)\p_{ij}^2
+\ds\sum_{i=1}^{n}f_{0i}(t,x)\p_{tx_i}^2
\endcases\tag 4.5$$
with
$$\cases
\ds P_1(s)=\hat u (s),\\
\ds P_2(s)=\hat u^2(s)-c^2(\hat\rho(s)),\\
\ds P_3(s)=c^2(\hat\rho(s)),\\
\ds P_4(s)=(\g-1)\biggl((n-1)\hat u(s)+s\hat u'(s)\biggr),\\
\ds P_5(s)=(n-1)(\g-1)\hat u^2(s)-(n-1)c^2(\hat\rho(s))-2s^2\hat
u'(s)+(\g+1)s \hat u(s)\hat u'(s)
\endcases\tag 4.6$$
and
$$\cases
f_{ij}(t,x)=\biggl(\p_i\hat\Phi\p_j(\hat\Phi
f_a+\dot\vp)+\p_j\hat\Phi\p_i(\hat\Phi f_a
+\dot\vp)+\p_i(\hat\Phi f_a+\dot\vp)\p_j(\hat\Phi+\dot\vp)\biggr)\\
\quad\quad\ds +(\g-1)\delta_{ij}\biggl(\p_t(\dot\vp+\hat\Phi
f_a)+\sum_{k=1}^{n}(\p_k(\Phi_a
+\f{1}{2}\dot\vp)\p_k\dot\vp+\p_k(\hat\Phi+\f{1}{2}
\hat\Phi f_a)\p_k(\hat\Phi f_a))\biggr),\\
\quad\quad\quad\quad\quad\quad\quad\quad\quad\quad\quad\quad\quad\quad\quad\quad
\quad\quad\quad\quad
\quad\quad\quad\quad\quad\quad\quad\quad i,j=1,\cdots, n,\\
f_{0i}(t,x)=-2\p_i(\dot\vp+\hat\Phi f_a),\quad i=1,\cdots, n
\endcases\tag 4.7$$and
$$\cases
R_0(t,x)=-2\ds\sum_{k=1}^{n}\p_{t k}^2(\hat\Phi
f_a)\p_k\dot\vp-\ds\sum_{i=1,j}^{n}\biggl(\p_{ij}^2\hat\Phi(\p_i(\hat\Phi
f_a)\p_j\dot\vp
+\p_j(\hat\Phi f_a)\p_i\dot\vp)\\
\quad\quad\quad\ds+\p_{ij}^2\hat\Phi\p_i\dot\vp\p_j\dot\vp+\p_{ij}^2(\hat\Phi f_a)
(\p_i\Phi\p_j\dot\vp+\p_j\Phi_a\p_i\dot\vp)\biggr)\\
\ds\quad\quad\quad\ds-(\g-1)\Delta\hat\Phi\biggl(\p_t(\hat\Phi
f_a)+\sum_{k=1}^{n}(\p_k(\hat \Phi
f_a+\f{1}{2}\dot\vp)\p_k\dot\vp+\p_k(\hat\Phi+\f{1}{2}\hat\Phi
f_a)\p_k(\hat\Phi f_a))\biggr)\\
\ds\quad\quad\quad\ds-(\g-1)\Delta(\hat\Phi
f_a)\biggl(\p_t(\dot\vp+\hat\Phi
f_a)+\sum_{k=1}^{n}(\p_k(\Phi_a+\f{1}{2}\dot\vp)\p_k\dot\vp+\p_k(\hat\Phi+\f{1}{2}\hat\Phi
f_a)\p_k(\hat\Phi f_a))\biggr)\\
\ds\quad\quad\quad\ds-\biggl(\p_t^2\Phi_a+2\sum_{k=1}^{n}\p_k\Phi_a\p_{tx_k}^2\Phi_a
+\sum_{i,j=1}^{n}\p_i\Phi_a\p_j\Phi_a\p_{ij}^2\Phi_a-c^2(\hat\rho)\Delta\Phi_a\biggr).
\endcases\tag 4.8$$

On the fixed boundary $r=\sigma(t,\o)$, we have
$$\Cal B_{\sigma}\dot\vp=0.\tag 4.9$$

On the free boundary $r=\zeta(t,\omega)$, by the continuity
condition (1.9), one can rewrite Rankine-Hugoniot condition (1.8) as
$$B_1\p_r\dot\vp+B_2\p_t\dot\vp+B_3\xi=\kappa(\xi,\na\dot\vp)+R_1(t,x)\quad
\text{on}\quad r=\zeta(t,\o),\tag 4.10$$ where
$$\cases
B_1&=2\hat\rho(s_0)\hat
u(s_0)-\ds\f{1}{c^2(\hat\rho(s_0))}\hat\rho(s_0)\hat
u(s_0)\biggl(\f{1}{2}\hat
u^2(s_0)-h(\hat\rho)(s_0)+h(\rho_0)\biggr),\\
B_2&=\hat\rho(s_0)-\rho_0-\ds\f{1}{c^2(\hat\rho(s_0))}\hat\rho(s_0)\biggl(\f{1}{2}\hat
u^2(s_0)-h(\hat\rho)(s_0)+h(\rho_0)\biggr),\\
B_3&=2\hat\rho(s_0)\hat u(s_0)\hat
u'(s_0)+\hat\rho'(s_0)\biggl(\f{1}{2}\hat
u^2(s_0)-h(\hat\rho)(s_0)+h(\rho_0)\biggr)\\
&\quad\quad-(\hat\rho(s_0)-\rho_0)\biggl(\hat u(s_0)\hat
u'(s_0)+\ds\f{c^2(\hat\rho)(s_0)}{\hat\rho(s_0)}\hat\rho'(s_0)\biggr),
\endcases
$$
the generic function $\kappa(\xi,\nabla\dot\vp)$ is used to denote
the quantity dominated by $C(b_0)|(\xi,\nabla\dot\vp)|^2$.

By Lemma 4.2  below, we know $B_1\not=0$ in (4.10) for large
$b_0$. Thus, the equation (4.10) can be rewritten as
$$\Cal
B_0\dot\vp+\mu_2\xi=\f{1}{B_1}\kappa(\xi,\nabla\dot{\vp})+\f{1}{B_1}R_1(t,x)
\quad \text{on}\quad r=\chi(t,\o),\tag4.11$$ where
$$\Cal
B_0\dot\vp=\p_r\dot{\vp}+\mu_1\p_t\dot{\vp}\quad\text{with}
\quad\mu_1=\ds\f{B_2}{
B_1}\quad\text{and}\quad\mu_2=\ds\f{B_3}{B_1}.$$

Besides, (1.9)  implies that
$$
\align \dot\varphi(t,\chi(t,\omega),\omega)= &\Phi
(t,\chi(t,\omega),\omega)-{{\Phi}_a} (t,s_0
t,\omega)-\bigl({{\Phi}_a}(t,\chi(t,\omega),\omega)-{{\Phi}_a}
(t,s_0 t,\omega)\bigr)
\\
=&-\bigl({{\Phi}_a}(t,\chi(t,\omega),\omega)-{{\Phi}_a} (t,s_0
t,\omega)\bigr)
\\
=&\mu_3(t,x) t\xi(t,\omega)\tag 4.12
\endalign
$$
with $\ds\mu_3(t,x)=-\int_0^{1}{\hat u}
(s_0+\tau\xi(t,\omega))(1+f_a(t,\chi(t,\o),\o))d\tau<0$

On the other hand, by the local existence and stability result established
in $\S 3$, we only need to solve problem (4.4) in the domain
$\{(t,r,\omega)\colon t\geq 1,\, \si(t,\omega)\leq
r\leq\chi(t,\omega),\, \omega \in \Bbb S^{n-1}\}$ with the boundary
conditions (4.9) and (4.11)-(4.12) and small initial data
$\dot\varphi(t,x)|_{t=1}$, $\p_t\dot\varphi(t,x)|_{t=1}$, and
$\xi(t,\omega)|_{t=1}$. Here, the smallness of initial data means that
$$\sum_{l\le k_0}|\na^l\dot\varphi(1, x)|
+\sum_{l\le k_0}|\na^l \xi(1,\omega)|\le C\ve, \tag 4.13$$ where
$k_0\in\Bbb N$, $k_0\ge 2 n+3$.\medskip

Under the preparations above, Theorem 1.1 is actually equivalent to\medskip

{\bf Theorem 4.1.} {\it For $n=2,3$, if $\sigma(t,\o)=t b(t,\o)$
satisfies the assumptions in \rom{Theorem 1.1}, then the problem $(4.4)$
with $(4.9)$ and $(4.11)$-$(4.13)$ has a unique global $C^{\infty}$
shock solution $\ds(\dot\vp, \xi)$. Moreover,
$(\na\dot\vp,\xi)$ approaches zero as $t$ tends to infinity with
rate $(1+t)^{-m_0}$ for any positive number
$m_0<\ds\f{5}{4}-\ds\f{1}{4}\sqrt{\ds\f{\g+1}{2}}$ if $n=2$  and
$m_0<\ds\f{3}{2}-\ds\f{1}{4}\sqrt{\ds\f{\g+7}{2}}$ if $n=3$
respectively.}
\medskip

For the later uses, we now list some elementary estimates on the
coefficients in (4.6) and (4.10)-(4.11) and the non-homogenous terms
$R_i (i=0,1)$ in (4.4) and (4.11). Since these estimates come from a
direct but tedious computation by Lemma 2.1 in $\S 2$ and Lemma B.1
in Appendix B, then we omit the details here.\medskip

{\bf Lemma 4.2.} {\it For $R_i (i=0,1)$ given in $(4.4)$ and
\rom{(4.11)} respectively, if $\na\dot\varphi\in O_0^{l}(\ve)$ and
$\xi\in O_0^{l}(\ve)$ for any $l\in\Bbb N$, then
$$R_0\in O_{-2}^{l}(\ve^2),\quad R_1\in O_{-1}^{l}(\ve).$$}

With respect to the coefficients of operator $\Cal L$ in (4.5), we
have\medskip

{\bf Lemma 4.3.} {\it If $b_0>0$ is large enough, $1<\gamma<3$ and $b_0\le
s\le s_0$, then
$$
\align
&P_1(s)=b_0(1+\a),\\
&P_2(s)=\f{3-\g}{2}b_0^2(1+\a),\\
&P_3(s)=\f{\g-1}{2}b_0^2(1+\a),\\
&P_4(s)=b_0\biggl(\a\biggr),\\
&P_5(s)=-\f{\g-1}{2}(n-1)b_0^2(1+\a),\\
&P_1'(s)=-(n-1)(1+\a),\\
&P_2'(s)=-2(n-1)b_0(1+\a),\\
&P_3'(s)=b_0\biggl(\a\biggr).
\endalign$$}\medskip

In addition, $B_i (i=1,2,3)$ in (4.10) and $\mu_j (j=1,2)$ in
(4.11) admit the following estimates:\medskip

{\bf Lemma 4.4.}  {\it For large $b_0$, we have
$$
\align &B_1=2\biggl(\f{\g-1}{2A\g}\biggr)^{\f{1}{\g-1}}b_0^{\f{\g+1}{\g-1}}
\biggl(1+\a\biggr)>0,\\
&B_2=\biggl(\f{\g-1}{2A\g}\biggr)^{\f{1}{\g-1}}b_0^{\f{2}{\g-1}}\biggl(1
+\a\biggr),\\
&B_3=-(n-1)\biggl(\f{\g-1}{2A\g}\biggr)^{\f{1}{\g-1}}
b_0^{\f{\g+1}{\g-1}}(1+\a),\\
&\mu_1=\f{1}{2b_0}\biggl(1+\a\biggr)>0,\\
&\mu_2=-\f{n-1}{2}\biggl(1+\a\biggr)<0.
\endalign
$$}

\medskip

For the computational requirements later on,
we list some properties
of the partial Klainerman vector fields, which can be verified directly.
\medskip

{\bf Lemma 4.5.} {The partial Klainerman vector fields given by (4.1)
satisfy

\medskip

(i). $[Z_1, Z_2]=Z_3; \ [Z_2, Z_3]=Z_1; \ [Z_3,Z_1]=Z_2; \
[Z_i,Z_0]=0,  \ 1\leq i\leq n$.
\medskip

(ii). $[Z_i, \p_r]=0, \ 1\leq i\leq n; \ [Z_0, \p_r]=-\p_r.$
\medskip

(iii). $\ds Z_i(r)=0, \ Z_i(\frac{r}{t})=0, \ 1\leq i\leq n;  \
Z_0(\frac{1}{r})=-\frac{1}{r}; \ Z_0(\frac{r}{t})=0.$
\medskip

(iv). $\ds\na_x f\cdot\na_x g=\p_r f\cdot\p_r
g+\f{1}{r^2}\ds\sum_{i=1}^{2n-3}Z_i f\cdot Z_i g$ for any $C^1$ smooth
functions $f$ and $g$.

\medskip

(v). $|Z_i v|\leq r|\na v|$ for any $C^1$ smooth function $v$,
here $0\leq i\leq n$.
\medskip

(vi). $\ds\p_1=\f{x_1}{r}\p_r-\f{x_2}{r^2}Z_1+\f{x_3}{r^2}Z_3; \
\p_2=\f{x_2}{r}\p_r-\f{x_3}{r^2}Z_2+\f{x_1}{r^2}Z_1;   \
\p_3=\f{x_3}{r}\p_r-\f{x_1}{r^2}Z_3+\f{x_2}{r^2}Z_2.$}

\medskip

In addition, we now give three basic but important equalities (Lemma 4.6 and
Lemma 4.7 below), which will be used
in $\S 5$ to look for the multipliers in establishing a priori energy estimates.
\medskip

{\bf Lemma 4.6.} {\it Set $\Cal L_{\Phi}=\Cal L+\Cal P$ and $\Cal
M=A\p_t+B\p_r$ with $C^{\infty}$ smooth coefficients $A$ and $
B$, then for any $C^2$ smooth function $G$, one has
$$\align
&\Cal L_{\Phi}G\cdot \Cal M G\\
=&\ds\p_t N_0+\sum_{i=1}^{n}\p_i N_i+E_{00}(\p_t
G)^2+\sum_{i=1}^{n}E_{0i}\p_t G\p_i
G+\sum_{i,j=1}^{n}E_{ij}\p_i G\p_j G,\tag 4.14
\endalign$$
where
$$\align
N_0(\na G)&=\ds\f{1}{2} A(\p_t G)^2+ B\p_t
G\p_r G+ B\p_r G\sum_{k=1}^{n}p_k\p_k G-\f{1}{2}
A\sum_{i=1}^{n}\sum_{j=1}^{n}p_i p_j
\p_i G\p_j G+\f{1}{2} Ap_0\sum_{k=1}^{n}(\p_k G)^2,\\
N_i(\na G)&=\ds\f{x_i}{r}\biggl[\undersetbrace N^0
\to{-\f{1}{2} B(\p_t G)^2- B\p_t
G\sum_{k=1}^{n}p_k\p_k G-\f{1}{2}
B\sum_{j=1}^{n}\sum_{k=1}^{n}p_j p_k\p_j G\p_k
G+\f{1}{2} Bp_0\sum_{k=1}^{n}(\p_k G)^2}\biggr]\\
&\quad+p_i\biggl[\undersetbrace N^1\to{A(\p_t G)^2+
B\p_t G\p_r G+B\p_r G\sum_{j=1}^{n}p_j\p_j
G+A\p_t G\sum_{j=1}^{n}p_j\p_j
G}\biggr]\\
&\quad \ds-\p_i
G\biggl[\undersetbrace N^2 \to{ Ap_0\p_t G- Bp_0\p_r G}\biggr],\qquad  i=1,\cdots, n\\
\endalign$$
with
$$p=\Phi,\quad p_i=\p_i\Phi\quad (i=1, \cdots, n),\quad p_0=c^2(\rho).\tag
4.15$$

And the coefficients $E_{00}$, $E_{0i}$ and $E_{ij}$ in \rom{(4.14)}
are smooth.}

{\bf Proof.} Indeed, in term of (1.6) and (4.15), we have
$$\Cal L_{\Phi}=\p_t^2+\ds\sum_{k=1}^{n}p_k\p_k+\sum_{i,j=1}^{n}p_i
p_j\p_{ij}-p_0\Delta+\text{first order operators} .$$

This, together with a direct computation, yields (4.14). \qquad\qquad\qquad\qquad\qquad
\qquad\qquad\qquad
\qed\medskip

{\bf Lemma 4.7.} {\it Under the assumptions of \rom{Lemma 4.6}, then
$$\align
&\Cal L G\cdot\Cal M G\\
=&\p_t N_{00}+\p_r N_{r0}+\ds\f{n-1}{r}N_{r0}+\ds\sum_{i=1}^{2n-3}Z_i N_{\o 0}^i\\
&+K_{00}(\p_t G)^2+K_{0r}\p_t G\p_r G+K_{rr}(\p_r
G)^2+K_{n n}\ds\sum_{k=1}^{n}(Z_i G)^2\\
&+\ds K_{0n}\p_t G\cdot\ds\sum_{k=1}^{2n-3}Z_i A\cdot
Z_i G+K_{rn}\p_rG\cdot\ds\sum_{k=1}^{2n-3}Z_i B\cdot Z_i
G\tag 4.16
\endalign$$
and
$$\align
&\Cal P G\cdot\Cal M G\\
=&\p_tM_{0}+\ds\sum_{k=1}^{n}\p_kM_k+K^{00}(\p_t G)^2+\ds\sum_{i=1}^{n}K^{0i}\p_t G\p_i
G+\sum_{i,j=1}^{n}K^{ij}\p_i G\p_j G\tag 4.17
\endalign$$
with
$$\cases
N_{00}(\na G)=\f{1}{2} A(\p_t G)^2+B\p_t G\p_r G+B P_1(\p_r
G)^2-\f{1}{2}A P_2(\p_r
G)^2+\ds\sum_{i=1}^{2n-3}\f{1}{2r^2}A P_3(Z_i G)^2,\\
N_{r0}(\na G)=-\f{1}{2}B(\p_t G)^2+A P_1(\p_t G)^2+A P_2\p_t G\p_r
G+\f{1}{2} B P_2(\p_r
G)^2+\ds\sum_{i=1}^{2n-3}\f{1}{2r^2} B P_3(Z_i G)^2,\\
N_{\o 0}^{i}(\na G)=\ds -\f{1}{r^2}A P_3 \p_t G Z_i G-\f{1}{r^2}B
P_3\p_r G Z_i G,\qquad \quad i=1,\cdots, 2n-3
\endcases\tag 4.18$$
and
$$\cases
K_{00}=\ds-\f{1}{2}\p_t A+\f{1}{2}\p_r B-\p_r(AP_1)+\f{1}{r}AP_4+\f{n-1}{2r}B-\f{n-1}{r}AP_1,\\
K_{0r}=\ds-\p_t B-\p_r(AP_2)+\f{1}{r} AP_5+\f{1}{r} BP_4-\f{n-1}{r}AP_2,\\
K_{rr}=\ds-\p_t(BP_1)+\f{1}{2}\p_t(AP_2)-\f{1}{2}\p_r(BP_2)+\f{1}{r}BP_5-\f{n-1}{2r}BP_2,\\
K_{n
n}=\ds-\f{1}{2r^2}\p_t(AP_3)-\f{1}{2}\p_r(\f{1}{r^2}BP_3)-\f{n-1}{2r^3}BP_3,\\
K_{0n}=\ds\f{1}{r^2}P_3,\\
K_{rn}=\ds\f{1}{r^2}P_3
\endcases\tag 4.19$$
and
$$\cases
2K^{00}=\ds -\sum_{i=1}^{n}\p_i(A f_{0i}),\\
2K^{0i}=\ds-\ds\sum_{k=1}^{n}\p_k(\f{x_i}{r}B f_{0k})+\sum_{k=1}^{n}\p_k(\f{x_k}{r}B f_{0i})-\sum_{j=1}^{n}\p_j(A(f_{ji}+f_{ij})),\qquad 1\leq i\leq n,\\
2K^{ij}=\ds-\p_t(\f{x_i}{r}B f_{0j})-\sum_{k=1}^{n}\p_k(\f{x_i}{r}B
(f_{kj}+f_{jk})) +\p_t(A f_{ij})+\sum_{k=1}^{n}\p_k(\f{x_k}{r}B
f_{ij}),\quad 1\leq i,j\leq n.
\endcases\tag 4.20$$
In addition, the explicit expressions of the terms $M_i$ $(0\le i\le
n)$ in \rom{(4.17)} are not given here since this is not required.}

Based on the local existence result given in
$\S 3$, we will use the continuous induction method to prove Theorem 4.1.
For this end, a priori estimates on the solution $(\dot\vp,\xi)$ are required
to be established.
We now introduce some notations in order to fulfill the requirements
in $\S 5$-$\S 6$
below. For any $T_0>1$,
set
$$\align
&D_{T_0}=\{(t,r,\o): 1<t<T_0, \sigma(t, \omega) <r<\zeta(t,\o),\o\in
\Bbb
S^{n-1}\},\\
&B_{T_0}=\{(t,r,\o): 1<t<T_0, r=\sigma(t, \omega), \o\in\Bbb
S^{n-1}\},\\
&\Gamma_{T_0}=\{(t,r,\o): 1<t<T_0, r=\zeta(t,\o),\o\in\Bbb
S^{n-1}\},
\endalign$$
where $B_{T_0}$ and $\Gamma_{T_0}$ are the lateral boundaries of
$D_{T_0}$.

\vskip 0.4 true cm \centerline{\bf \S 5. Proof of Theorem 1.1 in the case of
$n=3$} \vskip 0.3 true cm

In this section, we will establish a uniform weighted energy estimate on
$(\dot\vp,\xi)$ for the problem (4.4) together with (4.9) and
(4.11)-(4.13) for $n=3$ in the domain $D_{T_0}$, which is defined in the above. In $\S 5.1$,
we will establish the first order
weighted energy estimate of $\na\dot\vp$. Subsequently, in $\S 5.2$,
the higher-order weighted energy estimates of $\na\dot\vp$ are
derived by utilizing the modified Klainerman's vector fields. Based on such energy
estimates given in $\S 5.1-\S
5.2$, we can use continuous induction argument given in $\S 5.3$ to
obtain the global existence and behavior at infinity of the
solution $(\dot\vp,\xi)$ and then complete the proofs of Theorem 4.1 and Theorem
1.1 in the case of $n=3$.
\medskip

{\bf $\S 5.1.$ First order weighted energy estimates} \medskip

 {\bf Theorem 5.1.} {\it For $n=3$, if $\dot\vp\in C^{2}(D_{T_0})$
is a solution of equation $(4.4)$ with
the fixed boundary condition $(4.9)$, and
$$|\xi|+\ds\sum_{l=0}^{1}t^{l}\bigl(|\na^{l+1}\dot\vp|+|\na^{l+1}\xi|\bigr)
\le M\ve\tag 5.1$$ holds
for small $\ve$,  $t\in [1, T_0]$ and some positive constant $M$. Then for any
fixed constant $\mu<-1
-\f{1}{2} \sqrt{\f{\gamma+7}{2}}$, we have
$$ \align &C_1  T_0^{\mu+1}\int_{\sigma(T_0, \omega)\le
r\le\zeta(T_0,\o)}|\nabla \dot \varphi|^2(T_0, x) dx+C_2\iint_{D_{T_0}} t^\mu
|\nabla \dot \varphi|^2 dtdx\\
& +C_3\int_{\Gamma_{T_0}}t^{\mu+1}(\p_t\dot\vp)^2
dS+C_4\int_{\Gamma_{T_0}}t^{\mu+1}\f{1}{r^2}\ds\sum_{k=1}^{3}(Z_k\dot\vp)^2
dS\\
&\le C\ve^2+C_5\int_{\Gamma_{T_0}}{t}^{\mu+1}(\Cal
B_0\dot\vp)^2dS,\tag 5.2\endalign
$$
where $C_i (1\le i\le 5)$ are some positive constants depending on
$b_0$ and $\gamma$. In particular,
$$\cases
&C_3=\ds\f{(\g-1)b_0^2}{8}\biggl(1+\a)\biggr),\\
&C_5=\ds\f{(\g-1)b_0^4}{2}\biggl(1+\a)\biggr).
\endcases\tag 5.3$$}

{\bf Remark 5.1.} {\it The values of constants $C_3$ and $C_5$ will
play an important role in the energy estimates for the problem $(4.4)$
with $(4.9)$ and $(4.11)$-$(4.13)$ since the most troublesome term
$\ds\int_{\Gamma_{T_0}}{t}^{\mu+1}(\Cal B_0\dot\vp)^2dS$  in $(5.2)$
will be shown to be absorbed by the positive integrals in the left
hand side of $(5.2)$. The reason which the term
$\ds\int_{\Gamma_{T_0}}{t}^{\mu+1}(\Cal B_0\dot\vp)^2dS$ is most
troublesome is: due to the Neumann boundary condition $(4.11)$ other
than the artificial Dirichlet boundary condition as in $[11]$, the
usual Poincar\'e inequality does not hold for the solution $\dot\vp$
$($it is noted that the boundary condition $(4.11)$ contains the
function $\xi$, which is roughly equivalent to $\ds\f{\dot\vp}{t}$
in terms of $(4.12)$, then the estimate on $\dot\vp$ on the shock
surface must be done$)$, namely, the $L^2(\Gamma_{T_0})-$estimates of
$\na\dot\vp$ on the shock surface $\Gamma_{T_0}$ can not be obtained
directly.}

\medskip

{\bf Proof.} For $\Cal M\dot\vp=A(t,x) \p_t \dot \varphi +
B(t,x)\p_r \dot\varphi$, it follows from Lemma 4.6-Lemma 4.7 that
$$\align
&\iint_{D_{T_0}} R_0(t,x)\cdot \Cal M\dot\vp dtdx=
\iint_{D_{T_0}}\bigl(\Cal L\dot\vp+\Cal P\dot\vp\bigr)\cdot\Cal M\dot\vp dtdx\\
=&\ds\iint_{D_{T_0}}I_1 dtdx+\int_{\sigma(T_0 ,\o)<r<
\zeta(T_0 ,\o)}N_0(\na\dot\vp)(T_0 ,x)dx-\int_{\sigma(1,\o)<r
<\zeta(1,\o)}N_0(\na\dot\vp)(1,x)dx\\
&\ds+\int_{\Gamma_{T_0}}\biggl(\sum_{i=1}^{3}(\f{x_i}{r}N_i(\na\dot\vp)
-\p_i\chi\cdot
N_i(\na\dot\vp))-\p_t\chi N_0(\na\dot\vp)\biggr) dS\\
&\ds+\int_{B_{T_0}}\biggl(\p_t\sigma
N_0(\na\dot\vp)-\sum_{i=1}^{3}(\f{x_i}{r}N_i(\na\dot\vp)-\p_i\sigma
N_i(\na\dot\vp))\biggr)dS\tag 5.4
\endalign$$
with
$$\align
I_1&=K_{00}(\p_t\dot\vp)^2+K_{0r}\p_t\dot\vp\p_r\dot\vp+K_{rr}(\p_r\dot\vp)^2
+K_{33}\f{1}{r^2}\ds\sum_{i=1}^{3}(Z_i\dot\vp)^2\\
&+K_{03}\f{1}{r^2}\ds\sum_{i=1}^{3}Z_i A
Z_i\dot\vp\p_t\dot\vp+K_{r3}\f{1}{r^2}\ds\sum_{i=1}^{3}Z_i B
Z_i\dot\vp\p_r\dot\vp\\
&+K^{00}(\p_t\dot\vp)^2+\ds\sum_{k=1}^{3}K^{0i}\p_t\dot\vp\p_i\dot\vp
+\ds\sum_{i,j=1}^{3}K^{ij}\p_i\dot\vp
\p_j\dot\vp.\tag 5.5
\endalign$$

Our purpose is to choose suitable functions $A(t,x)$ and $
B(t,x)$ so that all integrals on $D_{T_0}, B_{T_0}$ and $t=T_0$ in
the right hand side of (5.4) are non-negative and the integral on
$\Gamma_T$ gives a ``good'' control on $\dot\vp$. From now on, we
will derive some sufficient conditions for the choices of $
A(t,x)$ and $B(t,x)$ in the process of analyzing each integral,
and then $A(t,x)$ and $B(t,x)$ can be determined. This
process is divided into the following five steps.
\smallskip

{\bf Step 1. The analysis on the term
$\ds\int_{B_{T_0}}\biggl(\p_t\sigma
N_0(\na\dot\vp)-\sum_{i=1}^{3}(\f{x_i}{r}N_i(\na\dot\vp)-\p_i\sigma
N_i(\na\dot\vp))\biggr)dS$}
\smallskip

Due to (1.7) and (4.9), using the notations in (4.14)-(4.15), we
have on $B_{T_0}$
$$\align
&\p_t\sigma N_0(\na\dot\vp)-\ds\sum_{i=1}^{3}(\f{x_i}{r}N_i(\na\dot\vp)
-\p_i\sigma N_i(\na\dot\vp))\\
=&\p_t\sigma N_0(\na\dot\vp)-N^0(\na\dot\vp)-\p_t\sigma N^1(\na\dot\vp)\\
=&(B-\p_t\sigma
A)\biggl[\f{1}{2}(\p_t\dot\vp)^2+\p_t\dot\vp\ds\sum_{i=1}^{3}p_i\p_i\dot\vp
+\f{1}{2}p_0|\na\dot\vp|^2
-\f{1}{2}\sum_{i=1}^{3}\sum_{j=1}^{3}p_i
p_j\p_i\dot\vp\p_j\dot\vp\biggr].
\endalign$$

It follows from the assumption (5.1), (4.15) and (iv) in Lemma 2.1
that $\ds\sum_{i=1}^{3}|p_i|^2>p_0$, and thus
$$p_0|\na\dot\vp|^2-\sum_{i=1}^{3}\sum_{j=1}^{3}p_i
p_j\p_i\dot\vp\p_j\dot\vp\leq
(p_0-\ds\sum_{i=1}^{3}|p_i|^2)|\na\dot\vp|^2<0,$$ and unfortunately,
$\ds\f{1}{2}(\p_t\dot\vp)^2+\p_t\dot\vp\ds\sum_{i=1}^{3}p_i\p_i\dot\vp
+\f{1}{2}p_0|\na\dot\vp|^2
-\f{1}{2}\sum_{i=1}^{3}\sum_{j=1}^{3}p_i p_j\p_i\dot\vp\p_j\dot\vp$
may change its sign on $B_{T_0}$.

So in order to control the term $\ds\int_{B_{T_0}}\biggl(\p_t\sigma
N_0(\na\dot\vp)-\sum_{i=1}^{3}(\f{x_i}{r}N_i(\na\dot\vp)-\p_i\sigma
N_i(\na\dot\vp))\biggr)dS$ well, $A$ and $B$ should
satisfy
$$B=\p_t\sigma A\quad \text{on}\quad r=\si(t,\o),\tag 5.6$$
and then
$$\ds\int_{B_{T_0}}\biggl(\p_t\sigma
N_0(\na\dot\vp)-\sum_{i=1}^{3}(\f{x_i}{r}N_i(\na\dot\vp)-\p_i\sigma
N_i(\na\dot\vp))\biggr)dS=0.\tag 5.7$$

In view of (5.6) and the self-similar property of the background
solution $\hat\Phi$, we set
$$A=t^{\mu}r a(\f{r}{t}),\quad B=t^{\mu+1}b_{\sigma}(t,r,w),\tag 5.8$$
with the functions $a$ and $b_\si$ to be determined later. By (5.6),
$a$ and $b_{\sigma}$ in (5.8) should  satisfy the following
restriction
$$b_{\sigma}=sa\p_t\sigma\quad \text{on}\quad r=\sigma(t,\o).\tag 5.9$$
\medskip

{\bf Step 2.} {\bf Positivity of
$\ds\int_{\sigma(T_0,\o)<r<\zeta(T_0,\o)}N_0(\na\dot\vp)(T_0,x)dx$}
\smallskip

Under the assumption (5.1), using (vi) in Lemma 4.5, one then has
$$\align
N_0(\na\dot\vp)&=N_{00}(\na\dot\vp)+t^{\mu+1}O(\ve)|\na\dot\vp|^2\\
&=t^{\mu+1}\biggl[\f{1}{2}sa(\p_t\dot\vp)^2+b_{\sigma}\p_t\dot\vp\p_r\dot\vp+b_{\sigma}
P_1(s)(\p_r\dot\vp)^2-\f{1}{2}sa P_2(s)(\p_r\dot\vp)^2\\
&\quad\quad\quad+\f{1}{2r^2}a
P_3(s)\sum_{i=1}^{3}(Z_i\dot\vp)^2+O(\ve)|\na\dot\vp|^2\biggr].\tag
5.10
\endalign$$

Due to $P_3(s)>0$ by Lemma 4.3, to ensure the positivity of
quadratic polynomial  of $\ds(\p_t\dot\vp, \p_r\dot\vp,
\f{1}{r}Z_1\dot\vp,$ $ \ds\f{1}{r}Z_2\dot\vp,\f{1}{r}Z_3\dot\vp)$ in
$N_0(\na\dot\vp)$, $a$ and $b_{\sigma}$ should fulfill
$$
\cases &\ds a(s)>0,\\
&\ds b_\sigma^2-2s P_1(s)a(s)b_\sigma+s^2 P_2(s) a^2(s)<0,
\endcases$$
which is equivalent to
$$a(s)>0,\quad \hat u(s)-c(\hat\rho)(s)<\f{b_\sigma}{sa(s)}<\hat u(s)+c(\hat\rho)(s)
\quad \text{in}\quad \O_+.\tag 5.11$$

In this case, we arrive at
$$\int_{\sigma(T_0, \omega)\le r\le \zeta(T_0,\o)} N_0(\na\dot\vp)(T_0,x) dS
\ge C(b_0,\g)T_0^{\mu+1}\int_{\sigma(T_0, \omega)<
r<\zeta(T_0,\o)}|\na\dot\vp|^2(T_0,x) dS.\tag 5.12
$$\medskip

{\bf Step 3.} {\bf Positivity of the integral on $D_{T_0}$}
\medskip

Under the constraints (5.9) and (5.11), we will choose $a$ and
$b_\sigma$ such that
$$I_1\ge 0,\tag 5.13$$

First, it can be verified directly that under the assumption (5.1), we have
$$f_{ij}=O_0^{1}(\ve),\qquad  0\leq i\leq 3; 1\leq j\leq 3,$$
then it follow from (4.20) that
$$\biggl|K^{00}(\p_t\dot\vp)^2+\ds\sum_{k=1}^{3}K^{0i}\p_t\dot\vp\p_i\dot\vp
+\ds\sum_{i,j=1}^{3}K^{ij}\p_i\dot\vp
\p_j\dot\vp\biggr|\leq t^{\mu}O(\ve)|\na\dot\vp|^2.\tag 5.14$$

On the other hand, according to (5.9) and (5.11), we choose
$$\cases
&a(s)=1,\\
&\ds b_\si(t,\o)=s^2\biggl(1+\f{\ve\p_t
b-\ds\f{\ve}{t}b}{s}+\f{e}{b_0}(s-b_0-\f{\ve}{t}b)\biggr),
\endcases\tag 5.15$$
where the constant $e$ will be determined later on.

It follows from Lemma 2.1, Lemma 4.3, (4.19), (5.15) and a direct
computation that
$$\cases
&K_{00}=t^{\mu}\biggl(\ds\f{1}{2}(2+e-\mu)b_0(1+\a)\biggr),\\
&K_{0r}=t^{\mu}\biggl(\ds\bigl(\f{\g+3}{2}+e-\mu\bigr)b_0^2(1+\a)\biggr),\\
&K_{rr}=t^{\mu}\biggl(\ds(\f{\g+1}{4}e-\f{\g+1}{4}\mu+1)b_0^3(1+\a)\biggr),\\
&K_{33}=t^{\mu}\biggl(\ds-\f{\g-1}{4}(2+e+\mu)b_0^3(1+\a)\biggr),\\
&K_{03}Z_iA=0,\qquad  i=1,2,3,\\
&K_{r3} Z_iB=t^{\mu}O(\ve),\qquad  i=1,2,3,\\
&K_{0r}^2-4K_{00}
K_{rr}=t^{2\mu}\biggl(\ds\f{\g-1}{4}\bigl(\g+7-2(e-\mu)^2\bigr)(1+\a)\biggr).
\endcases\tag
5.16$$

Therefore, for large $b_0$, with (5.5) and (5.14), the sufficient
conditions for (5.13) can be chosen as
$$2+e-\mu>0,\quad 2+e+\mu<0, \quad \g+7-2(e-\mu)^2<0 .\tag 5.17$$

If we let
$$\mu<-1-\f{1}{2}\sqrt{\f{\g+7}{2}}\tag 5.18$$
and
$$e=\f{1}{2}\sqrt{\f{\g+7}{2}}-1,\tag 5.19$$
then it can be verified from (5.16)-(5.17) directly that
$$K_{00}>0,\quad K_{0r}^2-4K_{00} K_{rr}<0,\quad K_{33}>0.$$

Combining this with (5.14)  shows (5.13).

In this case, we arrive at
$$\iint_{D_{T_0}} I_1 dtdx
\geq C(b_0, \gamma)\iint_{D_{T_0}}t^{\mu} \biggl((\p_{t}\dot\vp)^2+
(\p_{r}\dot\vp)^2+\f{1}{r^2}\ds\sum_{k=1}^{3}(Z_k\dot\vp)^2\biggr)dtdx.\tag
5.20$$
\medskip

{\bf Step 4. Estimate of
$\displaystyle\int_{\Gamma_{T_0}}\biggl(\sum_{i=1}^{3}(\f{x_i}{r}N_i(\na\dot\vp)-\p_i\chi
N_i(\na\dot\vp))-\p_t\chi N_0(\na\dot\vp)\biggr)dS$}
\medskip

By (iv) in Lemma 4.5, the assumption (5.1) and the definition of
$\dot\vp$ and $\xi$ (also see (4.18)),
$$\align
&\ds\sum_{i=1}^{3}(\f{x_i}{r}N_i(\na\dot\vp)-\p_i\chi N_i(\na\dot\vp))-\p_t\chi N_0(\na\dot\vp)\\
=&N^0(\na\dot\vp)+\p_r\Phi N^1(\na\dot\vp)-\p_r\dot\vp
N^2(\na\dot\vp)\\
&-\f{1}{r^2}\ds\sum_{i=1}^{3}Z_i\chi\cdot Z_i\Phi\cdot
N^1(\na\dot\vp)+\f{1}{r^2}\ds\sum_{i=1}^{3}Z_i\chi\cdot Z_i\dot\vp
N^2(\na\dot\vp)-\p_t\chi
N_0(\na\dot\vp)\\
=&N^0+P_1(s) N^1-\p_r\dot\vp N^2-s_0
N_0+t^{\mu+1}O(\ve)|\na\dot\vp|^2.\tag 5.21
\endalign$$

By the assumption (5.1),  (iv) in Lemma 4.5 and the
notations in (4.15), we have
$$\cases
&N^0=\ds t^{\mu+1}\biggl(-\f{1}{2}b_\si(\p_t\dot\vp)^2-b_\si\hat
u(s)\p_t\dot\vp\p_r\dot\vp-\f{1}{2}b_\si\hat
u^2(s)(\p_r\dot\vp)^2\\
&\quad\quad\quad\ds+\f{1}{2}b_\si
p_0(\p_r\dot\vp)^2+\f{1}{2r^2}b_\si
p_0\sum_{k=1}^{3}(Z_k\dot\vp)^2\biggr),\\
&N^1=t^{\mu+1}\biggl(s(\p_t\dot\vp)^2+b_\si\p_t\dot\vp\p_r\dot\vp+b_\si\hat
u(s)(\p_r\dot\vp)^2+s\hat
u(s)\p_t\dot\vp\p_r\dot\vp+O(\ve)|\na\dot\vp|^2\biggr),\\
&N^2=t^{\mu+1}\biggl(sp_0\p_t\dot\vp+b_\si p_0\p_r\dot\vp\biggr).
\endcases\tag 5.22$$

Combining (5.22) with (5.10), (5.15) and (5.21) yields
$$\align
&\ds\sum_{i=1}^{3}(\f{x_i}{r}N_i(\na\dot\vp)-\p_i\chi N_i(\na\dot\vp))-\p_t\chi N_0(\na\dot\vp)\\
=&t^{\mu+1}\biggl(\beta_{11}(\p_t\dot\vp)^2+\beta_{12}\p_t\dot\vp\p_r\dot\vp
+\beta_{13}(\p_r\dot\vp)^2+\beta_{14}\f{1}{r^2}\ds\sum_{k=1}^{3}(Z_k\dot\vp)^2\biggr)\tag
5.23
\endalign$$
with
$$\cases
&\beta_{11}=\ds b_0^2(1+\a),\\
&\beta_{12}=\ds -\f{\g-1}{2}b_0^3(1+\a),\\
&\beta_{13}=\ds -\f{\g-1}{2}b_0^4(1+\a),\\
&\beta_{14}=\ds \f{\g-1}{4}eb_0^3(s_0-b_0)(1+\a).
\endcases$$

Due to $\p_r\dot\vp=\Cal B_0\dot\vp-\mu_1\p_t\dot\vp$, from (5.23)
and (5.19) that
$$\align
&\ds\sum_{i=1}^{3}(\f{x_i}{r}N_i(\na\dot\vp)-\p_i\chi\cdot N_i(\na\dot\vp))
-\p_t\chi N_0(\na\dot\vp)\\
\geq&
t^{\mu+1}\biggl(\beta^{11}(\p_t\dot\vp)^2+\beta^{12}\p_t\dot\vp\cdot\Cal
B_0\dot\vp+\beta^{13}(\Cal
B_0\dot\vp)^2+\beta^{14}\ds\sum_{k=1}^{3}(Z_k\dot\vp)^2\biggr)
\endalign$$

with
$$\cases
&\ds\beta^{11}=\beta_{1}-\mu_1\beta_{12}+\mu_1^2\beta_{13}=\f{\g-1}{8}b_0^2(1+\a),\\
&\beta^{12}=\beta_{12}-2\mu_1\beta_{13}=b_0^3\biggl(\a\biggr),\\
&\beta^{13}=\beta_{13}=\ds-\f{\g-1}{2}b_0^4(1+\a),\\
&\beta^{14}=\beta_{14}=\ds\f{\g-1}{4}eb_0^3(s_0-b_0)(1+\a)>0.
\endcases$$

Consequently, one has
$$\align
&\int_{\Gamma_{T_0}}\biggl(\sum_{i=1}^{3}(\f{x_i}{r}N_i(\na\dot\vp)-\p_i\chi
N_i(\na\dot\vp))
-\p_t\chi N_0(\na\dot\vp)\biggr)dS\\
&\geq
\f{\g-1}{8}b_0^2\biggl(1+\a\biggr)\int_{\Gamma_{T_0}}t^{\mu+1}|\p_t\dot\vp|^2 dS\\
&\quad -\f{\g-1}{2}b_0^4\biggl(1+\a\biggr)\int_{\Gamma_{T_0}}t^{\mu+1}(\Cal
B_0\dot\vp)^2 dS\\
&\quad +C(b_0,
\gamma)\int_{\Gamma_{T_0}}t^{\mu+1}\f{1}{r^2}\ds\sum_{k=1}^{3}(Z_k\dot\vp)^2
dS.\tag5.24\endalign$$ \medskip

{\bf Step 5. The estimates on $\displaystyle\int_{\sigma(1, \omega)
\leq r\leq \zeta(1,\o)} N_0(\na\dot\vp)(1,x)dx$ and
$\ds\iint_{D_{T_0}}R_0(t,x)\cdot\Cal M\dot\vp dtdx$}

\medskip

It follows from (4.13) that
$$\align
&\biggl|\int_{\sigma(1, \omega)\leq r\leq
\zeta(1,\o)}N_0(\na\dot\vp)(1,x)dx\biggr| \leq C(b_0,
\gamma)\int_{\sigma(1, \omega)\leq r\leq
\zeta(1,\o)}|\na\dot\vp|^2(1,x) dx.\tag 5.25\endalign$$

Moreover, it follows from Lemma 4.2 and $\mu<-2$ in (5.19) that
$$\align
&\biggl|\iint_{D_{T_0}}R_0(t,x)\cdot \Cal M\dot\vp dtdx\biggr|\leq
C(b_0,\g)\ve\biggl(\iint_{D_{T_0}}t^{\mu}|\na\dot\vp|^2
dtdx+\ve\biggr).\tag 5.26
\endalign$$

Finally, substituting the estimates  (5.7), (5.12), (5.20) and
(5.24)-(5.26) into (5.4), then (5.2) and (5.3) can be obtained.
Therefore, Theorem 5.1 is proved.\qquad\qquad\qquad\qquad\qquad
\qquad\qquad\qquad\qquad\qquad\qquad\qquad\qed

\medskip

Based on Theorem 5.1, we will derive the first order uniform energy
estimate of $\na\dot\vp$. For this end, we require an improved
Hardy-type inequality on $\ds\int_{\Gamma_{T_0}}t^{\mu-1}|\dot\vp|^2
dS$ in terms
of the special structures of (4.11) and (4.12),
which is motivated by Theorem 330 in [12] and [19].

\medskip

{\bf Lemma 5.2. (Improved Hardy-type inequality)} {\it Under the
assumptions of \rom{Theorem 5.1},  for
$\mu<-1-\ds\f{1}{2}\sqrt{\ds\f{\g+7}{2}}$, we have
$$
\align \int_{\Gamma_{T_0}}t^{\mu-1}|\dot\vp|^2 dS& \le   C(b_0,
\gamma)\ve^2+\f{1}{\mu^2}\bigl(1+\a\bigr)
\int_{\Gamma_{T_0}}t^{\mu+1}(\p_t\dot\vp)^2dS\\
& +C(b_0, \gamma)\ve\biggl(\int_{\Gamma_{T_0}}t^{\mu+1}(\Cal
B_0\dot\vp)^2
dS+\f{1}{r^2}\int_{\Gamma_{T_0}}t^{\mu+1}\ds\sum_{k=1}^{3}(Z_k\dot\vp)^2
dS\biggr).\tag 5.27
\endalign
$$}\medskip

{\bf Remark 5.2.} {\it Classical Hardy-type inequality in \rom{[12]}
only means that
$$\int_{\Gamma_{T_0}}t^{\mu-1}|\dot\vp|^2 dS\le \ds\f{2}{\mu^2}
\bigl(1+O(b_0^{-\f{2}{\g-1}})
+O(b_0^{-2})\bigr)
\int_{\Gamma_{T_0}}t^{\mu+1}(\p_t\dot\vp)^2dS+\text{``some small
terms''}\tag5.28$$ holds. One should specially notice that the
coefficient of $\ds\int_{\Gamma_{T_0}}t^{\mu+1}(\p_t\dot\vp)^2dS$ is
$\ds\f{2}{\mu^2}$ in \rom{(5.28)} other than $\ds\f{1}{\mu^2}$ in
\rom{(5.27)}. If so, it will completely fail in deriving the first
order energy estimates on $\dot\vp$ from \rom{Theorem 5.1} since the
terms on the shock surface can not be absorbed by the left hind side
terms in (5.2).}\medskip

{\bf Proof.} It is noted that
$$\int_{\Gamma_{T_0}}t^{\mu-1}|\dot\vp|^2 dS
=\int_{\Bbb S^2}d\o\int_1^{T_0}
t^{\mu-1}|\dot\vp(t,\chi(t,\o),\o)|^2dt.$$

Set $m(\o)\equiv\ds \int_1^{T_0}
t^{\mu-1}|\dot\vp(t,\chi(t,\o),\o)|^2dt$, then by integration by
parts
$$
\align
m&(\o)=\ds\f{1}{\mu}t^{\mu}|\dot\vp(t,\chi(t,\o),\o)|^2\biggr|_{t=1}^{t=T_0}
-\f{2}{\mu}\int_1^{T_0} t^{\mu}\dot\vp(t,\chi(t,\o),\o)(\p_t\dot\vp
+\p_t\chi\p_r\dot\vp)(t,\chi(t,\o),\o)dt\\
&\le\f{1}{|\mu|}|\dot\vp(1,\chi(1,\o),\o)|^2-\f{2}{\mu}\int_1^{T_0}
t^{\mu}\dot\vp(t,\chi(t,\o),\o)\bigl(\p_t\chi\Cal
B_0\dot\vp+(1-\mu_1\p_t\chi)\p_t\dot\vp\bigr)(t,\chi(t,\o),\o)dt.\tag
5.29
\endalign$$

Due to (5.1), Lemma 2.1,  Lemma 4.4 and
$\p_t\chi=b_0(1+\a)$,  then
$$1-\mu_1\p_t\chi=\f{1}{2}(1+\a+O(\ve))>0.\tag 5.30$$

It follows from $\mu_3<0$ in (4.12), $\mu_2<0$ in Lemma 4.4, Lemma
4.2, (5.18), (5.29) and a direct computation that

$$
\align &m(\o)\\
&\le C(b_0, \gamma)\ve^2-\f{2}{\mu}\int_1^{T_0}
t^{\mu}\dot\vp(t,\chi(t,\o),\o)(1-\mu_1\p_t\chi)
\p_t\dot\vp(t,\chi(t,\o),\o)dt\\
&\quad +\f{2\mu_2}{\mu}\int_1^{T_0}
t^{\mu}\p_t\chi\dot\vp(t,\chi(t,\o),\o)\xi dt-\f{2}{\mu}\int_1^{T_0}
t^{\mu}\p_t\chi\dot\vp(t,\chi(t,\o),\o)\bigl(\f{\kappa(\xi,
\na\dot\vp)+R_1(t,x)}{
B_1}\bigr)dt\\
&=C(b_0, \gamma)\ve^2-\f{2}{\mu}\int_1^{T_0}
t^{\mu}\dot\vp(t,\chi(t,\o),\o)(1-\mu_1\p_t\chi)
\p_t\dot\vp(t,\chi(t,\o),\o)dt\\
&\quad +\f{2\mu_2}{\mu}\int_1^{T_0} t^{\mu}\p_t\chi
\bigl(\mu_3(t,x)t\xi\bigr)\xi dt-\f{2}{\mu}\int_1^{T_0}
t^{\mu}\p_t\chi\dot\vp(t,\chi(t,\o),\o)\bigl(\f{\kappa(\xi,
\na\dot\vp)+R_1(t,x)}{
B_1}\bigr)dt\\
&\le C(b_0, \gamma)\ve^2+\f{1}{2}\int_1^{T_0}
t^{\mu-1}|\dot\vp(t,\chi(t,\o),\o)|^2dt
+\f{2}{\mu^2}\int_1^{T_0} t^{\mu+1}(1-\mu_1\p_t\chi)^2|\p_t\dot\vp(t,\chi(t,\o),\o)|^2dt\\
&\quad +C(b_0, \gamma)\ve\biggl(\int_1^{T_0}
t^{\mu-1}|\dot\vp(t,\chi(t,\o),\o)|^2dt+
\int_1^{T_0} t^{\mu+1}|\xi|^2dt+\int_1^{T_0} t^{\mu+1}\bigl(|\na_{t,r}\dot\vp|^2\\
&\quad
+\f{1}{r^2}\ds\sum_{i=1}^{3}(Z_i\dot\vp)^2\bigr)(t,\zeta(t,\o),\o)dt\biggr).
\endalign
$$

This, together with (5.30), yields
$$\align
\int_{\Gamma_{T_0}}t^{\mu-1}|\dot\vp|^2 dS&\le C(b_0, \gamma)\ve^2
+(\f{4}{\mu^2}+O(\ve))\int_{\Gamma_{T_0}}t^{\mu+1}(1-\mu_1\p_t\chi)^2(\p_t\dot\vp)^2
dS\\
&\quad +C(b_0, \gamma)\ve\biggl(\int_{\Gamma_{T_0}}t^{\mu+1}(\Cal
B_0\dot\vp)^2
dS+\ds\int_{\Gamma_{T_0}}t^{\mu+1}\f{1}{r^2}\ds\sum_{k=1}^{3}(Z_k\dot\vp)^2 dS\biggr)\\
&\le C(b_0, \gamma)\ve^2+\f{1}{\mu^2}\bigl(1+\a\bigr)
\int_{\Gamma_{T_0}}t^{\mu+1}(\p_t\dot\vp)^2dS\\
&\quad +C(b_0, \gamma)\ve\biggl(\int_{\Gamma_{T_0}}t^{\mu+1}(\Cal
B_0\dot\vp)^2
dS+\ds\int_{\Gamma_{T_0}}t^{\mu+1}\f{1}{r^2}\ds\sum_{k=1}^{3}(Z_k\dot\vp)^2
dS\biggr).
\endalign
$$

Consequently, Lemma 5.2 is proved.\qquad\qquad\qquad\qquad\qquad\qquad\qquad
\qquad\qquad\qquad\qquad\qquad\qquad \qquad \qed
\medskip

Lemma 5.2 illustrates that using the special structures in (4.11)-(4.12),
we can obtain an improved Hardy-type inequality, which plays a crucial role
in establishing the first order weighted energy estimates of
$\na\dot\vp$. With Lemma 5,2, we have\medskip

{\bf Theorem 5.3.} {\bf (First Order Weighted Energy Estimate.)}
{\it Under the assumptions of \rom{Theorem 5.1}, for
$\mu<-1-\ds\f{1}{2}\sqrt{\ds\f{\g+7}{2}}$, one has
$$ \align
&T_0^{\mu+1}\int_{\sigma(T_0, \omega)\le
r\le\zeta(T_0,\o)}|\na\dot\vp|^2(T_0,x) dx
+\iint_{D_{T_0}}{t}^{\mu}|\na\dot\vp|^2
dtdx+\int_{\Gamma_{T_0}}{t}^{\mu+1}|\na\dot\vp|^2dS \leq
C\ve^2.\tag 5.31\endalign
$$}
\medskip

{\bf Proof.} To obtain (5.31), we require to give a
delicate estimate on the term
$C_5\displaystyle\int_{\Gamma_{T_0}}t^{\mu+1}(\Cal B_0\dot\vp)^2 dS$
in (5.2) so that it can be absorbed by the corresponding positive
terms in the left hand side of (5.2).

We now treat the term
$\displaystyle\int_{\Gamma_{T_0}}t^{\mu+1}(\Cal B_0\dot\vp)^2 dS$.

>From (4.11) and Lemma 4.2, we have
$$\align \ds\int_{\Gamma_{T_0}}t^{\mu+1}(\Cal
B_0\dot\vp)^2dS=&
\ds\int_{\Gamma_{T_0}}t^{\mu+1}\biggl(\f{\kappa(\xi,\nabla\dot\vp)+R_1(t,x)}{\Cal
B_1}-\mu_2\xi \biggr)^2dS\\
\leq&\mu_2^2(1+\a)\ds\int_{\Gamma_{T_0}}
t^{\mu-1}|t\xi|^2dS\\
&+C(b_0, \gamma)\int_{\Gamma_{T_0}}
t^{\mu+1}\kappa^2(\xi,\nabla\dot\vp)dS+C(b_0, \gamma) \ve^2.\tag
5.32
\endalign
$$

Due to the assumption (5.1) and the property of
$\kappa(\xi,\na\dot\vp)$ on $\Gamma_{T_0}$, then
$$
\align &\ds\int_{\Gamma_{T_0}}
t^{\mu+1}\kappa^2(\xi,\nabla\dot\vp)dS\\
&\ds\leq C(b_0,
\gamma)\ve^2\int_{\Gamma_{T_0}}t^{\mu+1}\biggl(|\xi|^2
+|\nabla\dot\vp|^2\biggr)dS\\
&\leq\ds C(b_0,
\gamma)\ve^2\int_{\Gamma_{T_0}}t^{\mu+1}\biggl(|\xi|^2+(\Cal
B_0\dot\vp)^2+(\p_{t}\dot\vp)^2
+\f{1}{r^2}\ds\sum_{i=1}^{3}(Z_i\dot\vp)^2\biggr)dS.\tag
5.33\endalign$$

In addition, due to the smallness of $\ve$ and the boundary
condition (4.12) together with Lemma 2.1 (ii), then we have from
(5.32)-(5.33) that
$$\align
\int_{\Gamma_{T_0}}t^{\mu+1}(\Cal B_0\dot\vp)^2 dS  &\leq
\f{\mu_2^2}{b_0^2}\biggl(1+\a\biggr)\int_{\Gamma_{T_0}}t^{\mu-1}|\dot\vp|^2
dS\\
&\quad +C(b_0,
\gamma)\ve^2\int_{\Gamma_{T_0}}t^{\mu+1}\biggl((\p_t\dot\vp)^2
+\f{1}{r^2}\ds\ds\sum_{k=1}^{3}(Z_k\dot\vp)^2\biggr)dS+C(b_0,
\gamma) \ve^2.
\endalign$$

Therefore, combining  this with (5.27) yields
$$\align
&\int_{\Gamma_{T_0}}t^{\mu+1}|\Cal B_0\dot\vp|^2 dS\\
&\leq
C(b_0, \gamma)\ve^2+\f{1}{\mu^2 b_0^2}\biggl(1+\a\biggr)
\int_{\Gamma_{T_0}}t^{\mu+1}|\p_t\dot\vp|^2 dS\\
&\quad +C(b_0,
\gamma)\ve\int_{\Gamma_{T_0}}t^{\mu+1}\f{1}{r^2}\ds\sum_{k=1}^{3}(Z_k\dot\vp)^2
dS.\tag 5.34
\endalign$$

Substituting (5.34) into (5.2), we obtain from the assumption in
Theorem 5.1 that
$$ \align &  T_0^{\mu+1}\int_{\sigma(T_0, \omega)\le
r\le\zeta(T_0,\o)}|\na\dot\vp|^2(T_0,x) dx+\iint_{D_{T_0}}t^{\mu}|\na\dot\vp|^2 dtdx\\
&+\f{\g-1}{8}(1-\f{4}{\mu^2})b_0^2\biggl(1+\a\biggr)
\int_{\Gamma_{T_0}}t^{\mu+1}(\p_t\dot\vp)^2
dS\\
&+\int_{\Gamma_{T_0}}t^{\mu+1}\f{1}{r^2}\ds\sum_{k=1}^{3}(Z_k\dot\vp)^2
dS \le C\ve^2.\tag 5.35\endalign
$$

Due to $\mu^2>4$ in (5.18), then (5.35) yields Theorem
5.3.\qquad\qquad\qquad\qquad\qquad\qquad\qquad\qquad\qquad\qquad \qed

\vskip 0.3 true cm

{\bf $\S 5.2.$ Higher order weighted energy estimate.}

\vskip 0.3 true cm

In this subsection, we will derive the higher-order energy estimates of
$\dot\vp$, so that the decay properties of
$\na\dot\vp$ and $\xi$ for large $t$ can be established.

Denote
$$ \t Z=\{\t Z_j:  0\leq j \leq 2n-3\},\quad \text{where}\quad \t Z_0
=t\p_t+ (r+t^2\p_t b) \p_r, \quad \t Z_k=Z_k + t Z_k b \cdot \p_r \
(1\leq j\leq 2n-3)\tag 5.36
$$ by the vector
fields which are tangent to the surface $r=\sigma(t,\o)$.

Under the
assumptions in Theorem 1.1, one has
$$\t Z_i-Z_i=O_1^{\infty}(\ve)\na,\qquad\quad 0\leq i\leq 2n-3.\tag 5.37$$

\medskip

{\bf Lemma 5.4.}  {\it For $n=2,3$, there exist functions
$\tau_{ij}=\tau_{ij}(t,r,\omega)\,(0 \le i,j\le 2n-3)$, such that
for $i,j=0,1, \cdots, 2n-3,$
$$\tau_{ij}\in O_0^{\infty}(1),\quad \text{and}\quad \Cal B_{\sigma}(\t Z_i
+\sum_{j=0}^{2n-3} \tau_{ij}\t Z_j)\dot\varphi=0 \quad
\text{on}\quad r=\si(t,\omega).\tag 5.38$$}

{\bf Proof.} It follows from a direct computation that for $0 \le i
\le 2n-3,$
$$\align
[\Cal B_{\sigma}, \t Z_i]=&h_{i0}\Cal B_{\sigma}+\sum_{j=1}^{2n-3} h_{ij}\t Z_j,\\
\endalign$$
with the smooth function $h_{ij}(t,r,\omega) (0\leq i,j\leq 2n-3)$
satisfying
$$h_{i0}=O_0^{\infty}(1),\quad h_{ij}=O_{-1}^{\infty}(1),\qquad 1\leq j\leq 2n-3.\tag 5.39$$

Due to $\Cal B_{\sigma}\dot\varphi=0$ and $\t Z_i \Cal
B_{\sigma}\dot\varphi=0$ $(0\leq i\leq 2n-3)$ on $r=\si(t,\omega)$
and the fact that $\t Z_i$ is tangent to the surface $r=\si(t,\o)$, then for
smooth function $\tau_{ij}$ $(0\leq i,j\leq 2n-3)$, we have on
$r=\si(t,\omega)$
$$\align
&\Cal B_{\sigma}(\t Z_i+\sum_{j=0}^{2n-3}\tau_{ij}\t Z_j)\dot\varphi
\\
=&[\Cal B_{\sigma}, \t Z_i]\dot\varphi+\sum_{j=0}^{2n-3} \Cal
B_{\sigma}\tau_{ij}\cdot\t Z_j\dot\varphi +\sum_{j=0}^{2n-3}
\tau_{ij}[\Cal B_{\sigma}, \t Z_j]\dot\varphi
\\
=&\sum_{j=0}^{2n-3} h_{ij}\t Z_j\dot\varphi+\sum_{j=0}^{2n-3} \Cal
B_{\sigma}\tau_{ij}\cdot\t Z_j\dot\varphi
+\sum_{j=0}^{2n-3} \tau_{ij}\sum_{k=0}^{2n-3} h_{jk}\t Z_k\dot\varphi\\
=&\sum_{j=0}^{2n-3} \biggl(\Cal
B_{\sigma}\tau_{ij}+\sum_{k=0}^{2n-3}
h_{kj}\tau_{ik}+h_{ij}\biggr)\t Z_j\dot\varphi.\tag 5.40
\endalign$$

For $0 \le i,j \le 2n-3$, let the functions $\tau_{ij}$ satisfy
$$\cases
&\ds \Cal B_{\sigma}\tau_{ij}+h_{i j}=0\quad \text{in}\quad \O_+,\\
&\tau_{ij}=0\qquad\text{on}\quad r=\si(t,\omega).
\endcases\tag 5.41$$

By (5.39) and the analogous proof procedure on Lemma B.1 in Appendix B, one can prove that
(5.41) has a $C^{\infty}$ solution $\tau_{ij}(0\leq i,j\leq 2n-3)$
in $\O_+$ which satisfy
$$\tau_{ij}=O_0^{\infty}(1). $$

Combining this  with (5.40)-(5.41) shows (5.38). Therefore, Lemma 5.4
is proved. \qquad\qquad\qquad\qquad\qquad \qed

\medskip

In the following, we denote
$$\Cal S=\{S_i:0\le i \le 2n-3\}, \quad where\quad S_i= \t Z_i +\sum_{j=0}^{2n-3}
\tau_{ij} \t Z_j,\tag 5.42$$
which are called as the {\bf modified Klainerman's vector fields}.

With (5.36)-(5.37) and (5.39), we have
$$S_i-Z_i=O_1^{\infty}(1)\na,\quad i=0,1,\cdots, 2n-3.\tag 5.43$$

Set
$$\Cal S_{\Gamma}=\{S_{i\Gamma}: 0\leq i\leq 2n-3\},\quad where\quad
S_{0\Gamma}=t\p_t+t\p_t\chi\p_r,\quad
S_{i\Gamma}=S_i+S_i\chi\cdot\p_r(1\leq i\leq 2n-3),\tag 5.44$$ which
are tangential to $\Gamma$.

\medskip

{\bf Remark 5.3.} {\it From \rom{Lemma 5.4}, we can also derive that
for $m\in\Bbb N\cup\{0\}$
$$\Cal B_{\sigma}S^m\dot\varphi=0\quad on \quad r=\sigma(t,\o). \tag 5.45$$}
\medskip

{\bf Lemma 5.5.}  {\it Let $\dot \varphi$ be a $C^{k_0}(D_{T_0})$
solution to \rom{(4.4)}, where $k_0\in\Bbb N$, and
$$\ds\sum_{0\leq l\leq [\f{k_0}{2}]+1}|\na S^l\dot\vp|\leq M\ve,$$
where $M>0$ is some constant and $\ve$ is sufficiently small. Then
$$\dsize C(b_0, \gamma,k_0)\dsize\sum_{0\le l\le k_0-1}|\na_x
S^l\dot\varphi|\leq\sum_{0\le l\le k_0-1}t^l|\na_x^{l+1}\dot\varphi|
\le C(b_0, \gamma,k_0)\bigl(\dsize\sum_{0\le l\le k_0-1}|\na_x
S^l\dot\varphi|+\ve\bigr)\quad \text{\rom{in} $D_{T_0}$.}\tag5.46$$}

{\bf Proof.} It only comes from a direct computation (or one can see
Lemma 5.1 of [19]), we omit it here.

\qquad\qquad\qquad\qquad\qquad\qquad\qquad\qquad
\qquad\qquad\qquad\qquad\qquad\qquad\qquad\qquad
\qquad\qquad\qquad\qquad\qquad\qed

\medskip

Based on (5.45), (4.13) and Remark 5.3, we can apply Theorem 5.1 for $S^m\dot\vp
(0\le m\le k_0-1)$ directly to yield\medskip

{\bf Lemma 5.6.} {\it Under the assumptions of \rom{Theorem 5.1}, if
$\dot\vp$ is a $C^{k_0}(\overline{D_{T_0}})$ solution ($k_0\geq 9$)
of the problem $(4.4)$ with $(4.9)$ and $(4.11)$-$(4.13)$, then for
$0\leq m\leq k_0-1$ and $\mu<-1-\ds\f{1}{2}\sqrt{\ds\f{\g+7}{2}}$
$$ \align &C_1  T_0^{\mu+1}\int_{\sigma(T_0, \omega)\le
r\le\zeta(T_0,\o)}|\na S^{m}\dot\vp|^2(T_0,x)
dx+C_2\iint_{D_{T_0}}t^{\mu}
|\na S^{m}\dot\vp|^2 dtdx\\
&+C_3\int_{\Gamma_{T_0}}t^{\mu+1}(\p_t S^{m}\dot\vp)^2
dS+C_4\ds\int_{\Gamma_{T_0}}t^{\mu+1}\f{1}{r^2}\ds\sum_{k=1}^{3}(Z_k
S^m\dot\vp)^2
dS\\
\le& C(b_0,\g)\ve^2+C_5\int_{\Gamma_{T_0}}t^{\mu+1}(\Cal B_0
S^{m}\dot\vp)^2dS+\iint_{D_{T_0}} I_{2m} dtdx,\tag 5.47\endalign
$$
where the constant $C_i$ ($1\le i\le 5$) is given in \rom{Theorem
5.1} and
$$I_{2m}=\bigl|[S^m,\Cal L+\Cal P]\dot\vp-S^m R_0\bigr|\cdot\bigl|\Cal
M S^m\dot\vp\bigr|.\tag 5.48$$}\medskip

As in Theorem 5.3, we should control the terms
$\displaystyle\int_{\Gamma_{T_0}}t^{\mu+1}(\Cal B_0 S^m\dot\vp)^2dS$
and $\ds\iint_{D_{T_0}} I_{2m} dtdx$ in (5.47) to obtain the related
higher-order weighted energy estimates. To this end, as in Lemma
5.2, we can establish the following improved Hardy-type inequality
on $\ds\int_{\Gamma_{T_0}}t^{\mu-1}|S^m\dot\vp|^2dS$ due to
(4.10)-(4.12).\medskip

{\bf Lemma 5.7.}  {\it Assume that $\dot\vp\in
C^{k_0}(\overline{D_{T_0}})$ and $\xi(t,\o)\in C^{k_0}([1,T_0]\times
\Bbb S^2)$ with $k_0\geq 9$ are the solution of $(4.4)$ with $(4.9)$
and $(4.11)$-$(4.13)$,
and further assume
$$\ds\sum_{0\leq l\leq [\f{k_0}{2}]+1}|S^l\xi|+\ds\sum_{0\leq l\leq
[\f{k_0}{2}]+1}|\na S^l\dot\vp|\leq M\ve.\tag 5.49$$ Then
$$\align
&\int_{\Gamma_{T_0}}t^{\mu-1}|S^m\dot\vp|^2 dS\\
\leq& C(b_0,
\gamma)\ve^2+\f{1}{\mu^2}(1+\a)\int_{\Gamma_{T_0}}t^{\mu+1}(\p_t
S^m\dot\vp)^2 dS\\
&\ds +C(b_0, \gamma)\ve\biggl(\int_{\Gamma_{T_0}}t^{\mu+1}(\Cal B_0
S^m\dot\vp)^2
dS+\ds\int_{\Gamma_{T_0}}t^{\mu+1}\f{1}{r^2}\ds\sum_{k=1}^{3}(Z_k S^m\dot\vp)^2 dS\biggr)\\
&+C(b_0, \gamma)\biggl(\ds\sum_{0\leq l\leq
m-2}\int_{\Gamma_{T_0}}t^{\mu+1}(\na S^{l}\dot\vp)^2
dS+\int_{\Gamma_{T_0}}t^{\mu+1}|\xi|^2 dS\biggr).\tag 5.50
\endalign$$

Moreover,
$$\align
\int_{\Gamma_{T_0}}t^{\mu+1}&(\Cal B_0 S^m\dot\vp)^2 dS \leq  C(b_0,
\gamma)\ve^2+\f{1}{\mu^2 b_0^2}(1+\a)\int_{\Gamma_{T_0}}(\p_t
S^m\dot\vp)^2 dS\\
&+C(b_0, \gamma)\ve\int_{\Gamma_{T_0}}t^{\mu+1}\f{1}{r^2}\ds\sum_{k=1}^{3}(Z_k S^m\dot\vp)^2 dS\\
&+C(b_0, \gamma)\biggl(\ds\sum_{0\leq l\leq
m-1}\int_{\Gamma_{T_0}}t^{\mu+1}(\na S^{l}\dot\vp)^2
dS+\int_{\Gamma_{T_0}}t^{\mu+1}|\xi|^2 dS\biggr).\tag 5.51
\endalign$$}

{\bf Proof.} Since the proof procedure is very similar to Lemma 5.2, we just
give some necessary descriptions here. It follows from (4.1) and
(4.12) that
$$\cases
\ds\Cal B_0 S^m\dot\vp+\mu_2
S_\Gamma^m\xi=\f{1}{B_1}S_{\Gamma}^m\left(\kappa(\xi,\na\dot\vp)+R_1(t,x)\right)+[\Cal
B_0, S_{\Gamma}^m]\dot\vp+\Cal
B_0(S^m-S_{\Gamma}^m)\dot\vp,\\
&\\
\ds S^m\dot\vp=\mu_3 t\cdot S^m_{\Gamma}\xi+\bigl[S^m_{\Gamma},
\mu_3 t\bigr]\xi+\bigl(S^m-S_{\Gamma}^m\bigr)\dot\vp\endcases$$ with
$\bigl[\cdot, \cdot \bigr]$ being the commutator.

Then due to Lemma 5.5, the properties of $\kappa(\xi,\na\dot\vp)$,
Lemma 4.2, the assumption (5.49) and (5.44)
$$\cases
\ds\Cal B_0 S^m\dot\vp+\mu_2 S^m_\Gamma\xi=\ds\sum_{0\leq l\leq m-1}O(b_0)\na S^l\dot\vp
+O(\ve)\na S^m\dot\vp+O_{-1}^{0}(\ve),\\
\ds S^m\dot\vp-\mu_3\cdot t \cdot
S^m_\Gamma\xi=O(b_0)\xi+\ds\sum_{0\leq l\leq m-2}O(b_0)\na
S^l\dot\vp+O(\ve)\na S^{m-1}\dot\vp+O_{-1}^{0}(\ve).
\endcases\tag 5.52
$$

Since the left sides of the two equalities in (5.52) admit the same
forms as in (4.11) and (4.12), similar to Lemma 5.2, one has (5.50).
Combining (5.50) with (5.52) and Lemma 5.5 yields (5.51).\qquad\qquad\qquad\qquad\qed

\medskip

{\bf Theorem 5.8.} {\it Assume that $\dot\vp\in
C^{k_0}(\overline{D_{T_0}})$ and $\xi(t,\o)\in C^{k_0}([1,
T_0]\times \Bbb S^2)$ with $k_0\ge 9$ are the solution of
\rom{(4.4)} with \rom{(4.9)} and \rom{(4.11)-(4.13)}, and assume
$$\ds\sum_{0\leq
l\leq [\f{k_0}{2}]+1}|S^l\xi|+\ds\sum_{0\leq l\leq
[\f{k_0}{2}]+1}|\na S^{l}\dot\vp|\leq M\ve.\tag 5.53$$ Then for
sufficiently small $\ve>0$ and
$\mu<-1-\ds\f{1}{2}\sqrt{\ds\f{\g+7}{2}}$,
$$
\align &\int_{\sigma(T_0, \omega)\le
r\le\zeta(T_0,\o)}\dsize\sum_{0\le l\le k_0-1}T_0^{2l+\mu+1}
|\na^{l+1}\dot\vp|^2(T_0,x)dx+\iint_{D_{T_0}}\dsize\sum_{0\le l
\le k_0-1}t^{2l+\mu}|\na^{l+1}\dot\vp|^2dtdx\\
&+\int_{\Gamma_{T_0}}\dsize\sum_{0\le l\le k_0-1}t^{2l+\mu+1}
|\na^{l+1}\dot\vp|^2dS\le C_0 \ve^2.\tag 5.54
\endalign
$$}

{\bf Proof.} First, by Lemma 5.6-Lemma 5.7 and Theorem 5.3, we have
$$ \align &T_0^{\mu+1}\int_{\sigma(T_0, \omega)\le
r\le\zeta(T_0,\o)}|\na S^m\dot\vp|^2(T_0,x) dx
+\iint_{D_{T_0}}t^{\mu}|\na S^m\dot\vp|^2 dtdx
+\int_{\Gamma_{T_0}}t^{\mu+1}|\na S^m\dot\vp|^2dS\\
\leq& C_0\biggl(\ve^2+\iint_{D_{T_0}} I_{2m} dtdx+\ds\sum_{0\leq
l\leq m-1}\int_{\Gamma_{T_0}}t^{\mu+1}|\na S^{l}\dot\vp|^2
dS+\int_{\Gamma_{T_0}}t^{\mu+1}|\xi|^2 dS\biggr),\tag
5.55\endalign$$

To prove Theorem 5.8, we should estimate $\ds\iint_{D_{T_0}} I_{2m}
dtdx$ and $\ds\int_{\Gamma_{T_0}}t^{\mu+1}|\xi|^2 dS$ in (5.55).

Due to Lemma 4.5, then
$$[Z_0, \Cal L]=2\Cal L,\quad [Z_1, \Cal L]=0,\qquad  i=1,2,3.$$

By (5.43), Lemma 5.5 and induction method that
$$[S^m, \Cal L]\dot\vp=\ds\sum_{0\leq l\leq m-1}C_{lm}S^{l}\Cal
L\dot\vp+\f{1}{t}\ds\sum_{k=1}^{m}O(\ve)\na S^{k}\dot\vp.\tag 5.56$$

On the other hand,
$$[S^m,\Cal P]\dot\vp=\ds\sum_{i,j=1}^{2}\biggl[S^m, f_{ij}\p_{ij}\biggr]
\dot\vp+\ds\sum_{i}^{2}\biggl[S^m,
f_{ij}\p_{ti}\biggr]\dot\vp.$$

With the assumption (5.53), it follows from (4.7) that
$$f_{ij}=O_0^{k_0-1}(\ve)+O_0^{k_0-1}(1)\na\dot\vp,\qquad 0\leq i\leq 3; 1\leq j\leq 3,$$
thus with Lemma 5.5,
$$\bigl|[S^m,\Cal P]\dot\vp\bigr|\leq \f{1}{1+t}\ds\sum_{l=0}^{m}O(\ve)|\na
S^m\dot\vp|.$$

Combining this with (5.18), (5.48), (5.56) and Lemma 4.2,
$$I_{2m}\leq t^{\mu}\biggl(O(\ve)|\na S^m\dot\vp|^2+\ds\sum_{0\leq l\leq m-1}\t
C_{lm}|\na S^l\dot\vp|^2+\f{O(\ve^2)}{(1+t)^2}\biggr).$$

This shows that
$$\iint_{D_{T_0}} I_{2m} dtdx\leq
C(b_0,\g)\biggl(\ve^2+\ve\iint_{D_{T_0}}t^{\mu}|\na
S^m\dot\vp|^2+\ds\sum_{l=0}^{m-1}\iint_{D_{T_0}}t^{\mu}|\na
S^m\dot\vp|^2 dtdx\biggr).\tag 5.57$$

Furthermore, by the boundary condition (4.12) and Lemma 5.2,
$$\int_{\Gamma_{T_0}}t^{\mu+1}|\xi|^2 dS\leq
C(b_0,\g)\int_{\Gamma_{T_0}}t^{\mu-1}|\dot\vp|^2 dS\leq
C(b_0,\g)\ve^2.\tag 5.58$$

Substituting (5.57) and (5.58) into (5.55) and using the induction
method and (5.46) in Lemma 5.5 yield (5.54). Thus the proof of
Theorem 5.8 is
completed.\qquad\qquad\qquad\qquad\qquad\qquad\qquad\qquad\qquad\qquad\qquad
\qquad \qed

\vskip 0.3 true cm

{\bf $\S 5.3.$ Proof of Theorem 1.1 for $n=3$}
\vskip 0.3 true cm

Based on the higher order energy estimate established in Theorem
5.8, we now prove the global existence of a shock wave in Theorem
1.1 by the local existence result in $\S 3$ and the continuous
induction method. For any given $t_0>0$, the solution of (1.6) with
the initial data given on $t=t_0$ and the boundary conditions
(1.7)-(1.10) in $[t_0, t_0+t^*]$ for some $t^*>0$ can be obtained by
the local existence of the solution in  $\S 3$ , provided that the
initial data are smooth and satisfy the compatibility conditions.
Moreover, if the perturbation of the initial data given on $t=t_0$
is small as $O(\ve)$, then the lifespan of the solution is at least
as large as $\ds\f{C}{\ve}$ with $C>0$. Therefore, as long as we can
establish that the maximum norm of $\dot\vp, \xi$ and their
derivatives decays with a rate in $t$, then the solution can be
extended continuously to the whole domain.  That is, by the local
existence result and the property of decay of the solution we can
obtain the uniform bound of $\dot\vp, \xi$ and their derivatives,
and then extend the solution continuously from $t_0\le t\le
t_0+\eta^*$ to $t_0+\eta^*\le t\le t_0+2\eta^*$ with $\eta^*>0$
being independent of $t_0$. Hence the key point to prove Theorem 1.1
is to give the decay of the maximum norm of $\dot \vp, \xi$ and
their derivatives.

To finish the proof of Theorem 1.1,  the following
Lemma is required.\medskip

{\bf Lemma 5.9.}  {\it Under the assumption $(5.53)$ in \rom{Theorem 5.8},
for $1\leq t\leq T_0$ we have
$$\dsize\sum_{0\le l\le k_0-4}|t^l\na^{l+1}\dot\vp|^2
\le C_0 t^{-3}\int_{\sigma(t, \omega)\le r\le\zeta(t,\o)}
\dsize\sum_{0\le l\le k_0-1} |t^l\na^{l+1}\dot\vp(t,r,\o)|^2 dx.\tag
5.59$$}

{\bf Proof.} We will apply  Sobolev's imbedding
theorem to establish (5.59).

For any $t_1\in [1,T_0]$,  set
$$(t',x')=\f{1}{t_1}(t,x).$$

Then one has
$$\na_{t,x}^k\dot\vp=\f{1}{t_1^k}\na_{t',x'}^k\dot\vp, \quad \forall\
k\in\Bbb N.\tag 5.60$$

Define $D_{*}=\{(t',x'): t'=1, \f{\si(t,\o)}{t}\leq |x'|\leq
\f{\zeta(t,\o)}{t}\}$, then by Sobolev's imbedding theorem in space dimensions 3 (since $D_{*}$
has the uniform interior cone condition),
$$|\na_{t',x'}\dot\vp|^2(1,x')|\leq C\int_{D_{*}}\ds\sum_{0\leq l\leq
3}|\na_{t',x'}^{l+1}\dot\vp|^2(1,x')dx'.$$

With (5.60),
$$\align
|\na_{t,x}\dot\vp|^2(t_1,x)=&\f{1}{t_1^2}|\na_{t',x'}\dot\vp|^2(1,x')\\
\leq &\f{C}{t_1^2}\int_{D_{*}}\ds\sum_{0\leq l\leq
3}|\na_{t',x'}^{l+1}\dot\vp|^2(1,x')dx'\\
=&\f{C}{t_1^2}\int_{\si(t,\o)\leq r\leq \zeta(t,\o)}\ds\sum_{0\leq
l\leq 3}|t_1^{l+1}\na_{t,x}^{l+1}\dot\vp|^2(t_1,x) \f{1}{t_1^3} dx\\
=&\f{C}{t_1^3}\int_{\si(t,\o)\leq r\leq \zeta(t,\o)}\ds\sum_{0\leq
l\leq 3}|t_1^{l}\na_{t,x}^{l+1}\dot\vp|^2(t_1,x)  dx.\endalign$$

This yields (5.59) for $l=0$.

In the same way, we can finish the proof of Lemma 5.9 in the case of $1\le l\le k_0-4$.
\qquad\qquad\qquad\qquad \qed\medskip

It follows from (5.54) that
$$\int_{\sigma(t, \omega)\le r\le\zeta(t,\o)}\dsize\sum_{0\le l\le k_0-1}
|t^l\na^{l+1}\dot\vp(t,x)|^2 dx\le C_0 \ve^2 t^{-\mu-1}.$$

Combining this with (5.59) yields  $\dsize\sum_{0\le l\le
k_0-4}|t^l\na^{l+1}\dot\vp|^2 \le C_0 \ve^2 t^{-\mu-4}$ for
$\sigma(t,\o)\le r\le\chi(t,\o)$ and $1\le t\le T$. For $k_0\ge 9$,
$k_0-4 \ge [\f {k_0}2]+1$, so one has
$$\dsize\sum_{0\le \ell\le [\f{k_0}{2}]+1} |t^l\na^{l+1}\dot\vp|\le
C_0 \ve{{t^{-\f{\mu}{2}-2}}}.\tag 5.61$$

In addition, the equations (4.11) and (4.12) yield $$\ds\sum_{0\leq
l\leq [\f{k_0}{2}]+1}|S^l\xi|\leq \f{M}{2}\ve,\tag 5.62$$ when
$\ds-\f{\mu}{2}-2<0$.

When we choose $\mu\in (-4,
-1-\ds\f{1}{2}\sqrt{\ds\f{\g+7}{2}})$, (5.61)-(5.62) shows the
induction assumption (5.53) in Theorem 5.8, then the proof of
Theorem 4.1 and furthermore Theorem 1.1 can be completed for $n=3$.
\qquad\qquad \qed

\vskip 0.3 true cm \centerline{\bf $\S 6$. Sketch on the proof of
Theorem 1.1 for $n=2$} \vskip 0.3 true cm

In this section, we only give the sketch on the proof of Theorem 1.1
for $n=2$ since the main procedure is same as the case $n=3$ and even simpler. In
addition, some related notations in this section admit the same meanings as in 3-D case.\medskip

{\bf Theorem 6.1.} {\it For n=2, assume that $\dot\vp\in
C^{k_0}(\overline{D_{T_0}})$ and $\xi(t,\o)\in C^{k_0}([1,
T_0]\times [0, 2\pi])$ with $k_0\ge 7$ are the solution of $(4.4)$ with
$(4.9)$ and $(4.11)$-$(4.13)$, and further assume
$$\ds\sum_{0\leq
l\leq [\f{k_0}{2}]+1}|S^l\xi|+\ds\sum_{0\leq l\leq
[\f{k_0}{2}]+1}|\na S^{l}\dot\vp|\leq M\ve.\tag 6.1$$ Then for
sufficiently small $\ve>0$ and
$\mu<-\ds\f{1}{2}-\ds\f{1}{2}\sqrt{\ds\f{\g+1}{2}}$,
$$
\align &\int_{\sigma(T_0, \omega)\le
r\le\zeta(T_0,\o)}\dsize\sum_{0\le l\le k_0-1}T_0^{2l+\mu+1}
|\na^{l+1}\dot\vp|^2(T_0,x)dx+\iint_{D_{T_0}}\dsize\sum_{0\le l
\le k_0-1}t^{2l+\mu}|\na^{l+1}\dot\vp|^2dtdx\\
&+\int_{\Gamma_{T_0}}\dsize\sum_{0\le l\le k_0-1}t^{2l+\mu+1}
|\na^{l+1}\dot\vp|^2dS\le C(b_0, \gamma)\ve^2.\tag 6.2
\endalign
$$}

\medskip

{\bf Sketch of the proof.}\medskip

We will divide the proof procedure into the following two steps.\medskip

\centerline{\bf Step 1. Establishing a priori estimate containing a shock boundary condition}
\medskip

At first, as in Theorem 5.1, we look for such an operator $\Cal M$
$$\align
\Cal M&=\Cal A(t,x)\p_t+\Cal
B(t,x)\p_r\\
&=t^{\mu}r\p_t+t^{\mu+1}b_\si(t,r,\o)\p_r,\tag 6.3\endalign$$ where
$$b_\si(t,r,\o)=s^2\biggl(1+\f{\ve\p_t b-\f{\ve}{t}b}{s}
+\f{\ve}{b_0}(e_1-b_0-\f{\ve}{t}b)\biggr)\tag 6.4$$
with the constant $e_1$ being determined later on.

For $0\leq m\leq k_0-1$, it follows from  Lemma 5.4 and Remarks 5.3
that $\Cal B_\si S^m\dot\vp=0$ holds on $r=\si(t,\o)$. In this situation,
we can claim that for any fixed constant
$-3<\mu<-\ds\f{1}{2}-\ds\f{1}{2}\sqrt{\ds\f{\g+1}{2}}$,
$$ \align &C_1  T_0^{\mu+1}\int_{\sigma(T_0, \omega)\le
r\le\zeta(T_0, \o)}|\na S^{m}\dot\vp|^2(T_0,x)
dx+C_2\iint_{D_{T_0}}t^{\mu}
|\na S^{m}\dot\vp|^2 dtdx\\
&+C_3\int_{\Gamma_{T_0}}t^{\mu+1}(\p_t S^{m}\dot\vp)^2
dS+C_4\ds\int_{\Gamma_{T_0}}t^{\mu+1}\f{1}{r^2}(Z_1 S^m\dot\vp)^2
dS\\
& \le C(b_0,\g)\ve^2+C_5\int_{\Gamma_{T_0}}t^{\mu+1}(\Cal B_0
S^{m}\dot\vp)^2dS+\iint_{D_{T_0}} I_{3m} dtdx,\tag 6.5\endalign
$$
where $C_i$ ($1\le i\le 5$) are just only the constants given in Theorem 5.1, and
$$I_{3m}=\bigl|[S^m,\Cal L+\Cal P]\dot\vp-S^m
R_0\bigr|\cdot\bigl|\Cal M S^m\dot\vp\bigr|.\tag 6.6$$

Indeed, it follows from the integration by parts that
$$\align
&\iint_{D_{T_0}}\bigl(S^m R_0-[S^m,\Cal L+\Cal P]\dot\vp\bigr)\cdot
\Cal M
S^m\dot\vp dtdx\\
=&\iint_{D_{T_0}}(\Cal L+\Cal P)S^m\dot\vp\cdot \Cal M S^m\cdot\dot\vp dtdx\\
=&\iint_{D_{T_0}} I_{4m} dtdx\\
&+\int_{\si(T_0,\o)<r<\zeta(T_0,\o)}N_0(S^m\dot\vp)(T_0,x)
dx-\int_{\si(1,\o)<r<\zeta(1,\o)}N_0(S^m\dot\vp)(1,x) dx\\
&+\int_{\Gamma_{T_0}}\biggl(\ds\sum_{i=1}^{2}(\f{x_i}{r}N_i(S^m\dot\vp)-\p_i\chi\cdot
N_i(S^m\dot\vp))-\p_t\chi N_0(S^m\dot\vp)\biggr)dS\\
&+\int_{B_{T_0}}\biggl(\p_t\si
N_0(S^m\dot\vp)-\ds\sum_{i=1}^{2}(\f{x_i}{r}N_i(S^m\dot\vp)-\p_i\si\cdot
N_i(S^m\dot\vp))\biggr) dS\tag 6.7
\endalign$$
with
$$\align
I_{4m}&=K_{00}(\p_t S^m\dot\vp)^2+K_{0r}\p_t S^m\dot\vp\p_r
S^m\dot\vp+K_{rr}(\p_r S^m\dot\vp)^2+\f{K_{22}}{r^2}(Z_1 S^m\dot\vp)^2\\
&+K_{01}Z_1\Cal A\cdot \p_t S^m\dot\vp Z_1 S^m\dot\vp+K_{r1}Z_1\Cal
B\cdot \p_r S^m\dot\vp Z_1 S^m\dot\vp\\
&+K^{00}(\p_t S^m\dot\vp)^2+\ds\sum_{i=1}^{2}K^{0i}\p_t
S^m\dot\vp\p_i S^m\dot\vp+\ds\sum_{i,j=1}^{2}K^{ij}\p_i S^m\dot\vp
\p_j S^m\dot\vp.\tag 6.8
\endalign$$

Similar to Step 3 in the proof of Theorem 5.1, under the assumption
(6.1), we can derive from (6.8) that
$$I_{4m}=K_{00}(\p_t S^m\dot\vp)^2+K_{0r}\p_t S^m\dot\vp \p_r
S^m\dot\vp+K_{rr}(\p_r S^m\dot\vp)^2+\f{K_{22}}{r^2}(Z_1
S^m\dot\vp)^2+t^{\mu}O(\ve)|\na S^m\dot\vp|^2$$ with
$$\cases
&K_{00}=t^{\mu}\biggl(\ds\f{1}{2}(1+e_1-\mu)b_0(1+\a)\biggr),\\
&K_{0r}=t^{\mu}\biggl(\bigl(\ds\f{\g+1}{2}+e_1-\mu\bigr)b_0^2(1+\a)\biggr),\\
&K_{rr}=t^{\mu}\biggl(\ds\f{\g+1}{4}(1+e_1-\mu)b_0^3(1+\a)\biggr),\\
&K_{22}=t^{\mu}\biggl(-\ds\f{\g-1}{4}(1+e_1+\mu)(1+\a)\biggr),\\
&K_{0r}^2-4K_{00}
K_{rr}=t^{2\mu}\biggl(\ds\f{\g-1}{4}\bigl(\g+1-2(e_1-\mu)^2\bigr)(1+\a)\biggr).
\endcases$$

The sufficient conditions for $I_{4m}\ge 0$ are
$$1+e_1-\mu>0,\quad 1+e_1+\mu<0,\quad \g+1-2(e_1-\mu)^2<0.\tag 6.9$$

Consequently, if we select
$$\mu<-\f{1}{2}-\f{1}{2}\sqrt{\f{\gamma+1}{2}}\qquad and \qquad e_1=\f{1}{2}
\sqrt{\f{\gamma+1}{2}}-\f{1}{2},\tag 6.10$$
then (6.9) stands. This implies
$$K_{00}>0,\quad K_{0r}^2-4K_{00} K_{rr}<0, \quad K_{22}>0,$$
and then
$$\iint_{D_{T_0}} I_{4m} dtdx\geq
C_2\iint_{D_{T_0}}t^{\mu}|\na\dot\vp|^2 dtdx.\tag 6.11$$

Moreover, similar to the proof of Theorem 5.1, under the constrain
(6.10), the definitions (6.3)-(6.4) and (4.13), we have
$$\cases
&\ds\int_{\Gamma_{T_0}}\biggl(\ds\sum_{k=1}^{2}(\f{x_i}{r} N_i(\na
S^m\dot\vp)-\p_i\chi\cdot N_i(\na S^m\dot\vp)-\p_t\chi N_0(\na
S^m\dot\vp)\biggr) dS\\
&\quad\ds\geq C_3\iint_{\Gamma_{T_0}}t^{\mu+1}(\p_t S^m\dot\vp)^2 dS
+C_4\int_{\Gamma_{T_0}}t^{\mu+1}\f{1}{r^2}(Z_1 S^m\dot\vp)^2 dS
 -C_5\int_{\Gamma_{T_0}}t^{\mu+1}(\Cal B_0 S^m\dot\vp)^2 dS,\\
&\ds\int_{B_{T_0}}\biggl(\p_t\si N_0(S^m\dot\vp)-\ds\sum_{k=1}^{2}(\f{x_i}{r}
N_i(\na S^m\dot\vp)-\p_i\si N_i(\na S^m\dot\vp)\biggr)dS=0,\\
&\ds\int_{\sigma(T_0, \omega)\le r \le\zeta(T_0,\o)}N_1(\na S^m\dot\vp)(T_0,x)dx
\geq C_1\iint_{\sigma(T_0, \omega)\le r \le \zeta(T_0,\o)}|\na S^m\dot\vp|^2(T_0,x)dx,\\
&\ds\int_{\sigma(1, \omega)\le r\le\zeta(1,\o)}N_1(\na
S^m\dot\vp)(1,x)dx\leq C(b_0,\g)\ve^2.\endcases\tag 6.12$$

Substituting (6.11)-(6.12) into (6.7) and subsequently combining with (6.6) yield (6.5).

\vskip 0.3 true cm

\centerline{\bf Step 2. Establishing a priori estimate on the solution}

\vskip 0.3 true cm

In order to obtain (6.2), we should deal with the term
$\ds\int_{\Gamma_{T_0}}t^{\mu+1}(\Cal B_0 S^m\dot\vp)^2 dS$ in
(6.5).

Due to $\mu_2<0$ and $\ds \mu_3(t,x)<0$ in (4.11) and (4.12), we can
obtain the following improved Hardy-type inequality
$$\align
\ds\int_{\Gamma_{T_0}}t^{\mu-1}(S^m\dot\vp)^2 dS&\leq C(b_0,
\gamma)\ve^2+\f{1}{\mu^2}(1+\a)\int_{\Gamma_{T_0}}t^{\mu+1}(\p_t
S^m\dot\vp)^2
dS\\
&\ds+C(b_0, \gamma)\ve\biggl(\int_{\Gamma_{T_0}}t^{\mu+1}(\Cal B_0
S^m\dot\vp)^2 dS+\int_{\Gamma_{T_0}}t^{\mu+1}\f{1}{r^2}(Z_1
S^m\dot\vp)^2 dS\biggr)\\
&+C(b_0, \gamma)\Big(\ds\sum_{0\leq l\leq
m-1}\int_{\Gamma_{T_0}}t^{\mu-1}(\na S^l\dot\vp)^2 dS +
\int_{\Gamma_{T_0}} t^{\mu+1} |\xi|^2 dS\Big).\tag 6.13
\endalign$$

So with (6.13), it follows from (5.52) and Lemma 4.4 that
$$\align
\ds\int_{\Gamma_{T_0}}t^{\mu+1}(\Cal B_0 S^m\dot\vp)^2 dS \leq&
C(b_0, \gamma)\ve^2+\f{1}{4\mu^2
b_0^2}(1+\a)\int_{\Gamma_{T_0}}t^{\mu+1}(\p_t S^m\dot\vp)^2 dS\\
&+\ds C(b_0, \gamma)\ve\int_{\Gamma_{T_0}}t^{\mu+1}\f{1}{r^2}(Z_1
S^m\dot\vp)^2 dS\\
&+C(b_0, \gamma)\Big(\ds\sum_{0\leq l\leq
m}\int_{\Gamma_{T_0}}t^{\mu-1}(\na S^l\dot\vp)^2 dS +
\int_{\Gamma_{T_0}} t^{\mu+1} |\xi|^2 dS\Big).\tag 6.14
\endalign$$

Substituting (6.13)-(6.14) and (4.12) into (6.5) yields
$$ \align &T_0^{\mu+1}\int_{\sigma(T_0, \omega)\le
r\le\zeta(T_0,\o)}|\na S^m\dot\vp|^2(T_0,x) dx
+\iint_{D_{T_0}}t^{\mu}|\na S^m\dot\vp|^2 dtdx\\
&+\f{\g-1}{8}b_0^2\biggl(1-\f{1}{\mu^2}\biggr)(1+\a)\int_{\Gamma_{T_0}}t^{\mu+1}
|\p_t S^m\dot\vp|^2dS\\
&+\int_{\Gamma_{T_0}}t^{\mu+1}\f{1}{r^2}(Z_1 S^m\dot\vp)^2 dS\\
&\leq C(b_0,\g)\biggl(\ve^2+\iint_{D_{T_0}}I_{3m}dtdx+\ds\sum_{0\leq
l\leq m-1}\int_{\Gamma_{T_0}}t^{\mu+1}(\na S^l\dot\vp)^2
dS\biggr).\tag 6.15\endalign$$

Next, we deal with $\ds\iint_{D_{T_0}} I_{3m}dtdx$
in (6.15).

Due to Lemma 4.4,
$$[Z_0, \Cal L]=2\Cal L,\quad [Z_1, \Cal L]=0 .$$

By (5.43), Lemma 5.5 and induction method, one has
$$[S^m, \Cal L]\dot\vp=\ds\sum_{0\leq l\leq m-1}C_{lm}S^{l}\Cal
L\dot\vp+\f{1}{t}\ds\sum_{k=1}^{m}O(\ve)\na S^{k}\dot\vp.\tag 6.16$$

On the other hand,
$$[S^m,\Cal P]\dot\vp=\ds\sum_{i,j=1}^{2}\biggl[S^m, f_{ij}\p_{ij}\biggr]
\dot\vp+\ds\sum_{i}^{2}\biggl[S^m,
f_{ij}\p_{ti}\biggr]\dot\vp.$$

With the assumption (5.53), it follows from (4.7) that
$$f_{ij}=O_0^{k_0-1}(\ve)+O_0^{k_0-1}(1)\na\dot\vp\quad (0\leq i\leq 2; 1\leq j\leq 2),$$
thus with Lemma 5.5,
$$\bigl|[S^m,\Cal P]\dot\vp\bigr|\leq \f{1}{1+t}\ds\sum_{l=0}^{m}O(\ve)|\na
S^m\dot\vp|.$$

Combining this with (5.56), (5.18), (6.16) and Lemma 4.2 yields
$$I_{3m}\leq t^{\mu}\biggl(O(\ve)|\na S^m\dot\vp|^2+\ds\sum_{0\leq l\leq m-1}\t
C_{lm}|\na S^l\dot\vp|^2+\f{O(\ve^2)}{(1+t)^2}\biggr).$$

This shows that
$$\iint_{D_{T_0}} I_{3m} dtdx\leq
C(b_0,\g)\biggl(\ve^2+\ve\iint_{D_{T_0}}t^{\mu}|\na
S^m\dot\vp|^2+\ds\sum_{l=0}^{m-1}\iint_{D_{T_0}}t^{\mu}|\na
S^m\dot\vp|^2 dtdx\biggr).$$

Substituting this into (6.15) yields (6.2) and then the proof of Theorem
6.1 is completed.\qquad\qquad\qquad\qquad\qed\medskip

Based on Theorem 6.1, we can prove Theorem 4.1 and Theorem 1.1 in
the case of $n=2$ by the local existence result and the continuous
induction method. Similar to the case of $n=3$, we just need to verify
the induction assumption (6.2). Similar to Lemma 5.9, it follows
from the Sobolev's imbedding theorem and the
assumptions of Theorem 6.1 that for $\sigma(t, \omega)\le
r\le\chi(t,\o)$ and $1\le t\le T_0$, one has
$$\dsize\sum_{0\le l\le k_0-3}|t^l\na^{l+1}\dot\vp|^2
\le Ct^{-2}\int_{\sigma(t, \omega)\le r\le\chi(t,\o)}
\dsize\sum_{0\le l\le k_0-1} |t^l\na^{l+1}\dot\vp(t,r,\o)|^2 dx.$$

On the other hand, (6.2) shows that
$$\int_{\sigma(t, \omega)\le r\le\chi(t,\o)}\dsize\sum_{0\le l\le k_0-1}
|t^l\na^{l+1}\dot\vp(t,x)|^2 dx\le C_0\ve^2 t^{-\mu-1}.$$

Hence $\dsize\sum_{0\le l\le k_0-3}|t^l\na^{l+1}\dot\vp|^2 \le
C_0\ve^2 t^{-\mu-3}$ for $\sigma(t, \omega)\le r\le\chi(t,\o)$ and
$1\le t\le T_0$. For $k_0\ge 7$, then $k_0-3\geq [\f{k_0}{2}]+1$, so one
has $$\dsize\sum_{l\le [\f{k_0}{2}]+1} |t^l\na^{l+1}\dot\vp|\le
C(b_0)\ve{{t^{-\f{\mu}{2}-\f32}}}.\tag 6.17$$

In addition, the equations (4.11) and (4.12) yield
$$\ds\sum_{0\leq l\leq [\f{k_0}{2}]+1}|S^l\xi|\leq \f{M}{2}\ve,\tag
6.18$$ when $-\mu-3<0$. When we choose $\mu\in (-3,
-\f{1}{2}-\f{1}{2}\sqrt{\f{\gamma+1}{2}})$, (6.17)-(6.18) shows the
assumption 6.1 in Theorem 6.1, then the proof of Theorem 4.1 and
meanwhile Theorem 1.1 can be completed for $n=2$.

\qquad\qquad\qquad\qquad\qquad\qquad\qquad\qquad\qquad\qquad\qquad\qquad\qquad\qquad
\qquad\qquad\qquad\qquad\qquad\qquad
\qquad\qed\medskip

\newpage \centerline{\bf Appendix A. Some basic computations} \vskip 0.3 true cm

In this Appendix, at first, we will give  the detailed derivations on the coefficients
in the equation (3.10) under
the transformations (3.5) and (3.6). In terms of the expressions
of $A_{k,0}, A_{k,1}, ..., A_{6,2}^{ij}$,
which are  given in (3.10), we have

{\bf Lemma A.1.}
{\it
$$\cases
A_{1,0}&=\psi,\\
A_{2,0}&=(R-2)\p_T b-a_1\bigl[(R-1)a_3+\psi\p_T a_0\bigr],\\
A_{3,0}^i&=0, \qquad i=1,2,3,\\
A_{4,0}&=\p_T
a_0\cdot a_1\biggl((R-1) a_1\cdot a_3-(R-2)\p_T b\biggr),\\
A_{5,0}^i&=0, \qquad i=1,2,3,\\
A_{6,0}^{ij}&=0, \qquad i,j=1,2,3\\
A_{7,0}&=(\p_T\psi)^2+\psi\p_T^2 b+\p_T\psi\cdot\p_T b+(R-2)\p_T^2
b\cdot\p_R\psi\\
&\quad -a_1\bigl(\p_T\psi\cdot a_3+\p_T
a_0(\p_T\psi\cdot\p_R\psi+2\p_R\psi\cdot\p_T b)\bigr)+2\p_T
a_0\cdot a_1^2\cdot a_3\cdot \p_R\psi,
\endcases$$
where $a_0, a_1, a_2$ and $a_3$ are given in \rom{(3.8)}.}

{\bf Proof.} Under the transformations (3.5) and (3.6),  a direct computation yields
$$\align
\p_t^2\phi
&=-b_0(\p_T-a_1\cdot \p_T a_0\p_R)(a_1\cdot a_3)\\
&=b_0 a_1^2\cdot a_3\bigl(\p_T\psi+(R-1)\p_{TR}^2\psi\bigr)\\
&\quad -b_0 a_1\bigl((\p_T\psi)^2+\psi\p_T^2\psi+\p_T\psi\p_T b+\psi\p_T^2
b+(R-2)\p_T^2 b\p_R\psi+(R-2)\p_T b\p_{TR}^2\psi\bigr)\\
&\quad -b_0 a_1^3\cdot a_3\cdot \p_T a_0\bigl(2\p_R\psi+(R-1)\p_R^2\psi\bigr)\\
&\quad +b_0 a_1^2 \p_T
a_0\bigl(\p_T\psi\p_R\psi+\psi\p_{TR}^2\psi+2\p_R\psi\p_T
b+(R-2)\p_T b\p_R^2\psi\bigr).\tag A.1\\
\endalign
$$

On the other hand, by comparing the coefficients of power $t^0$ in
the equation of (3.2), we also have

$$
\align
\p_t^2\phi=&-b_0
a_1\bigl(A_{1,0}\p_T^2\psi+A_{2,0}\p_{TR}^2\psi+\sum_{i=1}^{3}A_{3,0}^{i}\p_T
Z_i\psi+A_{4,0}\p_R^2\psi+\sum_{i=1}^{3}A_{5,0}^i\p_R
Z_i\psi\\
&+\sum_{i=1}^{3}\sum_{j=1}^{3}A_{6,0}^{ij}Z_i Z_j\psi+A_{7,0}\bigr)
.\tag A.2
\endalign$$

Thus, by (A.1) and (A.2), we can complete the proof on Lemma A.1.\qquad\qquad\qquad\qquad
\qquad\qquad\qquad \qed\medskip

Analogously, a direct but tedious computation yields

{\bf Lemma A.2.}  {\it
$$\cases
A_{1,1}&=0,\\
A_{2,1}&=2\ds(b_0 a_1 a_2-a_0)\bigl(1-(R-1)a_1\p_R\psi\bigr)\\
&\quad \ds+\f{2b_0 a_1}{a_0^2}\biggl((R-1)a_1\sum_{i=1}^{3}(a_4^i)^2-(R-2)\sum_{i=1}^{3}
Z_i b\cdot a_4^i\biggr),\\
A_{3,1}^i&=\ds -\f{2b_0}{a_0^2}a_1\cdot a_4^i\cdot\psi, \qquad i=1,2,3,\\
A_{4,1}&=2\p_T a_0 a_1\biggl(b_0 a_1 a_2-a_0\biggr)\biggl((R-1)a_1\p_R\psi-1\biggr)\\
&\quad\ds+\f{2b_0}{a_0^2}\p_T a_0 a_1^2\biggl((R-2)\sum_{i=1}^{3}Z_i
b\cdot a_4^i-(R-1)a_1\sum_{i=1}^{3}(a_4^i)^2\biggr),\\
A_{5,1}^i&=\ds\f{2b_0}{a_0^2}\p_T a_0 a_1^2\psi a_4^i, \qquad i=1,2,3,\\
A_{6,1}^{ij}&=0, \qquad\qquad\qquad i,j=1,2,3,\\
A_{7,1}&=\ds 2a_1\biggl(\f{b_0}{a_0^2} a_1\sum_{i=1}^{3}(a_4^i)^2-\p_R\psi(b_0 a_1 a_2
-a_0)\biggr)(\p_T\psi-2a_1\p_T a_0\p_R\psi)+2a_3\\
&\quad\ds-\f{2b_0}{a_0^2}a_1\sum_{i=1}^{3}a_4^i\biggl(\p_T\psi
Z_i\psi+\p_T\psi Z_i b+\psi\p_T Z_i b+(R-2)\p_T Z_i b\p_R\psi\\
&\quad\ds-a_1\p_T a_0\bigl[\p_R\psi Z_i\psi+2\p_R\psi Z_i b\bigr]\biggr).
\endcases$$}

{\bf Lemma A.3.}  {\it
$$
A_{1,2}=0,\quad
A_{2,2}=0,\quad
A_{3,2}^i=0, \quad i=1,2,3,$$
and
$$\align
A_{4,2}&=a_1\biggl((b_0 a_1 a_2-a_0)^2-c^2(\rho)\biggr)(1-(R-1)a_1\p_R\psi)\\
&\quad\ds+\f{2b_0 a_1}{a_0^2}\bigl(a_0-b_0 a_1 a_2\bigr)\sum_{i=1}^{3}a_4^i\biggl((R-2)
a_1 Z_i b-(R-1)a_1^2 a_4^i\biggr)\\
&\quad+\ds\f{1}{a_0^2}\sum_{i=1}^{3}\sum_{j=1}^{3}\biggl(c^2(\rho)\delta_{ij}
-\f{(b_0 a_1)^2}{a_0^2}a_4^i a_4^j\biggr)\biggl((R-2)Z_j b-(R-1)a_1 a_4^j\biggr) a_1 Z_i a_0,\\
A_{5,2}^i&=\ds\f{2b_0 a_1^2}{a_0^2}\bigl(a_0-b_0 a_1 a_2\bigr) a_4^i\psi\\
&\quad\ds+\f{1}{a_0^2}\sum_{j=1}^{3}\biggl(c^2(\rho)\delta_{ij}
-\f{b_0^2 a_1^2}{a_0^2}a_4^i a_4^j\biggr)\biggl(a_1 Z_j a_0\vp+(R-1)a_1 a_4^j-(R-2)Z_j b\biggr),
\quad i=1,2,3,\\
A_{6,2}^{ij}&=\ds\f{1}{a_0^2}\biggl(\f{b_0^2 a_1^2}{a_0^2}a_4^i a_4^j-c^2(\rho)\delta_{ij}\biggr)\psi,
\qquad\qquad \quad i,j=1,2,3,\\
\endalign$$
$$\align
A_{7,2}&=\ds 2(a_1\p_R\psi)^2\biggl(c^2(\rho)-(b_0 a_1
a_2-a_0)^2\biggr)+\f{2a_2}{a_0}c^2(\rho)+\f{b_0 a_1}{a_0^3}(b_0 a_1 a_2-2a_0)\sum_{i=1}^{3}(a_4^i)^2\\
&\ds+\f{2b_0 a_1^2}{a_0^2}(a_0-b_0 a_1
a_2)\p_R\psi\sum_{i=1}^{3}a_4^i\biggl(Z_i\psi+2Z_i b-2a_1
a_4^i\biggr)\\
&\ds+\f{1}{a_0^2}\sum_{i=1}^{3}\sum_{j=1}^{3}\biggl(\f{b_0^2
a_1^2}{a_0^2}a_4^i
a_4^j-c^2(\rho)\delta_{ij}\biggr)\biggl(Z_i\psi
Z_j\psi+Z_i\psi Z_j b+\psi Z_i Z_j b+(R-2)Z_i Z_j b\p_R\psi\\
&\ds-a_1 Z_i a_0\p_R\psi\bigl(Z_j\psi+2Z_j b\bigr)-a_1
a_4^j\bigl(Z_i\psi-2 a_1 Z_i a_0\p_R\psi\bigr)\biggr).
\endalign$$}
Next, we estimate $E_k(d_0,\psi,\na\psi)$ $(k=0,1,2,\cdots K)$, which are defined in (3.21) and
(3.32) respectively.

{\bf Lemma A.4.} {\it If $\|\psi-\hat\psi\|_{C^2}\leq C\ve$, then for large $b_0$ and small $\ve$
$$\ds E_k(d_0,\psi,\na\psi)\geq
\f{(\g-1)b_0}{2}>0,\quad k=0,1,2,\cdots K.\tag A.3$$
}

{\bf Proof.}  It follows from the expressions of $E_k(d_0,\psi,\na\psi)$, Lemma 3.2, Lemma 2.1
and  a direct computation that for $k=0,1,2,\cdots K$
$$\align
&E_k(d_0,\psi,\na\psi)\\
&=\biggl(k(k+1)(s_0-b_0)-2(\g-1)k(s_0-b_0)+(\g-1)b_0\biggr)(1+\a)\\
&\quad +(O((s_0-b_0)^3)+O_k(\ve).\tag A.4
\endalign$$

It is noted that $s_0-b_0=O(b_0^{-\f{3-\g}{\g-1}})+O(b_0^{-1})$
holds. This, together with (A.4), yields (A.3). \qquad \qquad \qquad \qed\medskip

For the functions $D_{21}^0, D_{22}^0$ in (3.21) and $D_{21}^k, D_{22}^k$ ($k\ge 1$) in (3.33), we have

{\bf Lemma A.5.} {\it Under the assumptions in \rom{Lemma A.4}, one
has on $R=2$ and for $k=0, 1, ..., K$
$$\cases
D_{21}^k=-\hat\rho_0(1+\a)<0,\\
\\
D_{22}^k\leq
-\ds\f{(1+k)(s_0-b_0)}{b_0}\hat\rho(1+\a)<0.\endcases$$}\medskip

{\bf Proof.} It follows from (3.9) that
$$\cases
\ds \p_{\psi}\Bbb H(b,\psi,\na\psi)&=\ds\f{\rho}{c^2(\rho)}\p_{\psi}A_0(b,\psi,\na\psi),\\
\ds \p_{\p_R\psi}\Bbb H(b,\psi,\na\psi)&=\ds\f{\rho}{c^2(\rho)}\p_{\p_R\psi}A_0(b,\psi,\na\psi),\\
\ds \p_{\p_T\psi}\Bbb H(b,\psi,\na\psi)&=\ds\f{\rho}{c^2(\rho)}\p_{\p_T\psi}A_0(b,\psi,\na\psi).
\endcases$$

This, together with Lemma 3.2, yields for $k=0,1,2,\cdots, K$ and small $\ve$
$$\align
D_{21}^k(d_0,\psi,\na\psi)=&\p_{\p_R\psi}\biggl(\Bbb H(\cdot)\hat\psi-\f{1}{b_0
a_1}(\Bbb H(\cdot)-\rho_0)
(T\p_T a_0+a_0)\biggr)(b_0,\hat\psi,\na\hat\psi)+O_k(\ve)\\
=&\biggl(\f{\hat\rho}{c^2(\hat\rho)}\hat\psi\p_{\p_R\psi}A_0-\f{1}{b_0
a_1}\f{\hat\rho_0}{c^2(\hat\rho_0)}\p_{\p_R\psi}A_0
(T\p_T a_0+a_0)\\
&-\f{1}{b_0}(\Bbb H(\cdot)-\rho_0)(T\p_T
a_0+a_0)\biggr)(b_0,\hat\psi,\na\hat\psi)+O_k(\ve)\\
=&-\hat\rho_0(1+\a)<0
\endalign$$
and
$$\align
D_{22}^k(d_0,\psi,\na\psi)=&\biggl(\p_{\psi}+\f{1}{T}k\p_{\p_T\psi}\biggr)
\biggl(\Bbb H(\cdot)\psi-\f{1}{b_0
a_1}(\Bbb H(\cdot)-\rho_0)
(T\p_T a_0+a_0)\biggr)(b_0,\hat\psi)+O_k(\ve)\\
=&\biggl(\hat\psi\f{\hat\rho}{c^2(\hat\rho)}\p_{\psi}A_0+\Bbb H(\cdot)-\f{1}{b_0
a_1}\f{\hat\rho}{c^2(\hat\rho)}\p_{\psi}A_0(T\p_T
a_0+a_0)-\f{1}{b_0}(\Bbb H(\cdot)-\rho_0)(T\p_T a_0+a_0)\\
&-\f{1}{b_0
a_1}(\Bbb H(\cdot)-\rho_0)+\f{k}{T}\bigl[\hat\psi\f{\hat\rho}{c^2(\hat\rho)}\p_{\p_T\psi}A_0-\f{1}{b_0
a_1}\f{\hat\rho}{c^2(\hat\rho)}\p_{\p_T\psi}A_0(T\p_T
a_0+a_0)\\
&-\f{T}{b_0
a_1}(\Bbb H(\cdot)-\rho_0)\bigr]\biggr)(b_0,\hat\psi)+O_k(\ve)\\
\leq&-\f{(1+k)(s_0-b_0)}{b_0}\hat\rho_0(1+\a).
\endalign$$

Thus, we complete the proof on Lemma A.5.\qquad\qquad\qquad\qquad\qquad\qquad\qquad\qquad\qquad
\qquad\qquad \qed\medskip

{\bf Lemma A.6.}  {\it In \rom{(3.48)}, if the variables $(b,
\psi^k, \dot\psi_{n-1})$ are replaced by $(b_0, \hat\psi, 0)$
correspondingly, then we have
$$\align
&\ds \Cal B_{20}^{n}=-\bigl(\f{\g-1}{2A\g}\bigr)^{\f{1}{\g-1}}
b_0^{\f{3-\g}{\g-1}}(s_0-b_0)(1+\a),\\
&\ds \Cal B_{21}^{n}=-\bigl(\f{\g-1}{2A\g}\bigr)^{\f{1}{\g-1}}b_0^{\f{2}{\g-1}}(1+\a),\\
&\ds \Cal B_{22}^{i, n}=\a,\qquad \qquad i=1,2,3.\\
\endalign$$
}

{\bf Proof.} For $(b, \psi^k, \dot\psi_{n-1})=(b_0, \hat\psi, 0)$, we have
$$\cases
&\ds\Cal B_{20}^{n}=-\f{1}{b_0 a_1}(\Bbb H(b_0,\hat\psi,\na\hat\psi)-\rho_0)
+e^{-X}(\psi-\f{1}{b_0 a_1}(\p_X a_0+a_0))\cdot\p_{\p_T\psi}\Bbb H(\cdot)\biggl|_{b=b_0,\psi=\hat\psi},\\
&\ds\Cal B_{21}^{n}=-\f{1}{b_0}(\Bbb H(\cdot)-\rho_0)(\p_X a_0+a_0)
+(\hat\psi-\f{1}{b_0 a_1}(\p_X a_0+a_0))\p_{\p_R\psi}\Bbb H(\cdot)\biggl|_{b=b_0,\psi=\hat\psi}\\
&\ds\Cal B_{22}^{i, n}=0, \qquad\qquad \qquad i=1,2,3.
\endcases$$

This, together with Lemma 2.1 and Lemma 3.2, yields for $(b, \psi^k, \dot\psi_{n-1})=(b_0, \hat\psi, 0)$
$$\cases
\ds \Cal B_{20}^{n}&=\ds-\f{1}{b_0}\hat\rho_0\hat\psi(1+\a)\\
&=-\bigl(\f{\g-1}{2A\g}\bigr)^{\f{1}{\g-1}}b_0^{\f{3-\g}{\g-1}}(s_0-b_0)(1+\a),\\
\ds \Cal B_{21}^{n}&=\ds-\hat\rho_0\hat\psi(1+\a)\\
&=-\bigl(\f{\g-1}{2A\g}\bigr)^{\f{1}{\g-1}}b_0^{\f{2}{\g-1}}(1+\a),\\
\Cal B_{22}^{i, n}&=\a, \qquad \qquad i=1,2,3.\\
\endcases$$

Therefore, we complete the proof of Lemma A.6. \qquad\qquad\qquad\qquad
\qquad\qquad\qquad\qquad\qquad\qquad \qed

\newpage

\vskip 0.4 true cm \centerline{\bf Appendix B. Modified background
solution} \vskip 0.4 true cm

In this Appendix, we look for a modified  background solution
$\Phi_{a}(t,x)$ such that it satisfies the Neumann-type boundary condition (1.7) since
the background solution $\hat\Phi(t,x)$ given in
Remark 2.1 does not satisfy (1.7), which will be crucial in looking for the
multiplier to derive the energy estimate for the problem $(4.4)$
with $(4.9)$ and $(4.11)$-$(4.13)$.

{\bf Lemma B.1.} {\it Assume $\ds|\f{\zeta}{t}-s_0|\leq \f{s_0-b_0}{2}$, then there exists
a smooth function
$f_a=f_{a}(t,r,\o)$ such that
$$f_a\in O_{-1}^{\infty}(\ve)\tag B.1$$
and the function
$\Phi_a(t,r,\o)=(1+f_a(t,r,\o))\hat\Phi(t,x)$
satisfies
$$\Cal B_{\sigma}\Phi_{a}=\p_t\sigma\quad \text{on}\quad r=\sigma(t,\o),\tag B.2$$
where the operator $\Cal B_\si$ is given in \rom{(1.7)}, and the
meaning of the notation $O_{-1}^{\infty}(\ve)$ can be referred to
\rom{(4.2)}.}\medskip

{\bf Proof.} We look for $f_a(t ,r,\o)$ to satisfy
$$\cases
&\ds \Cal B_{\sigma}f_a=\f{1}{\hat\Phi(\f{\sigma}{t})}\biggl(\p_t\sigma-\hat u(\f{\sigma}{t})\biggr)
\quad \text{in}\quad\Omega_+,\\
& f_a=0\quad \text{on}\quad r=\sigma(t,\o),
\endcases\tag B.3
$$
which implies that (B.2) holds.

Obviously, under the assumptions in Theorem 1.1, $r=\sigma(t,\o)$ is
not a characteristic surface of the first order linear equation in (B.3). This means
that (B.3) has a unique smooth solution $f_a(t,r,\o)$ in
$\O_+$. We now analyze the properties of the solution
$f_a(t,r,\o)$.

In terms of Lemma 4.5 (iv), $\Cal B_\sigma$ has such a form
$\Cal B_\sigma=\ds\sum_{i=1}^{n}(\f{x_i}{r}-\p_i\sigma)\cdot\p_i$.
We assume that $$L_{t,r_0,\o_0}^{r}:\quad x=x(r; t,r_0,\o_0)$$ is an integral curve
of $\Cal
B_{\sigma}$, which starts from the point $(t,r_0,\o_0)$ with
$r_0=\si(t,\o_0)$. That is,  $x=x(r;t,r_0,\o_0)$ satisfies
$$\cases
&\ds\f{dx_i(r;t,r_0,\o_0)}{dr}=\biggl(\f{x_i}{r}-\p_i\si\biggr)\biggl|_{L_{t,r_0,\o_0}^r},\qquad
i=1,\cdots,n,\\
&x(r_0;t,r_0,\o_0)=r_0\o_0.
\endcases\tag B.4$$

Under the assumptions in Theorem 1.1, the ODE system (B.4) has a
unique smooth solution in $\O_+$, which is written as $$x=x(r;t,r_0,\o_0).$$

It follows from (B.3) and (B.4) that
$$\cases
&\ds\f{df_a}{dr}=
\f{1}{\hat\Phi(\f{\sigma}{t})}(\p_t\sigma-\hat
u(\f{\sigma}{t}))\quad \text{on}\quad L_{t,r_0,\o_0}^r\cap\O_+,\\
&f_a(t,r_0,\o_0)=0.
\endcases$$

>From this,  $f_a(t,r,\o)$ admits an explicit expression as follows
$$f_a(t,r,\o)=\int_{r_0}^{r}\f{1}{\hat\Phi(\f{\sigma}{t})}(\p_t\sigma-\hat
u(\f{\sigma}{t}))\biggl|_{L_{t,r_0,\o_0}^\tau\cap\O_+}d\tau.\tag
B.5$$

By the assumptions of Theorem 1.1 and $\hat u(b_0)=b_0$, one has
$$E(t,\o)\triangleq\f{1}{\hat\Phi(\f{\sigma}{t})}(\p_t\sigma-\hat
u(\f{\sigma}{t}))=O_{-2}^{\infty}(\ve).$$

This, together with $|\zeta-\si|\leq 2(s_0-b_0)t$ and (B.5), yields
$$|f_a|\leq \f{C\ve}{t^2}|\zeta-\si|\leq \f{C\ve}{t}.\tag B.6$$

For any vector field $S\in \Cal S$ (here $\Cal S$ denotes the
modified Klainerman's vector fields in (5.42)), in terms of  Remark 5.3, we
have
$$\cases
&\Cal B_\si S f_a=S E(t,\o)+[\Cal B_\si, S]f_a\quad \text{in}\quad
\O_+,\\
\\
&S f_a=0\quad \text{on}\quad r=\si(t,\o).
\endcases\tag B.7$$

By (B.6) and some related computations in the proof procedure of Lemma 5.4,  (B.7) can be reduced into
$$\cases
&\Cal B_\si S f_a=O_{-2}^{\infty}(\ve)+O_{-1}^{\infty}(\ve)\t Z f_a\quad \text{in}\quad \O_+,\\
\\
& S f_a=0\quad \text{on}\quad r=\si(t,\o),
\endcases\tag B.8$$
where $\t Z$ denotes a linear combination of some components in
$\Cal S$.

It follows
from (B.8) that
$$|S f_a|\leq \f{C\ve}{t}.\tag B.9$$

>From this,  by the equation in
(B.3) and the definition of $S$, one has
$$|\p_r f_a|\leq \f{C\ve}{t}|S f_a|+\ve|\p_r
f_a|+\f{C\ve}{t^2}$$ and further
$$|\p_r f_a|\leq \f{C\ve}{t^2}.\tag B.10$$

Consequently, it follow Lemma 4.5 (vi) and (B.9)-(B.10) that
$$|\na f_a|\leq \f{C\ve}{t^2}.$$

Analogously, by induction method, we can arrive at
$$|\na^l f_a|\leq \f{C_l\ve}{t^{1+l}},\quad  l\in\Bbb N.\tag B.11$$

Combining (B.11) with (B.6) shows (B.1) and the proof of Lemma B.1
is completed. \qquad \qquad \qquad \qquad  \qed

\medskip

{\bf Acknowledgments.} The third author Yin Huicheng would like to thank Professor Wu Sijue
for her interests and some fruitful discussions on the topics in this paper when
she visited the Institute of Mathematical
Sciences at Nanjing University in December of 2011.


\vskip 0.6 true cm

Li Jun

\vskip 0.2 true cm

Department of Mathematics and IMS, Nanjing University,
Nanjing 210093, P.R.~China

Email address: lijun\@nju.edu.cn

\vskip 0.6 true cm

Witt Ingo

\vskip 0.2 true cm

Mathematical Institute, University of G\"{o}ttingen,
Bunsenstr.~3-5, D-37073 G\"{o}ttingen, Germany

Email address: iwitt\@uni-math.gwdg.de

\vskip 0.6 true cm

Yin Huicheng
\vskip 0.2 true cm

Department of Mathematics and IMS, Nanjing University,
Nanjing 210093, P.R.~China

Email address: huicheng\@nju.edu.cn

\bye

\Refs \refstyle{C}





\ref\key 1\by M. M. Ad$^{\prime\prime}$yutov, Yu. A. Klokov,
A. P. Mikhailov, V. V. Stepanova, A. A. Shamrai\paper The existence and
properties of the solutions of a self-similar problem on a flat
piston\jour Mat. Model. 5, no. 7, 71-85\yr 1993\endref




\ref\key 2\by S. Alinhac\paper Blowup of small data solutions for a
quasilinear wave equation in two space dimensions\jour  Ann. of Math. (2)
149, no. 1, 97-127\yr 1999\endref


\ref\key 3\by S. Alinhac\paper Blowup of small data solutions for a
class of quasilinear wave equations in two space dimensions. II
\jour Acta Math. 182, no. 1, 1-23\yr 1999\endref


\ref\key 4\by Chen Shuxing\paper  A singular multi-dimensional
piston problem in compressible flow\jour J.Differential Equations 189,
no. 1, 292-317\yr 2003\endref

\ref\key 5\by Chen Shuxing\paper Existence of stationary supersonic flows past a pointed body
\jour Arch. Ration. Mech. Anal. 156, no. 2, 141-181\yr 2001\endref


\ref\key 6\by D. Christodoulou\paper Global solutions of nonlinear
hyperbolic equations for small initial data\jour  Comm. Pure Appl. Math. 39,
no.2, 267-282\yr 1986
\endref


\ref\key 7\by R. Courant, K. O. Friedrichs\paper Supersonic flow and
shock waves \publ Interscience Publishers Inc., New York, 1948
\endref

\ref\key 8\by C. M. Dafermos\paper Hyperbolic conservation laws in
continuum physics\publ Springer, Berlin, Heidelberg, New York, 2000
\endref


\ref\key 9\by D. Gilbarg, N. S. Tudinger\paper Elliptic partial
differential equations of second order\publ Second edition. Grundlehren
der Mathematischen Wissenschaften, 224, Springer, Berlin-New York\yr
1983
\endref

\ref\key 10\by P. Godin \paper Long time existence of a class of
perturbations of planar shock fronts for second order hyperbolic
conservation laws\jour Duke Math. J. 60 (2), 425-463\yr 1990
\endref

\ref\key 11\by P. Godin \paper Global shock waves in some domains for
the isentropic irrotational potential flow equations\jour
Comm. P. D. E., Vol. 22, no. 11-12, 1929-1997 \yr 1997
\endref




\ref\key 12\by G. H. Hardy, J. E. Littlewood, G. Polya \book Inequality
\publ Cambridge University Press, London, New York \yr 1964
\endref


\ref\key 13\by L. H$\ddot o$rmander\book Lectures on nonlinear
hyperbolic differential equations\publ Springer-Verlag\yr 1997
\endref

\ref\key 14\by M. Ikawa\paper A mixed problem for hyperbolic
equations of second order with a first order derivative boundary
condition\jour  Publ. Res. Inst. Math. Sci., Kyoto Univ. 5, 119-147\yr
1969
\endref


\ref\key 15\by F. John\book Nonlinear wave equations, formation of
singularities, University Lecture Series 2\publ American
Mathematical Society, Providence, RI\yr 1990
\endref

\ref\key 16\by S. Klainerman\paper The null condition and global
existence to nonlinear wave equations, Nonlinear systems of partial
differential equations in applied mathematics, Part 1 (Santa Fe,
N.M., 1984), 293-326\publ  Lectures in Appl. Math., 23, Amer. Math. Soc.,
Providence, RI, 1986
\endref

\ref\key 17\by S. Klainerman, T.Sideris \paper On almost global
existence for nonrelativistic wave equations in 3D \jour Comm. Pure
Appl. Math., Vol. 49, 307-321 \yr 1996
\endref

\ref\key 18\by E. Lefrancois, J. P. Boufflet\paper An introduction to
fluid-structure interaction: application to the piston problem\jour
SIAM Rev. 52, no. 4, 747-767\yr 2010
\endref

\ref\key 19\by Li Jun, Witt Ingo, Yin Huicheng\paper On the global
existence and stability of a three-dimensional supersonic conic
shock wave\jour Comm. Math. Phys. (to appear), arXiv: 1110. 0626\yr 2011\endref

\ref\key 20\by T. T. Li\paper Global classical solutions for
quasilinear hyperbolic systems\publ Wiley, Masson, New York, Paris,
1994\endref


\ref\key 21\by H. Lindblad, I. Rodnianski\paper Global stability of
Minkowski space-time in harmonic gauge\jour  Annals of
Mathematics, 171, 1401-1477\yr 2010
\endref


\ref\key 22\by Liu Tai-Ping\paper The free piston problem for gas
dynamics\jour  J. Differential Equations 30, no. 2, 175-191\yr 1978
\endref




\ref\key 23\by A. Majda\paper Compressible fluid flow and systems of
conservation laws in several space variables\publ Applied Mathematical
Sciences 53(1984), Springer, New York, Berlin, Heidelberg,
Tokyo\endref

\ref \key 24\by A. Majda, E. Thomann \paper Multi-dimensional shock
fronts for second order wave equations \jour Comm. P. D. E., Vol. 12,
777-828 \yr 1987
\endref



\ref \key 25\by G. M\' etivier\paper Interaction de deux chocs pour un syst$\acute e$me
de deux lois de conservation,
en dimension deux d'espace\jour Trans. Amer. Math. Soc. 296, no. 2, 431-479\yr
1986
\endref


\ref \key 26\by W. Paul\paper The periodic oscillation of an
adiabatic piston in two or three dimensions\jour  Comm. Math. Phys. 275,
no. 2, 553-580\yr 2007
\endref


\ref \key 27\by J. Rauch\paper BV estimates fail for most quasilinear
hyperbolic systems in dimension greater than one\jour
Comm. Math. Phys., 106, 481-484\yr 1986
\endref

\ref\key 28\by Y. Shibata\paper On the Neumann problem for some
linear hyperbolic systems of second order\jour  Tsukuba J. Math. 12,
149-209\yr 1998
\endref


\ref \key 29\by T. Sideris \paper Formation of singularities in three
dimensional compressible fluids \jour Comm. Math. Phys., 101, 475-487
\yr 1985
\endref


\ref \key 30\by M. Slemrod\paper Resolution of the spherical piston
problem for compressible isentropic gas dynamics via a self-similar
viscous limit\jour  Proc. Roy. Soc. Edinburgh Sect. A 126, no. 6,
1309-1340\yr 1996
\endref

\ref \key 31\by J. Smoller\paper Shock waves and reaction-diffusion
equations, Springer, Berlin, New York, 1983
\endref

\ref \key 32\by G. B. Whitham\paper Linear and nonlinear waves\publ Wiley,
New York, London, Sydney, Toronto, 1974\endref

\ref \key 33\by Wu Sijue\paper Global wellposedness of the 3-D full water wave problem
\jour Invent. Math. 184, no. 1, 125-220\yr 2011\endref

\ref \key 34\by Wu Sijue\paper Almost global wellposedness of the 2-D full water wave problem
\jour Invent. Math. 177, no. 1, 45-135\yr 2009\endref




\ref\key 35\by Xin Zhouping\book Some current topics in nonlinear
conservation laws\publ  Some current topics on nonlinear conservation
laws, xiii-xxxi, AMS/IP Stud. Adv. Math., 15, Amer. Math. Soc.,
Providence, RI, 2000\endref


\ref\key 36\by Xin Zhouping, Yin Huicheng\paper Global
multidimensional shock wave for the steady supersonic flow past a
three-dimensional curved cone\jour Anal. Appl. (Singap.) 4, no. 2,
101-132\yr 2006\endref



\ref\key 37\by Yin Huicheng\paper Formation and construction of a
shock wave for 3-D compressible Euler equations with the spherical
initial data\jour Nagoya Math. J. 175, 125-164 \yr 2004
\endref


\ref\key 38\by Zheng Yuxi\paper Systems of conservation laws.
Two-dimensional Riemann problems\publ Progress in Nonlinear Differential
Equations and their Applications, 38, Birkh$\ddot a$user Boston,
Inc., Boston, MA, 2001\endref
\enddocument